\documentclass[11pt]{article}
\usepackage[utf8]{inputenc}
\usepackage[english]{babel}
\usepackage{amsmath}
\usepackage{amsthm}
\usepackage{amsfonts}
\usepackage{amssymb}
\usepackage{appendix}
\usepackage{graphicx}
\usepackage{color}
\usepackage{rotating}
\usepackage{authblk}
\usepackage{url}
\usepackage{tcolorbox}
\usepackage{fullpage}
\usepackage[linesnumbered,ruled, noend]{algorithm2e}
\usepackage{tabularx}
\usepackage[caption=false]{subfig}
\usepackage{multirow, float, makecell}
\usepackage{mathrsfs}
\usepackage[labelfont=bf,labelsep=period,justification=raggedright]{caption}
\graphicspath{{images/}}

\usepackage{natbib}
\setcitestyle{authoryear,open={((},close={))}}
 \bibpunct[, ]{(}{)}{,}{a}{}{,}%

\newtheorem{proposition}{Proposition}	

\allowdisplaybreaks

\begin{document}
\title{A K-Nearest Neighbor Heuristic for Real-Time DC Optimal Transmission Switching}
\author{Emma~S.~Johnson$^{1,2}$, Shabbir~Ahmed$^{1*}$, Santanu~S.~Dey$^{1}$, Jean-Paul~Watson$^{3}$\\
	\small $^{1}$School of Industrial and Systems Engineering, Georgia Institute of Technology, Atlanta, GA, USA. \\
	\small ejohnson335@gatech.edu, santanu.dey@isye.gatech.edu \\
	\small $^{2}$Sandia National Laboratories, Albuquerque, NM, USA \\
	\small $^{3}$Lawrence Livermore National Laboratory, Livermore, CA, USA. watson61@llnl.gov \\
	\small $^*$This author was deceased at the time of submission. \\
}

\maketitle

\begin{abstract}
While transmission switching is known to reduce power generation costs, the difficulty of solving even DC optimal transmission switching (DCOTS) has prevented optimal transmission switching from becoming commonplace in real-time power systems operation.
In this paper, we present a k-nearest neighbors (KNN) heuristic for DCOTS which relies on the insight that, for routine operations on a fixed network, the DCOTS solutions for similar load profiles and generation cost profiles will likely turn off similar sets of lines.
We take a data-driven approach and assume that we have DCOTS solutions for many historical instances, which is realistic given that the problem is solved every 5 minutes in practice.
Given a new instance, we find a set of ``close" instances from the past and return the best of their solutions for the new instance.
We present a case study on 7 test networks ranging in size from 118 to 3,375 buses. 
We compare the proposed heuristic to DCOTS heuristics from the literature, commercial solver heuristics, and a simple greedy local search algorithm.
In most cases, we find better quality solutions in less computational time.
In addition, the computational time is within the limits imposed by real-time operations, even on larger networks.
Last, we present an empirical study of our training data to understand why the heuristic works well.
\end{abstract}


\section{Introduction}

\label{sec:intro}
Transmission switching is an inexpensive way to reduce generation costs in a congested power system. 
Ideally, we want to always dispatch the cheapest generators first, which we call a merit-order dispatch.
However, network constraints such as line flow limits and Ohm's law can make this dispatch infeasible.
Due to Braess' Paradox, by not allowing flow through certain lines, we can reduce network congestion and approach the desired merit-order dispatch (\cite{BlumsackLI2007}).
\cite{FisherOF2008} show that an optimal network topology for one load profile is not necessarily optimal for another, presenting a need for real-time switching as load fluctuates.

Transmission switching in real-time poses a challenge since the DC optimal power flow (DCOPF) problem is solved every 5 minutes by most independent system operators (ISOs) to determine generator dispatch for the next time interval. 
This means that the DCOPF problem itself should be solved in significantly less than 5 minutes to allow time in the remainder of the 5-minute interval for feasibility checks and other post-processing.
Since DC optimal transmission switching requires co-optimizing the generator dispatch and the network topology, to be useful in practice, DC optimal transmission switching (DCOTS) needs to be solved in less than 5 minutes.
This is computationally challenging since DCOTS is NP-hard, and also cannot be approximated within a constant factor (\cite{LehmannGH2014, KocukJDLLS2016}).

Many mixed-integer programming (MIP) formulations and heuristics for DCOTS have been proposed, but none scale well enough that they could be implemented in real-time.
There are a variety of formulations in the literature, but, because they scale poorly with network size, they can at best be used as heuristics by limiting the set of candidate lines for switching.
\cite{KocukJDLLS2016} proposes a cycle-based formulation for DCOTS and use it to derive valid inequalities leveraged in a cutting plane solution approach.
The cuts improve computational time, but the approach does not scale such that it could be used in real-time.
\cite{RuizGRCPF2017} proposes a smaller formulation for security-constrained DCOTS based on the Power Transfer Distribution Factor (PTDF) power flow equations.
By security-constrained, we mean the solution is also feasible under a set of predefined contingency scenarios.
This formulation can be used as a heuristic since it is computationally efficient when the number of contingencies is small and the set of candidate lines for switching is limited.
\cite{Jabernejad2018} proposes an approximate model for DCOTS which is guaranteed to yield feasible solutions with the same generation costs but lower numbers of switched lines.
However, using this approximation for large networks still entails solving a large mixed integer program, which does not yet scale well.

The difficulty of MIP formulations of DCOTS has so far prohibited exactly optimizing the topology on large networks within the time limit imposed by real-time.
However, numerous heuristics have been proposed.
\cite{RuizFRC2012} use four different metrics based on sensitivity analysis of DCOPF in order to estimate the cost benefit of switching each line of the network.
They use these estimates to prioritize the lines.
They then iterate through this priority list, switching off lines that result in cost savings.
Similarly, \cite{FullerRC2012} uses one of the same sensitivity criteria to choose a subset of high priority lines. 
They then run a greedy algorithm over just this set of lines, switching off the line resulting in the most cost savings until there are no more cost-saving lines or until they reach a maximum cardinality of switched lines.
The sensitivity-based heuristics are effective in terms of solution quality.
However, the number of linear programs they solve scales with the product of the number of lines in the system and the budget of how many lines to switch.
This makes them impractical for large networks and for switching large numbers of lines.

Others have proposed using sensitivity information or historical data to limit the set of switchable lines so as to reduce the number of binaries in the DCOTS problem, e.g. \cite{LiuWO2012}. 
In an effort to understand the economic impacts of switching via a case study on a congested network, \cite{HedmanOFO2011} suggests two heuristics for DCOTS.
One involves a variation of the greedy algorithm used in \cite{FullerRC2012}, but uses a partitioning technique to explore multiple disjoint sets of feasible solutions in parallel.
The second uses data from past instances to select a limited set of switchable lines which were cost-beneficial to switch off in the past.
They then solve the DCOTS problem allowing only these lines to be opened.

The idea of using data from past solves has recently begun to take hold in the form of learning-based heuristics.
Since the dispatch problem is solved so frequently, and under normal conditions demands and generation costs for a given time period do not vary beyond around 10\% and 5\% respectively (\cite{AlinsonQA2019}), a likely avenue for scalable DCOTS heuristics is to harness off-line computational power and develop a heuristic based on the solutions to historical instances of the problem.

There has been recent interest in applying machine learning methods to power systems problems.
In \cite{AlinsonQA2019}, the authors present three learning-based methods which use historical data to solve the security constrained unit commitment problem.
\cite{YangO2019} applies k-nearest neighbors, an artificial neural network, and decision tree regression to learn sets of high-priority lines to consider for switching.
They then use a greedy algorithm based on \cite{FullerRC2012}, which uses this line prioritization to generate a topology.
In addition, they use machine learning methods to train an algorithm selection oracle which, given an instance, chooses which among these algorithms to run.

We propose a k-nearest neighbors approach different than that of \cite{YangO2019}.
We apply the method for learning initial feasible solutions from \cite{AlinsonQA2019} to DCOTS.
Rather than training an oracle to map parameter vectors to sets of high priority lines, as in \cite{YangO2019}, we instead use k-nearest neighbors to learn topologies directly from the instance data.
More specifically, we assume that we have a large collection of solved instances of DCOTS.
Given a new instance, our heuristic selects the nearest $k$ instances in parameter space out of this solved collection.
We then test the quality of each of the optimal toplogies of the $k$ nearest instances and return the lowest-cost topology as our switching solution.
We show through case studies on seven different transmission networks ranging from 118 to 3,375 buses that we achieve, on average, better solutions than other heuristics in less time, making our heuristic potentially practical for real-time operations. 
Even for larger networks, the number of solved training instances needed is moderate.
A collection of 270 training instances yields solutions that are usually within 2\% of the best known solution and, in most cases, are within 1\% of the best known solution.
Another advantage of the method is that it scales well for larger networks, never taking more than 40 seconds on any of our test instances.
This is largely because we never iterate through lists of lines.
Instead, the effect of the size of the network is on training time, which is done offline, and on the computational time for the linear program (LP) solves which check the cost of the near topologies.
However, this effect is minimal since there are only $k$ such solves.

In the remainder of the paper, we give an overview of our heuristic approach in Section \ref{sec:approach}, details on the networks used in our case study in Section \ref{sec:test-cases}, computational results benchmarking our heuristic against others from the literature in Section \ref{sec:results}, an analysis to give some intuition behind the algorithm in Section \ref{sec:why}, and we conclude in Section \ref{sec:conclusions}.

\section{Heuristic Approaches}\label{sec:approach}
The DCOTS problem is given formally in Appendix \ref{sec:model}.
In this section, we will first present our proposed heuristic in Section \ref{sec:learning-heuristic}.
We then present five heuristics that will be used as benchmarks.

\subsection{KNN Heuristic}\label{sec:learning-heuristic}
Our proposed heuristic relies on the assumption that the DCOTS problem has likely been solved on the same network many times.
We assume the variable data are the load profiles and the generation costs.
All other parameters are known, and are constant on a given network.
We can therefore characterize an instance of DCOTS as the vector of generation costs appended to the vector of demands.
If $q$ is such a vector, we will denote an instance of DCOTS as $I(q)$.
Suppose we have a collection of solved DCOTS instances $\mathcal Q$, each with an $\epsilon$-optimal transmission switching solution.
Given a new instance $I(q)$, we propose a heuristic based on k-nearest neighbors to find a transmission switching solution:
Among the set of solved instances, we find the nearest $k$ instances as measured by their parameter vectors' distances from $q$ using a $p$-norm for some $p$.
We then test the transmission switching solutions from each of these $k$ closest instances.
We return the transmission switching solution which has the best objective value.
The algorithm is presented formally in Algorithm \ref{learning-alg} in Appendix \ref{sec:algorithms-formally}.

Note that the solution returned by Algorithm \ref{learning-alg} will respect the cardinality constraint since all the solutions to the training instances do.
However, the solution is not guaranteed to be feasible, which is why we have included slacks on the nodal balance constraint in the DCOTS formulation given in Appendix \ref{sec:model}.
Since we are already using the DC approximation for power flow, the solution would already have to be corrected for feasibility in practice, so we allow these small violations.

The scalability of this algorithm relies mainly on the value of $k$ and the number of training instances in $\mathcal Q$.
Thus, though the size of the network can dramatically increase the training time, the only effect it has on the algorithm's computational time is to increase the time spent in the $k$ DCOPF solves in the loop beginning at line \ref{lp-loop} of Algorithm \ref{learning-alg}.
The time taken in the loop beginning at line \ref{dist-loop} depends on $|\mathcal Q|$.
In our experiments, we found that $|\mathcal Q| = 270$ was sufficient.
We tested Algorithm \ref{learning-alg} with both Euclidean and $\ell_\infty$-norms.

\subsection{Greedy Local Search}\label{sec:local-search}

We compare the above heuristic to the four heuristics from \cite{RuizFRC2012}, to a greedy local search algorithm, presented as the line enumeration algorithm in \cite{YangO2019}, and to Gurobi's primal heuristics.
For completeness, we describe the greedy algorithm here.
We first calculate the cost, fixing all lines closed.
Then, for each line, we fix only that line open and calculate the cost. 
If none of the lines improve the cost when opened, we are done.
Else, we fix open the line that improves the cost the most.
We repeat this process on the remaining set of lines that could be opened, terminating either when we see no improvement from any of the lines or when we have opened as many lines as the cardinality constraint will allow.
The algorithm is given formally in Algorithm \ref{local-search-alg} in Appendix \ref{sec:algorithms-formally}.

Note that this algorithm scales poorly since it requires solving nearly $K\cdot|\mathcal L|$ linear programs.
Also note that, in a very congested system, where the solution with all lines closed is not feasible, this algorithm also does not guarantee a feasible solution, and the heuristics from \cite{RuizFRC2012} do not either.

\subsection{Sensitivity-Based Heuristics}\label{sec:ruiz-heuristics}

We also compare to the four sensitivity-criteria-based heuristics from \cite{RuizFRC2012}.
For brevity, we do not describe the algorithms in detail here, but note that each of these heuristics uses sensitivity information from the DCOPF problem in order to calculate criteria to indicate lines which are likely to be cost-beneficial to open.
Each heuristic uses a different such criterion to order the lines for consideration in a greedy algorithm.
That is, we follow the procedure from Algorithm \ref{local-search-alg}, but the set $S$ is ordered based on the criterion.
The four criteria are the Line Profits criterion, the Price Difference criterion, the Total Cost criterion, and the PTDF-Weighted Cost criterion.
Note that, even without a cardinality restriction, the four different criteria can find different solutions since each switching decision in the greedy algorithm is conditioned on the ones before it, so the ordering of $S$ does change the heuristic solution.

\subsection{Gurobi Heuristics}\label{sec:gurobi-heuristics}

Last, we compare to Gurobi's performance when we run it with the Heuristics parameter set to 1, meaning that it spends all its time on primal heuristics.
We also warm start these runs with the solution where all lines are closed, which is an obvious feasible solution in all but very congested cases.

\section{Test Cases}\label{sec:test-cases}

We test the heuristic on the 118-bus test system as modified in \cite{blumsack}, on the 300-bus test system as modified in \cite{KocukJDLLS2016}, on the 1354 and 2869 bus PEGASE systems (\cite{pegase}), and on modified versions of the 1951 RTE, 2869 PEGASE, and 3375 Polish system from \cite{pglib}. 
For these three instances we used the heavily loaded versions of the systems, which we have labeled with the postfix ``api," as in the \cite{pglib} library.
All but the 118-bus and 300-bus sytems were downloaded from the IEEE PES Power Grid Library (\cite{pglib}): The heavily loaded versions are from v20.07 and the 1354 and 2869 PEGASE systems are from v19.05.
The instances have varying levels of congestion, but we chose only instances which have a cost benefit from transmission switching. 
For example, the original version of the 1951 RTE case was not selected because, with the nominal demands and generation costs, there is no benefit to opening lines. 
Details of these instances are shown in Table \ref{instance-info}.
\begin{table}
	\centering
	\caption{Test instance sizes\label{instance-info}}
	\begin{tabular}{l |  r  r  r  r }
			\hline
			Test Case & Number of Buses & Number of Generators & Number of Lines \\
			\hline
			118blumsack  & 118 & 19 & 186  \\
			300kocuk & 300 & 61 & 411 \\
			1354pegase & 1,354 & 260 & 1,991 \\
			1951rte\_api & 1,951 & 391 & 2,596 \\
			2869pegase & 2,869 & 510 & 4,582 \\
			2869pegase\_api & 2,869 & 510 & 4,582 \\
			3375wp\_api & 3,374 & 596 & 4,161 \\
			\hline
	\end{tabular}
	{}
\end{table}

For each of these systems, we generate 300 instances following the methodology from \cite{AlinsonQA2019}.
For completeness, we describe the process here:
Suppose that $d^0$ is the original vector of demands and $c^0$ is the original vector of generation costs.
For $i \in \{1, 2, \dots, 300\}$, for each bus $b \in \mathcal B$, draw $\beta_b^i$ from a uniform distribution on the interval $[0.9, 1.1]$.
Then let $d_b^i = \beta_b^id_b^0$.
This process results in 300 demand profiles $d^1, d^2, \dots, d^{300}$.
We generate the generation costs in exactly the same way except that we only allow 5\% variation around the nominal generation cost, so we draw from a uniform distribution on $[0.95, 1.05]$.
We solve all 300 of the instances to 1\%-optimality or to a time limit of 0.5 hour, whichever comes first.
We randomly select 30 instances of these 300 as test instances and leave the other 270 for training.
In Appendix \ref{sec:test-instance-congestion}, using these test and training sets, we give further detail regarding the congestion of these test instances and their potential cost benefit from transmission switching.

\section{Computational Results}\label{sec:results}

We report results running the KNN heuristic with $k=10$ and using both the $\ell_2$-norm and the $\ell_\infty$-norm to measure the distance from the training instances.
We give our heuristic and the Gurobi heuristics a time limit of 5 minutes.
We give the four heuristics from \cite{RuizFRC2012} a time limit of 10 minutes to forgive inefficiencies in our implementation.
For all the heuristics, we test for cardinality limits of 5 and 10, and also with no cardinality constraint.
For the 118blumsack case, we also compare to the greedy local search algorithm given in Section \ref{sec:local-search}, but we find it to be the worst-performing heuristic in this case and intractable for larger test instances.

Our model and all our heuristics are implemented in Pyomo (\cite{pyomo-book, pyomo-paper}), relying on Gurobi version 9.0.2 as the solver (\cite{gurobi}).
We solve on a server with 96 Intel Xeon 2.30GHz processors and 529GB RAM.
For all our experiments, we constrain Gurobi to 8 threads.

\begin{figure*}
	\begin{tabular}{c@{\hskip 0in}c}
		\includegraphics[width=0.5\linewidth]{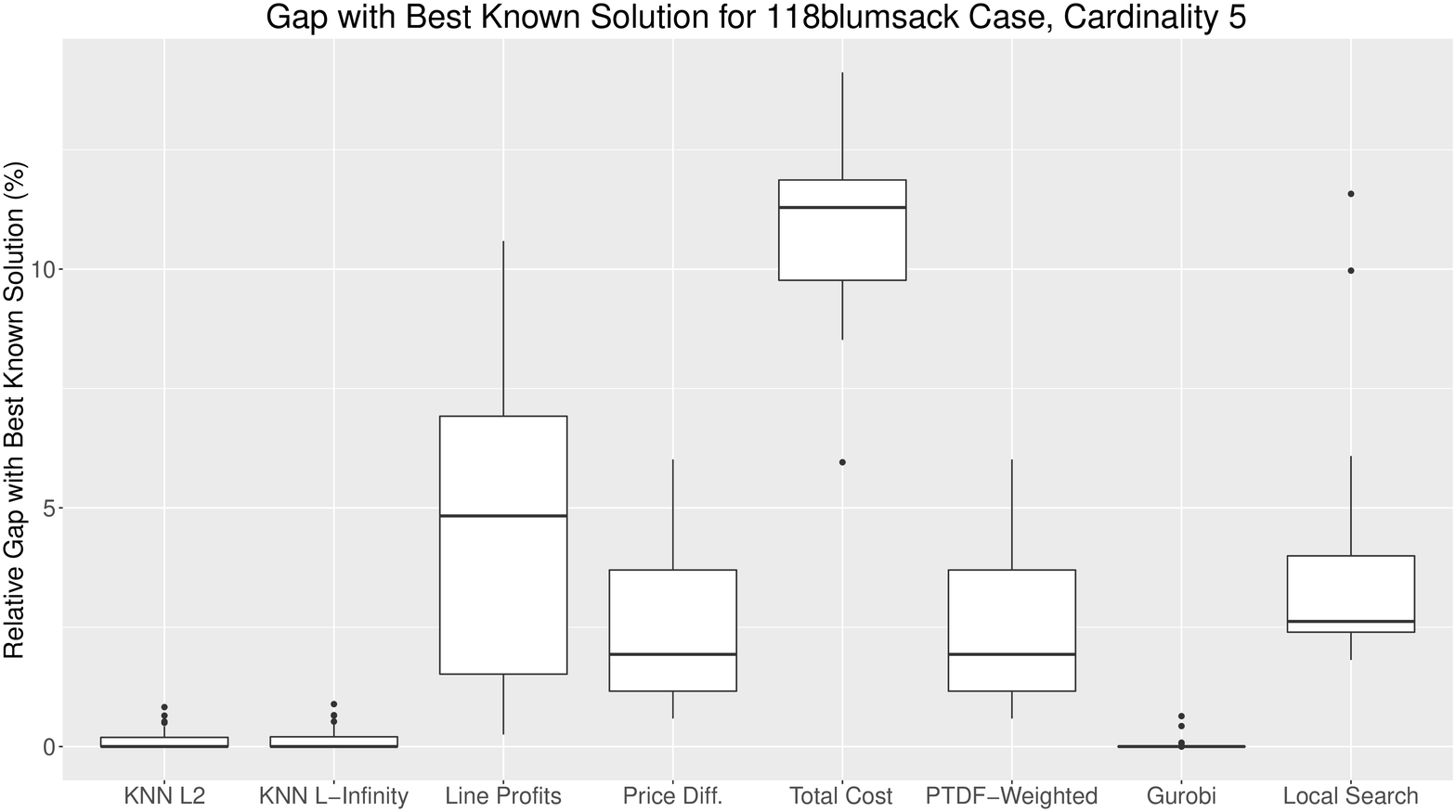} &
		\includegraphics[width=0.5\linewidth]{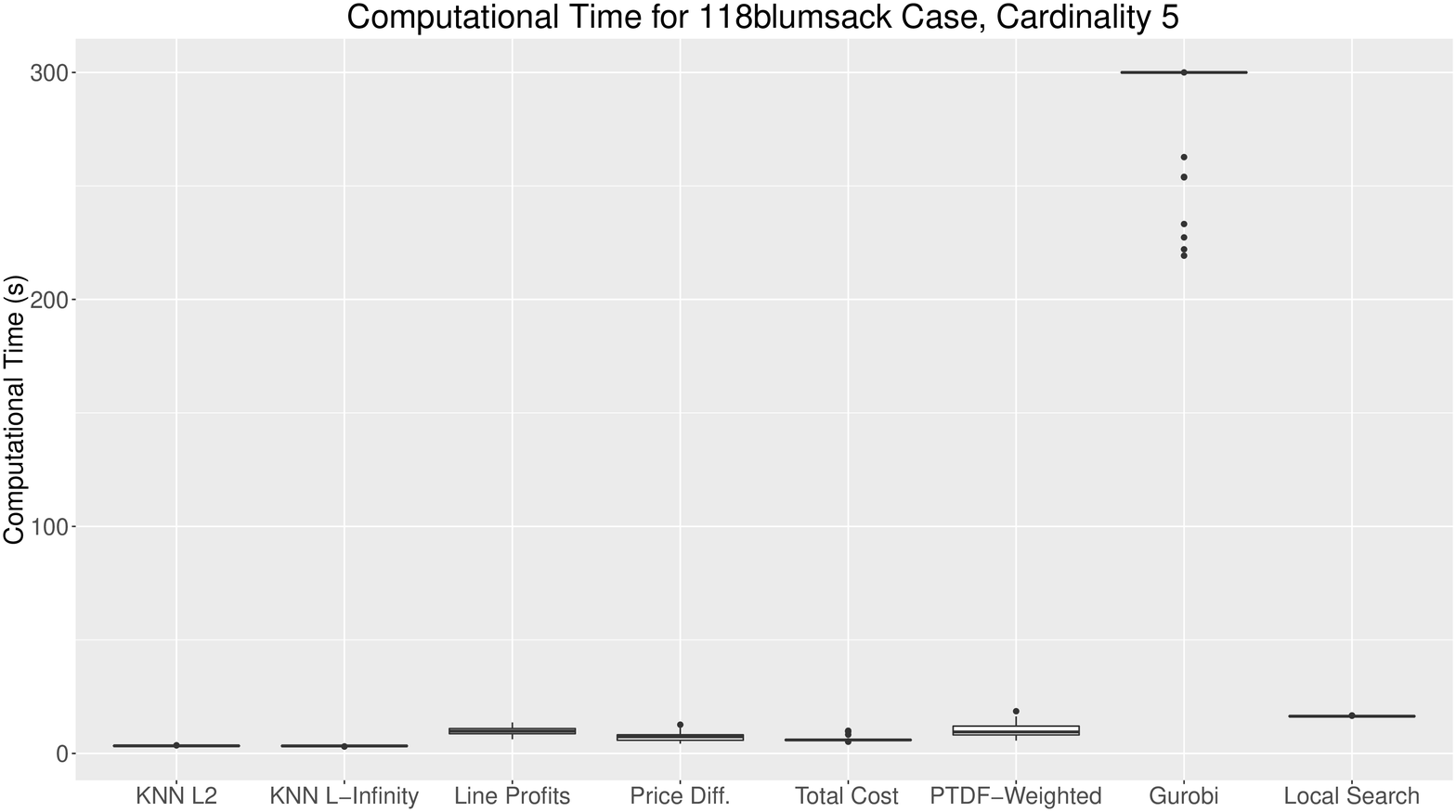} \\
		\includegraphics[width=0.5\linewidth]{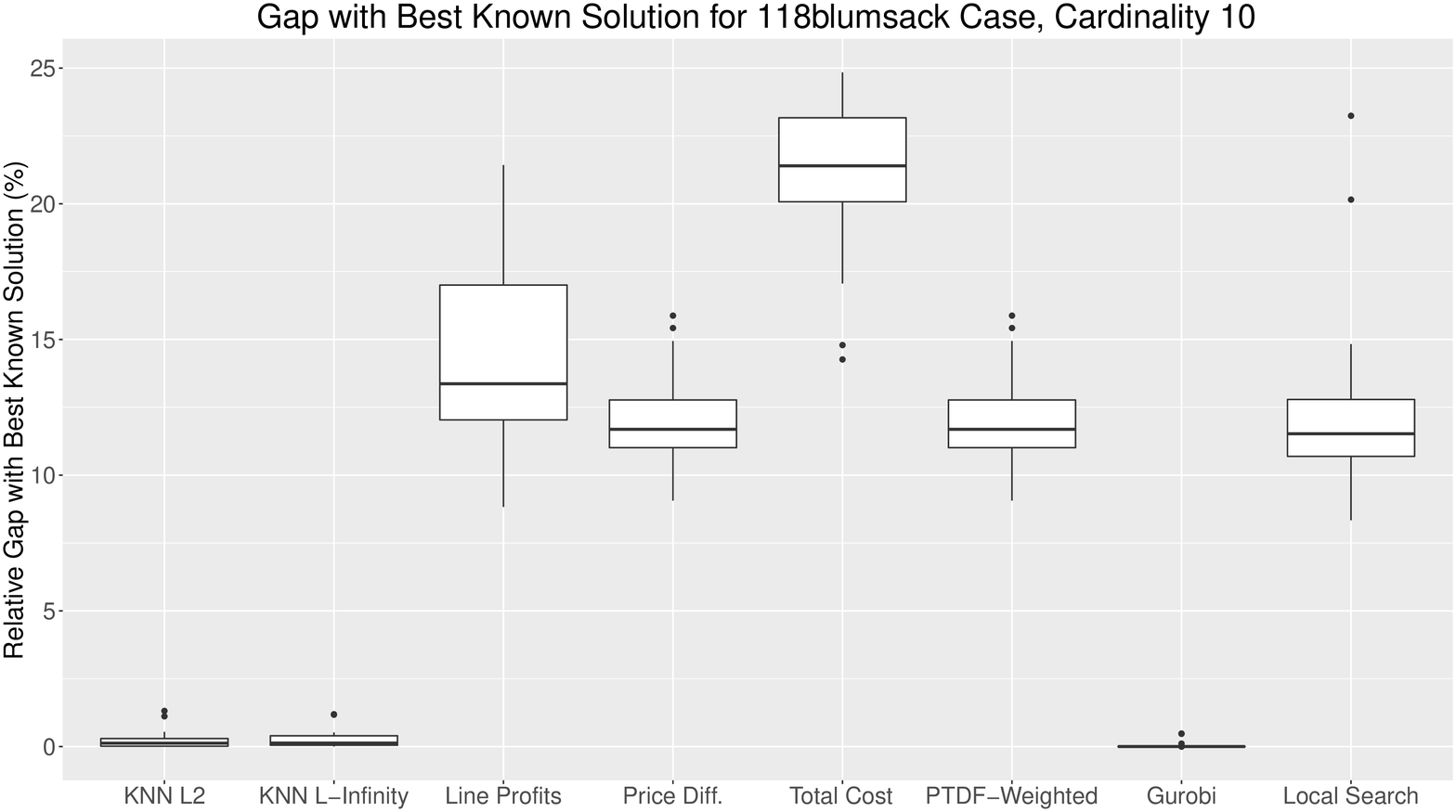} &
		\includegraphics[width=0.5\linewidth]{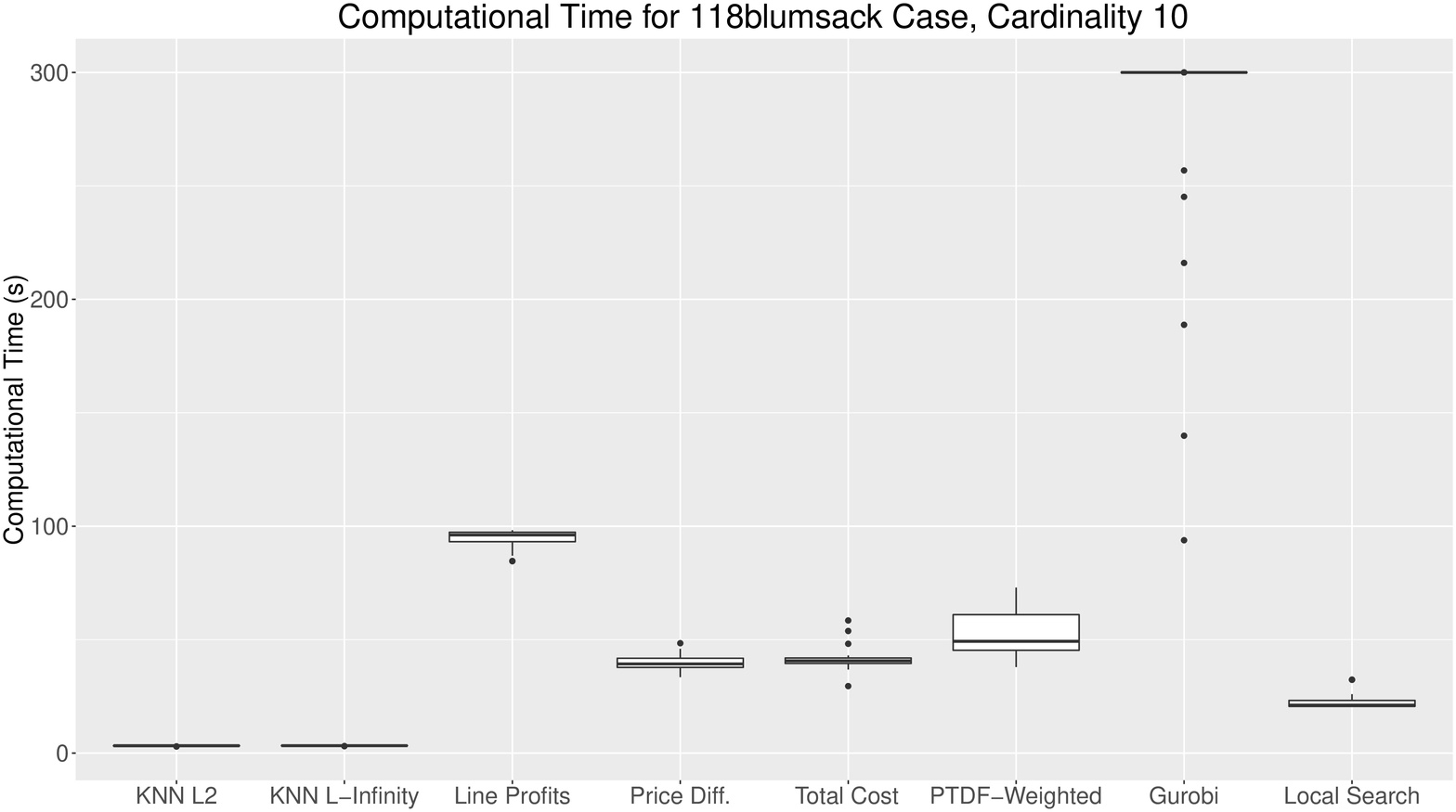} \\
		\includegraphics[width=0.5\linewidth]{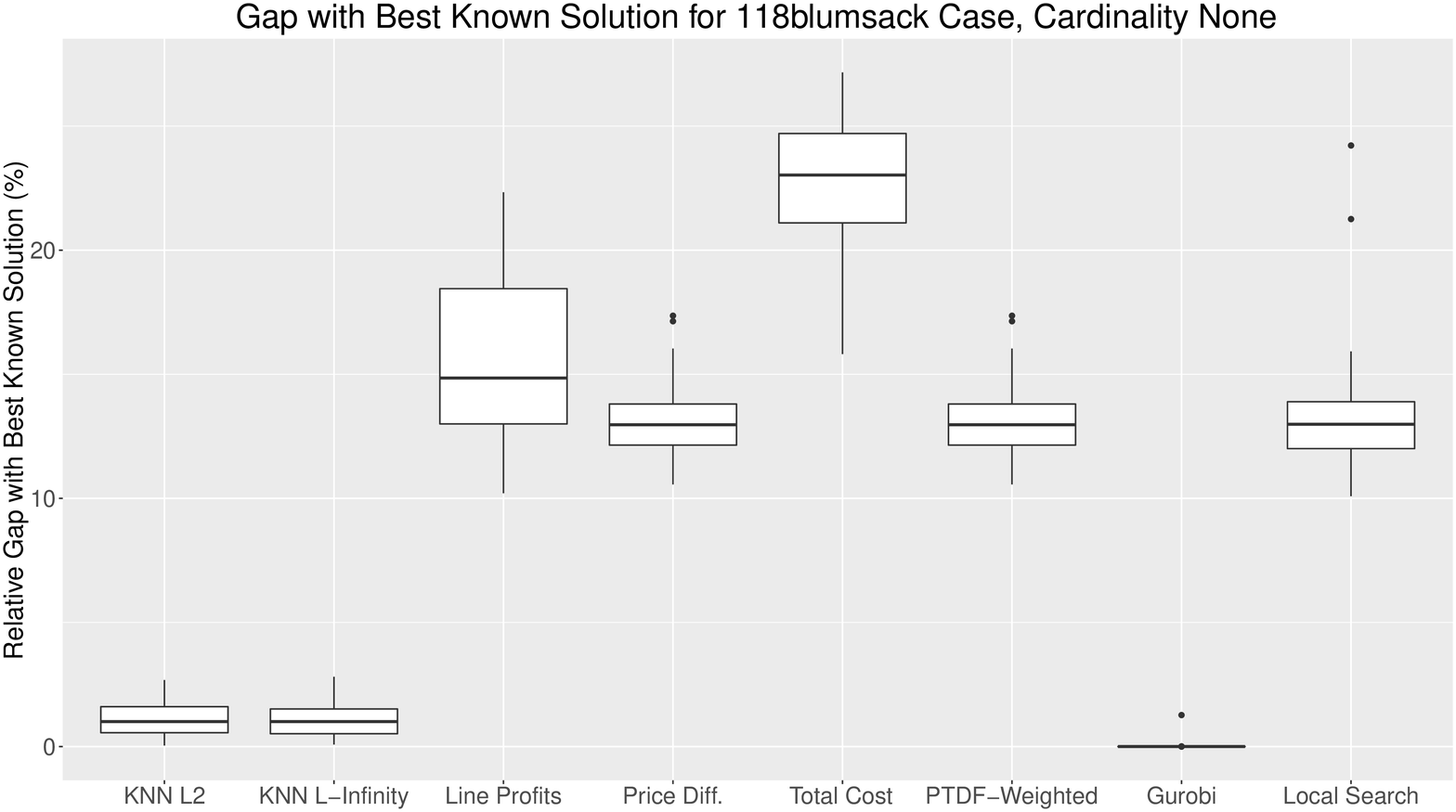} &
		\includegraphics[width=0.5\linewidth]{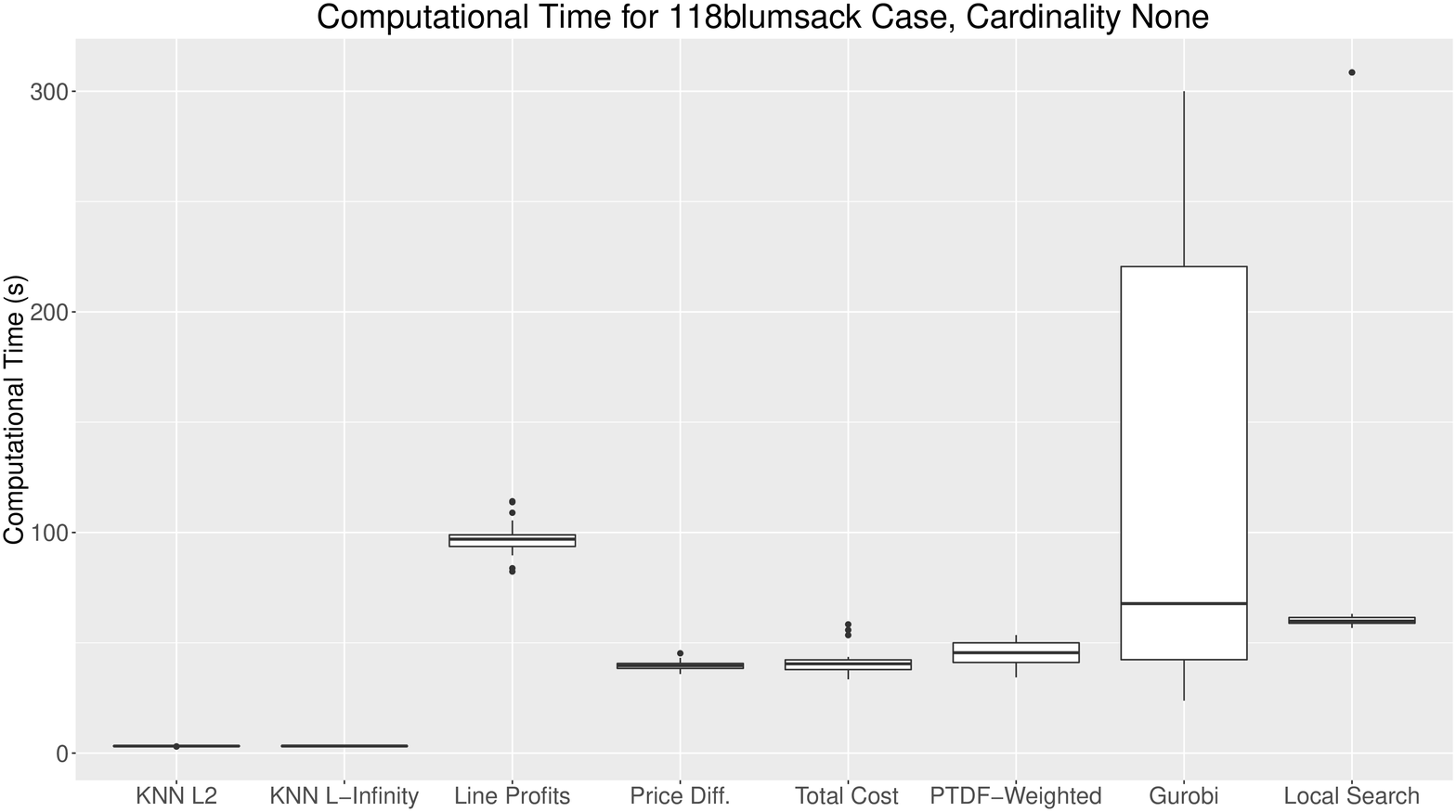}
	\end{tabular}
	\vspace{-0.4cm}
	\caption{Solution quality and computational time results for the 118 bus test case for the three different cardinality options.}
	\label{fig:118-results}
\end{figure*}

\begin{figure*}
	\begin{tabular}{c@{\hskip 0in}c}
		\includegraphics[width=0.5\linewidth]{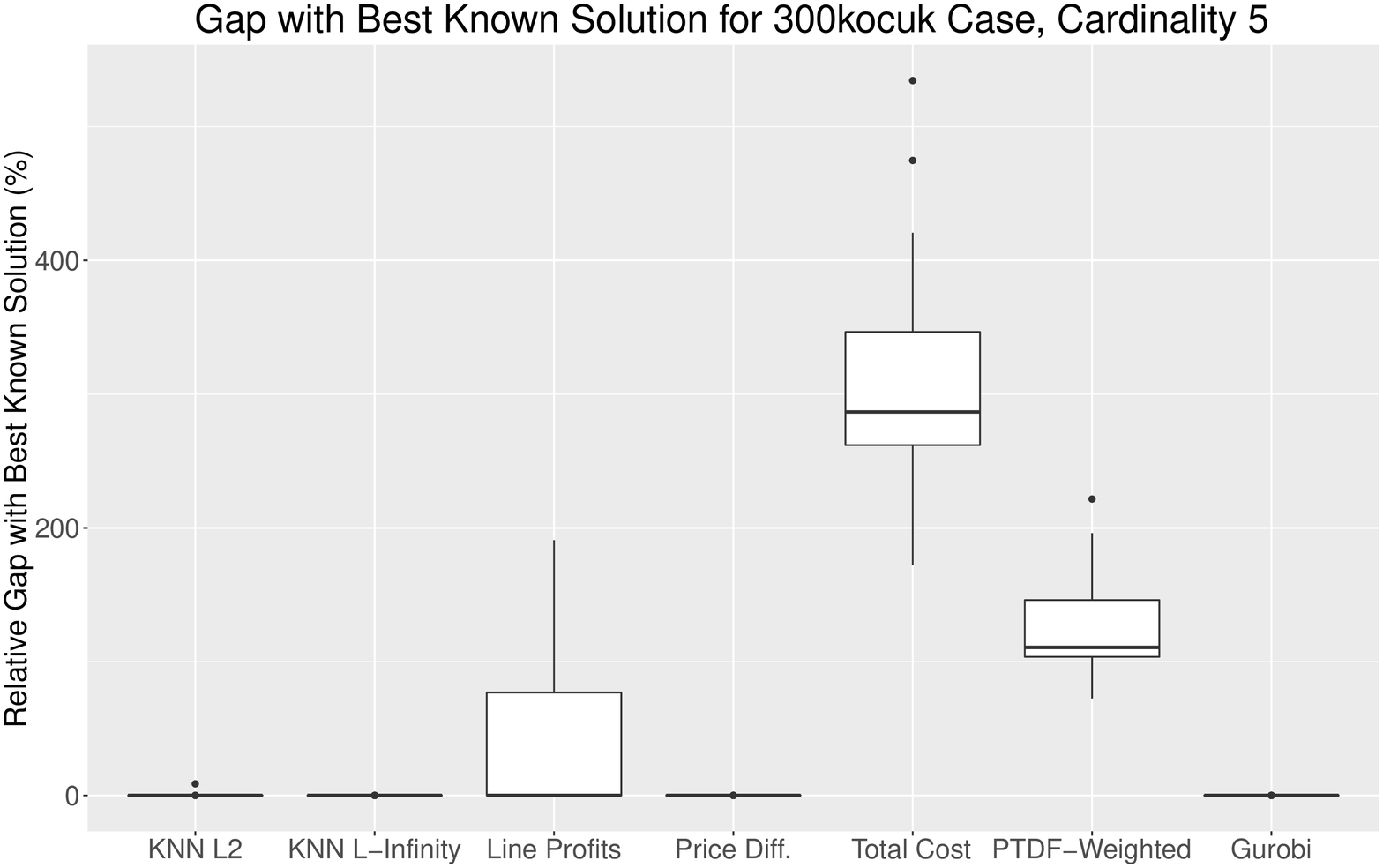} &
		\includegraphics[width=0.5\linewidth]{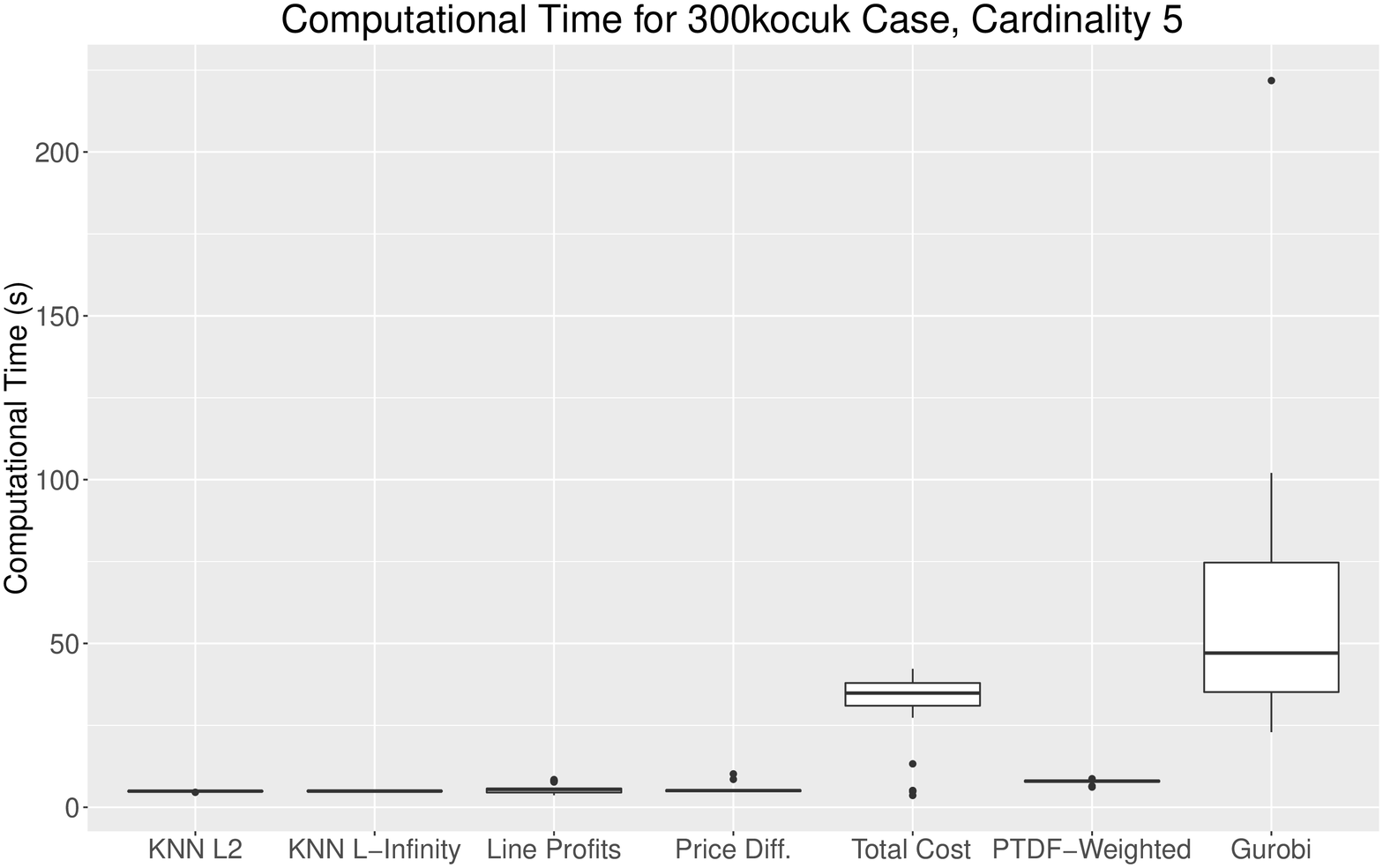} \\
		\includegraphics[width=0.5\linewidth]{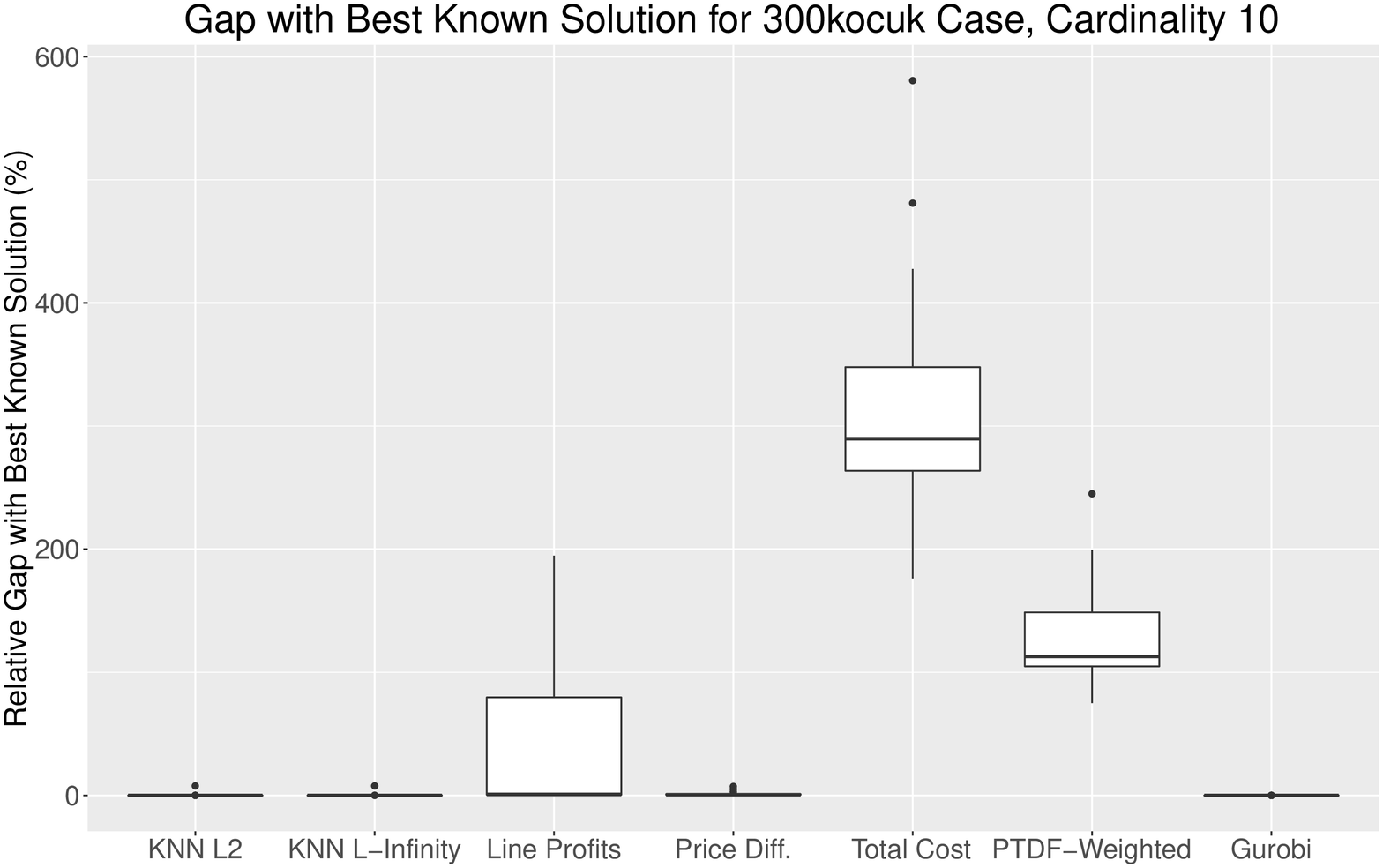} &
		\includegraphics[width=0.5\linewidth]{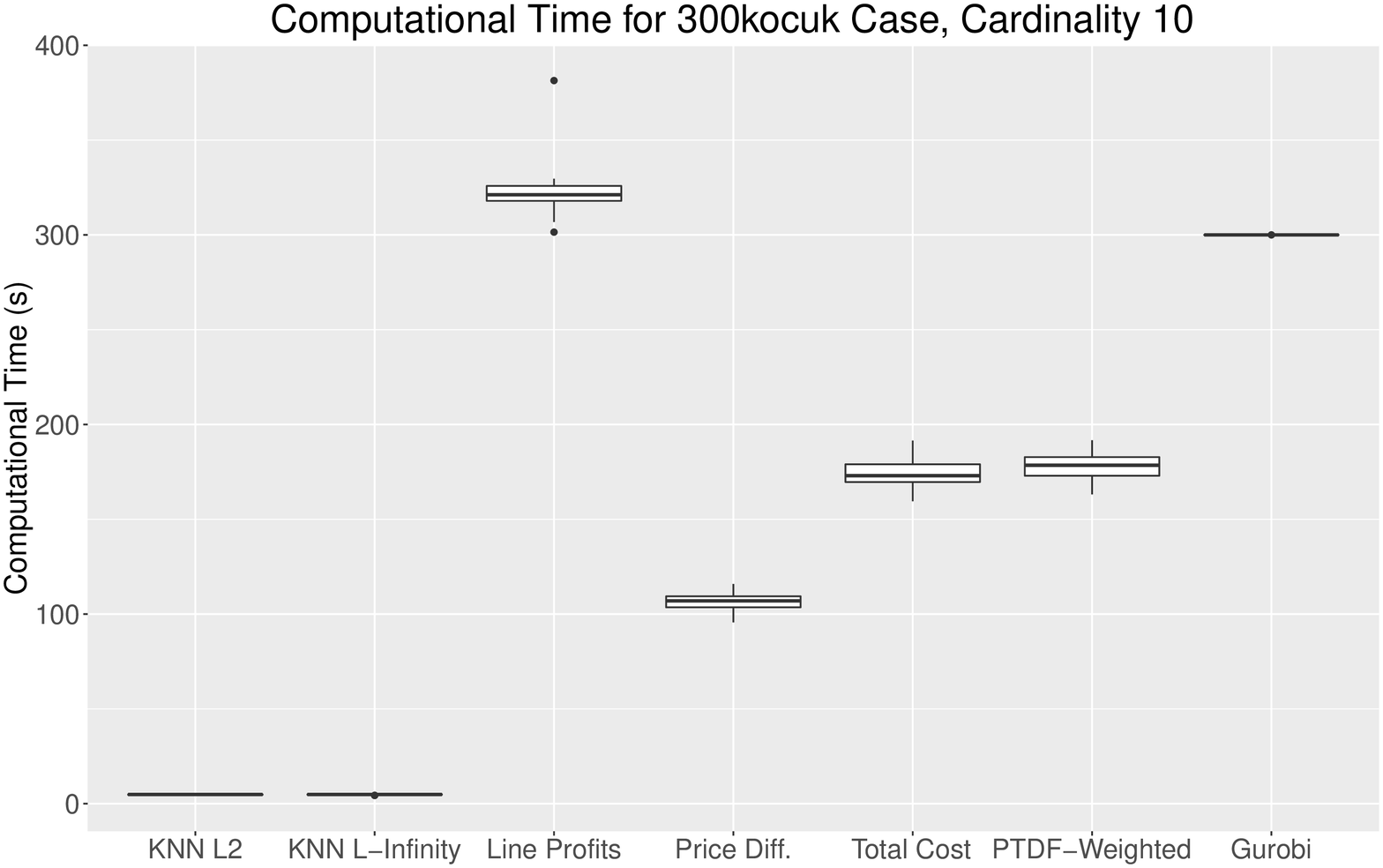} \\
		\includegraphics[width=0.5\linewidth]{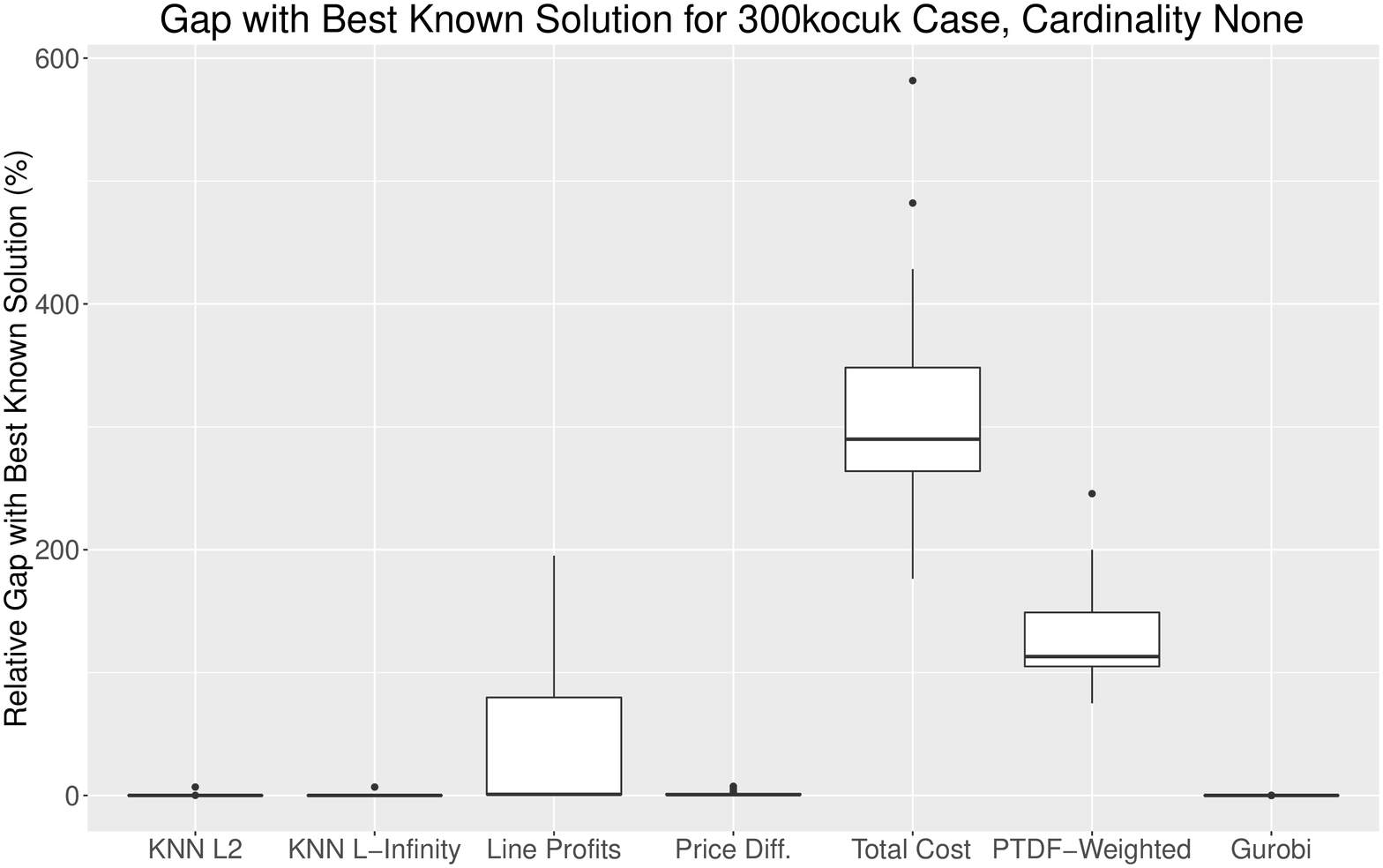} &
		\includegraphics[width=0.5\linewidth]{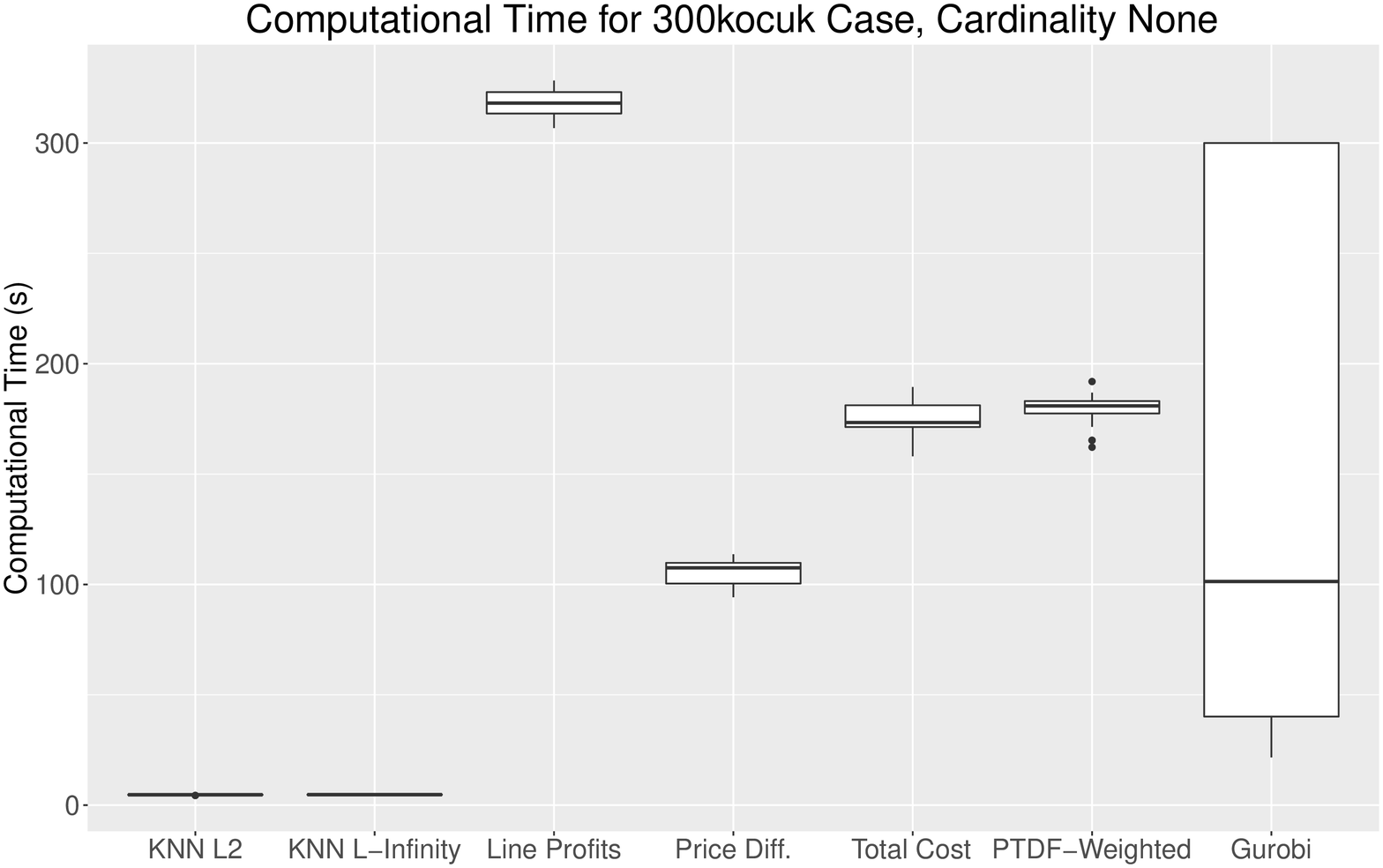}
	\end{tabular}
	\vspace{-0.4cm}
	\caption{Solution quality and computational time results for the 300 bus test case for the three different cardinality options.}
	\label{fig:300-results}
\end{figure*}

\begin{figure*}
	\begin{tabular}{c@{\hskip 0in}c}
		\includegraphics[width=0.5\linewidth]{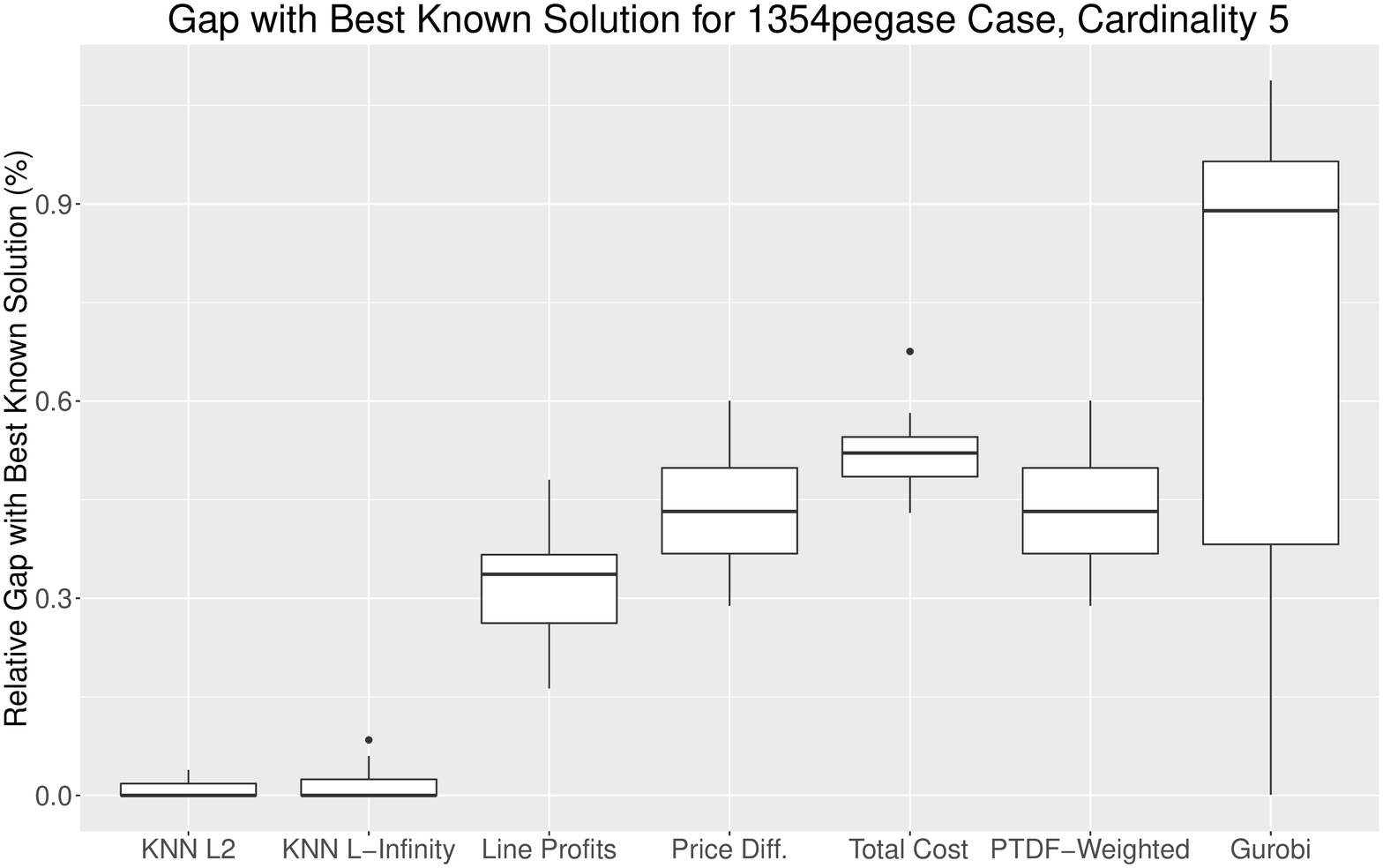} &
		\includegraphics[width=0.5\linewidth]{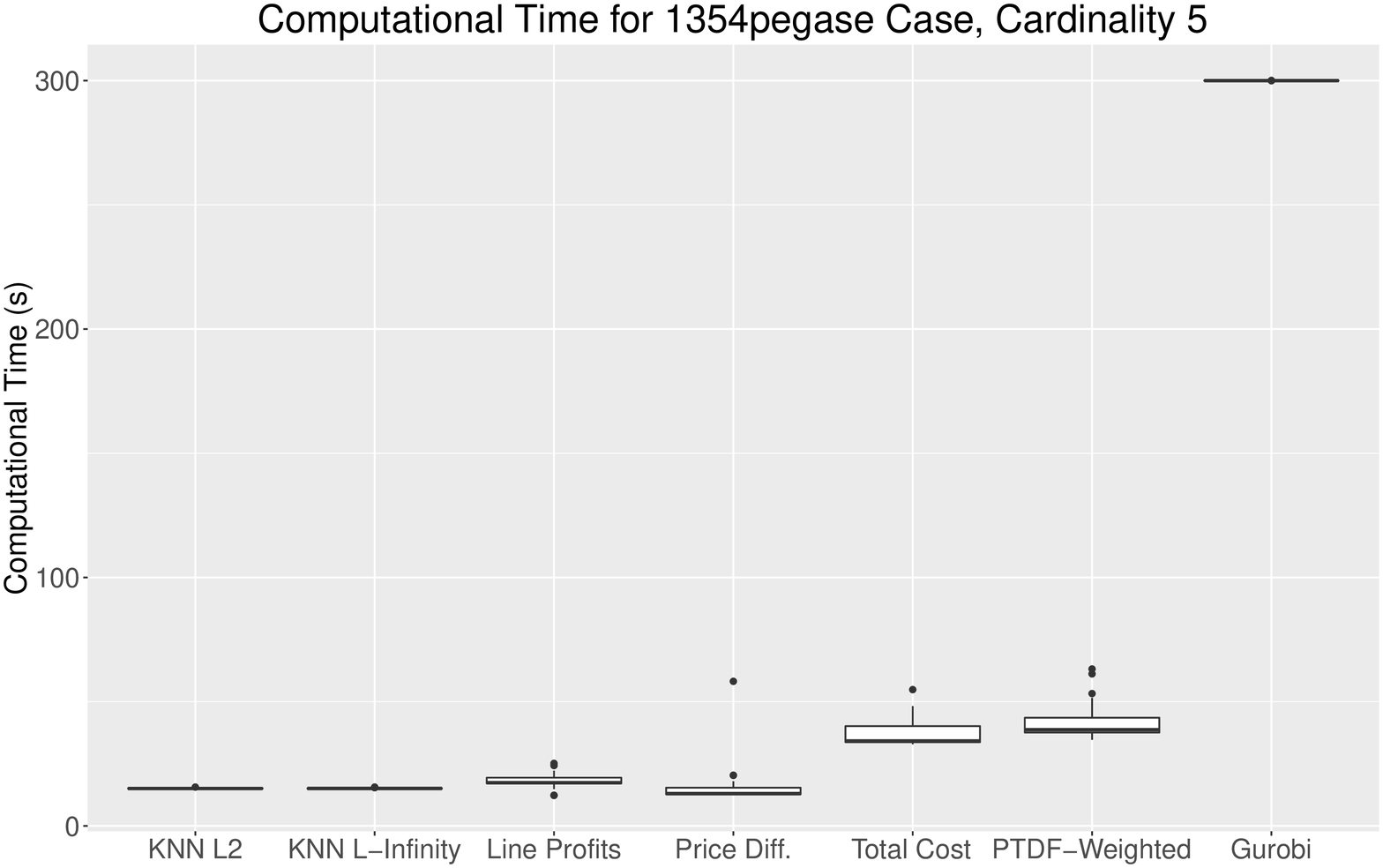} \\
		\includegraphics[width=0.5\linewidth]{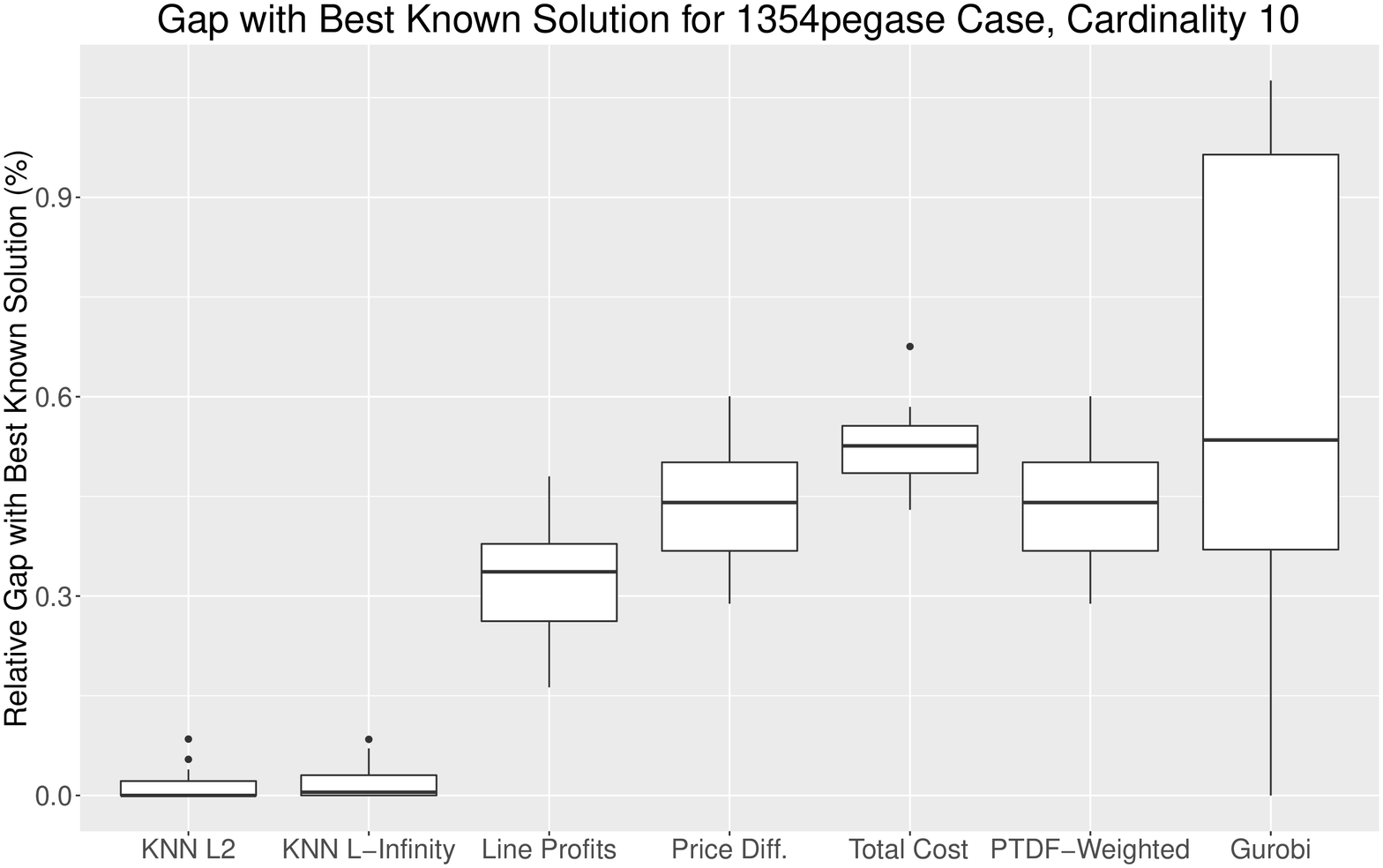} &
		\includegraphics[width=0.5\linewidth]{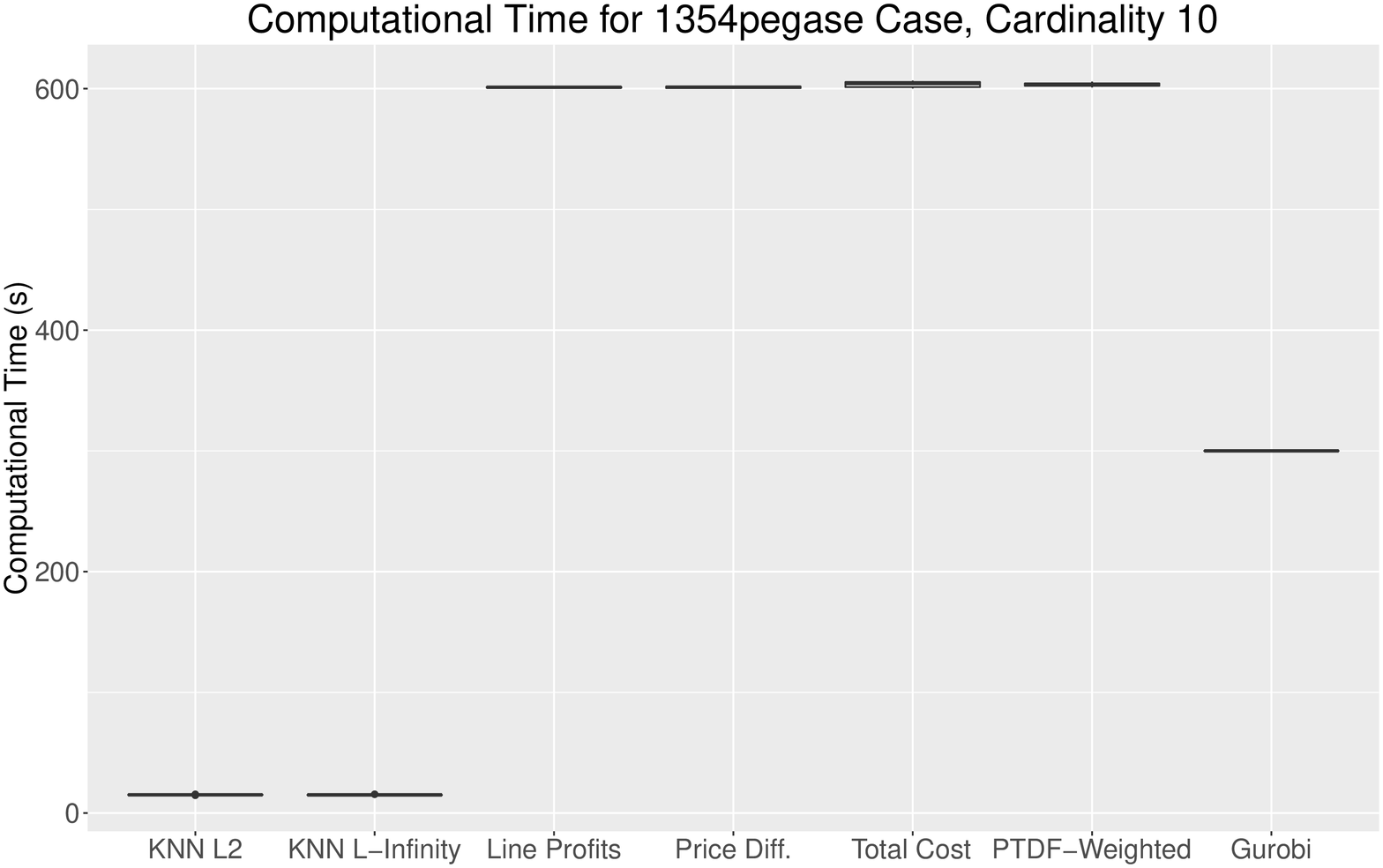} \\
		\includegraphics[width=0.5\linewidth]{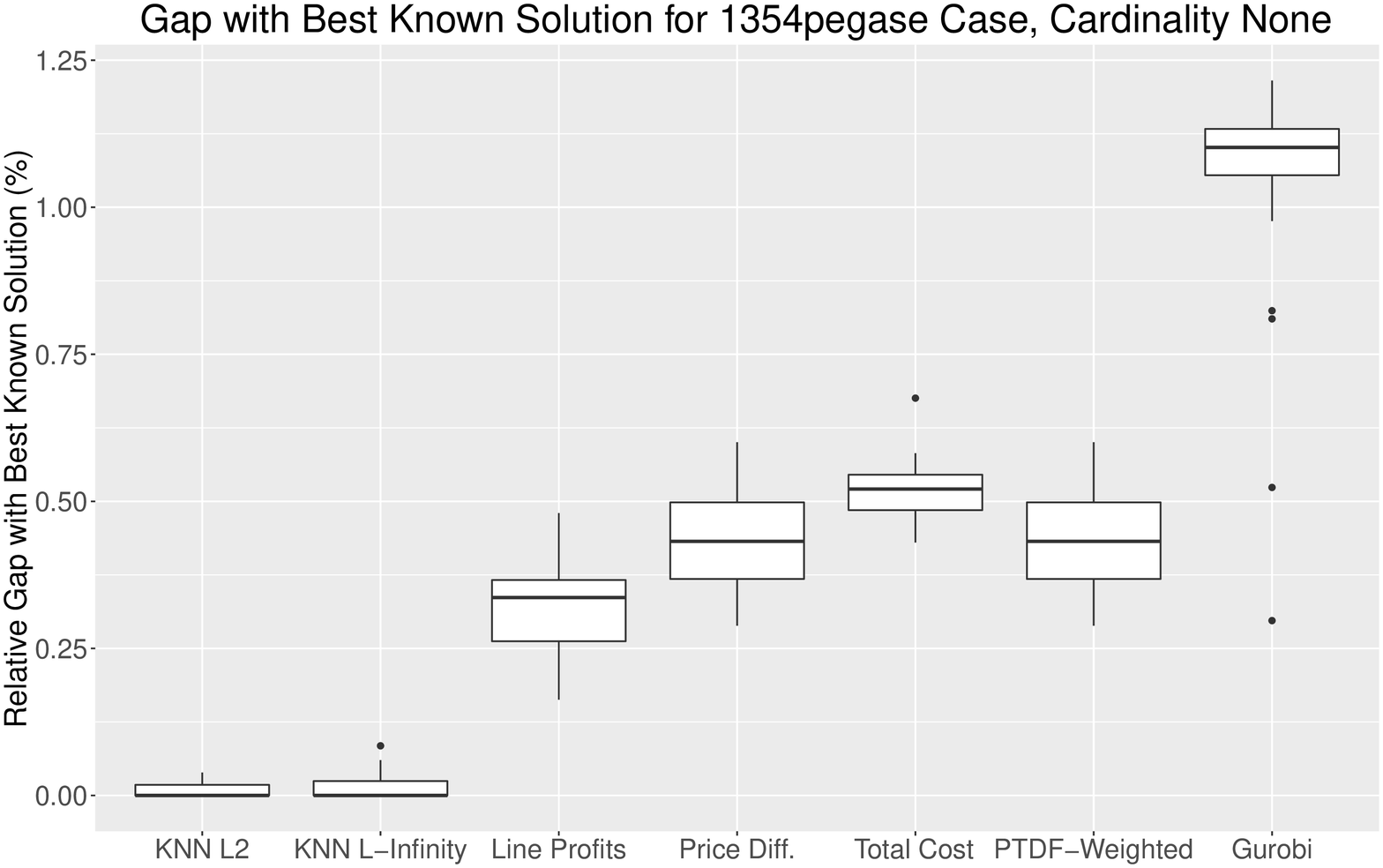} &
		\includegraphics[width=0.5\linewidth]{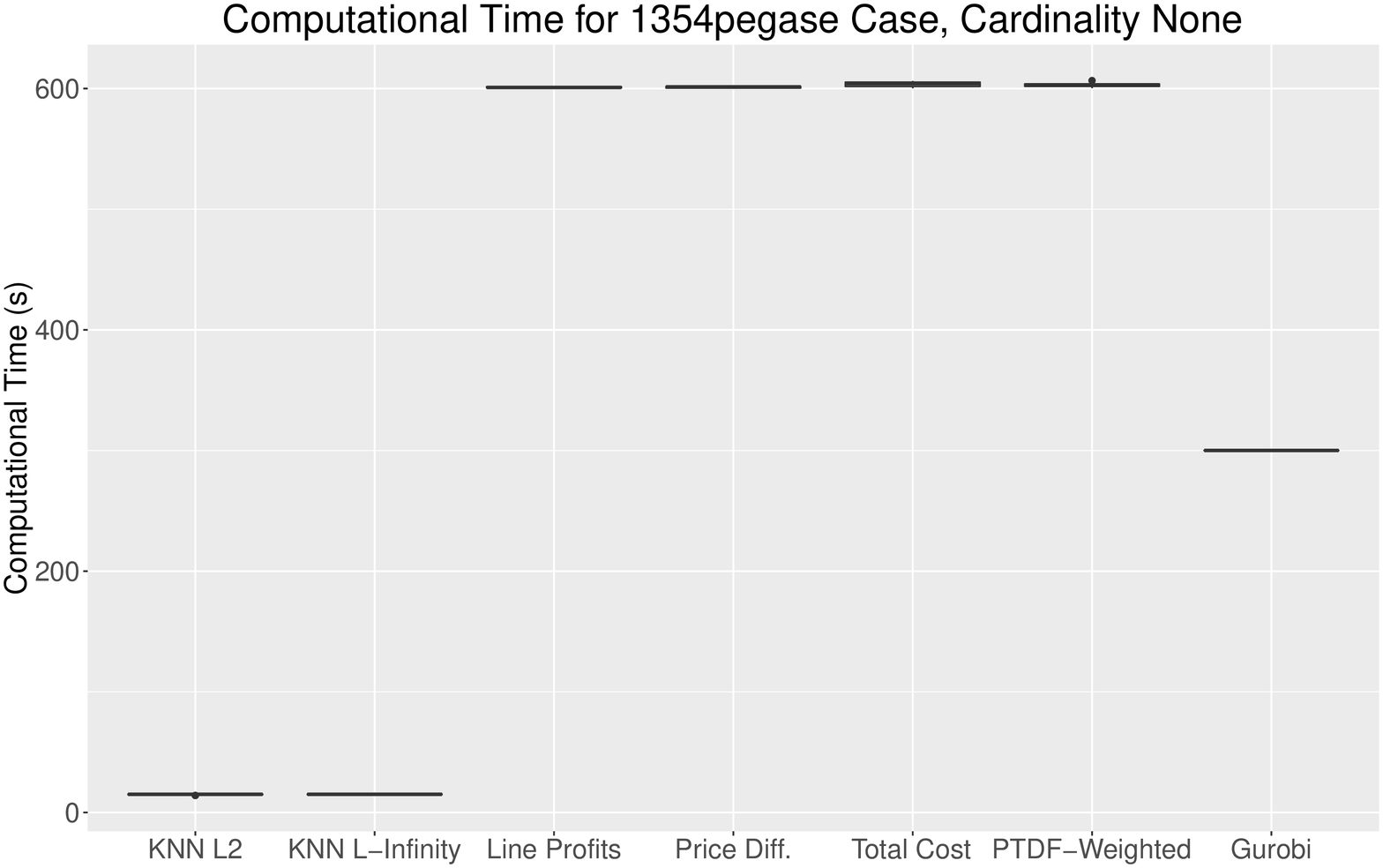}
	\end{tabular}
	\vspace{-0.4cm}
	\caption{Solution quality and computational time results for the 1354 bus test case for the three different cardinality options.}
	\label{fig:1354-results}
\end{figure*}

\begin{figure*}
	\begin{tabular}{c@{\hskip 0in}c}
		\includegraphics[width=0.5\linewidth]{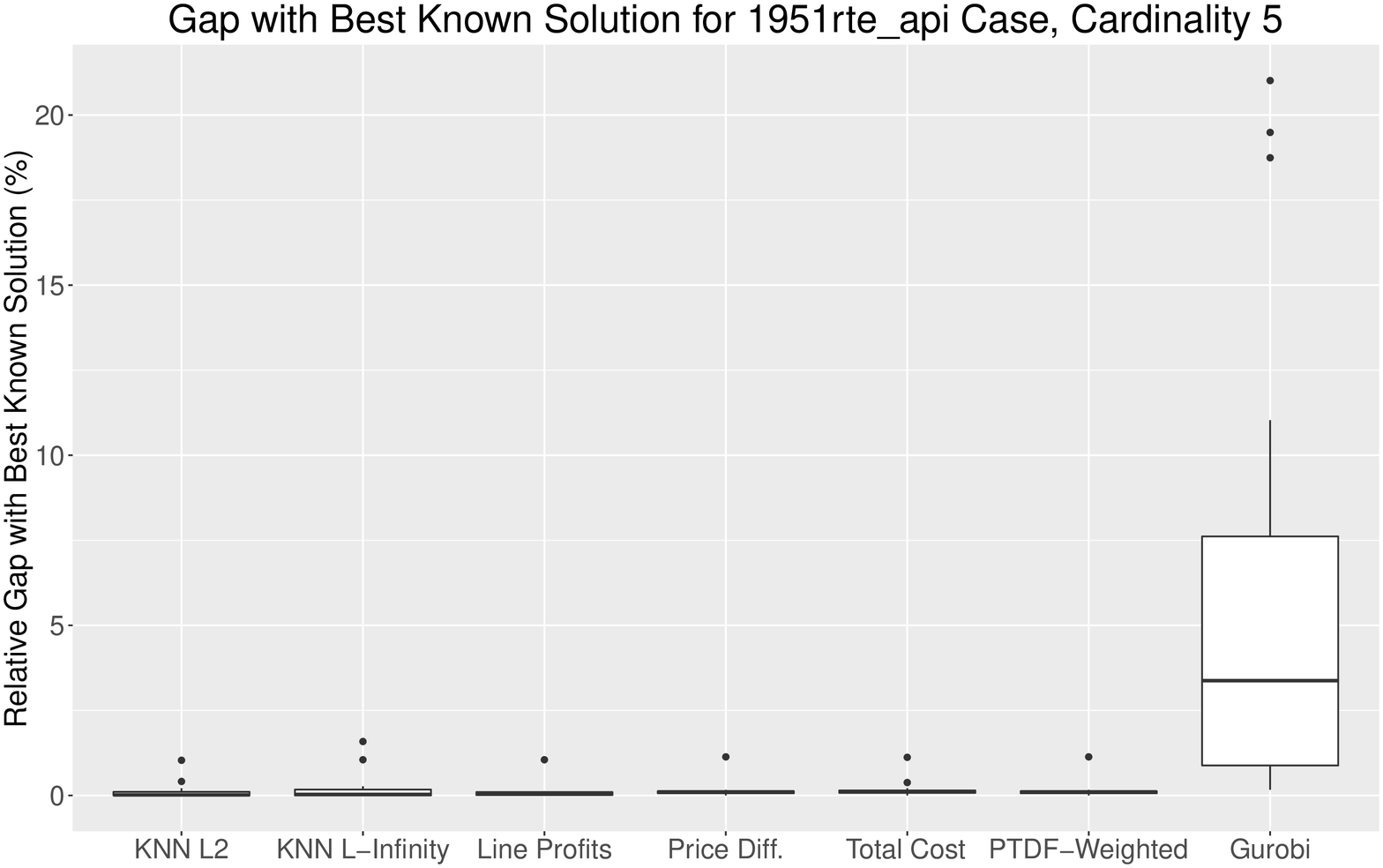} &
		\includegraphics[width=0.5\linewidth]{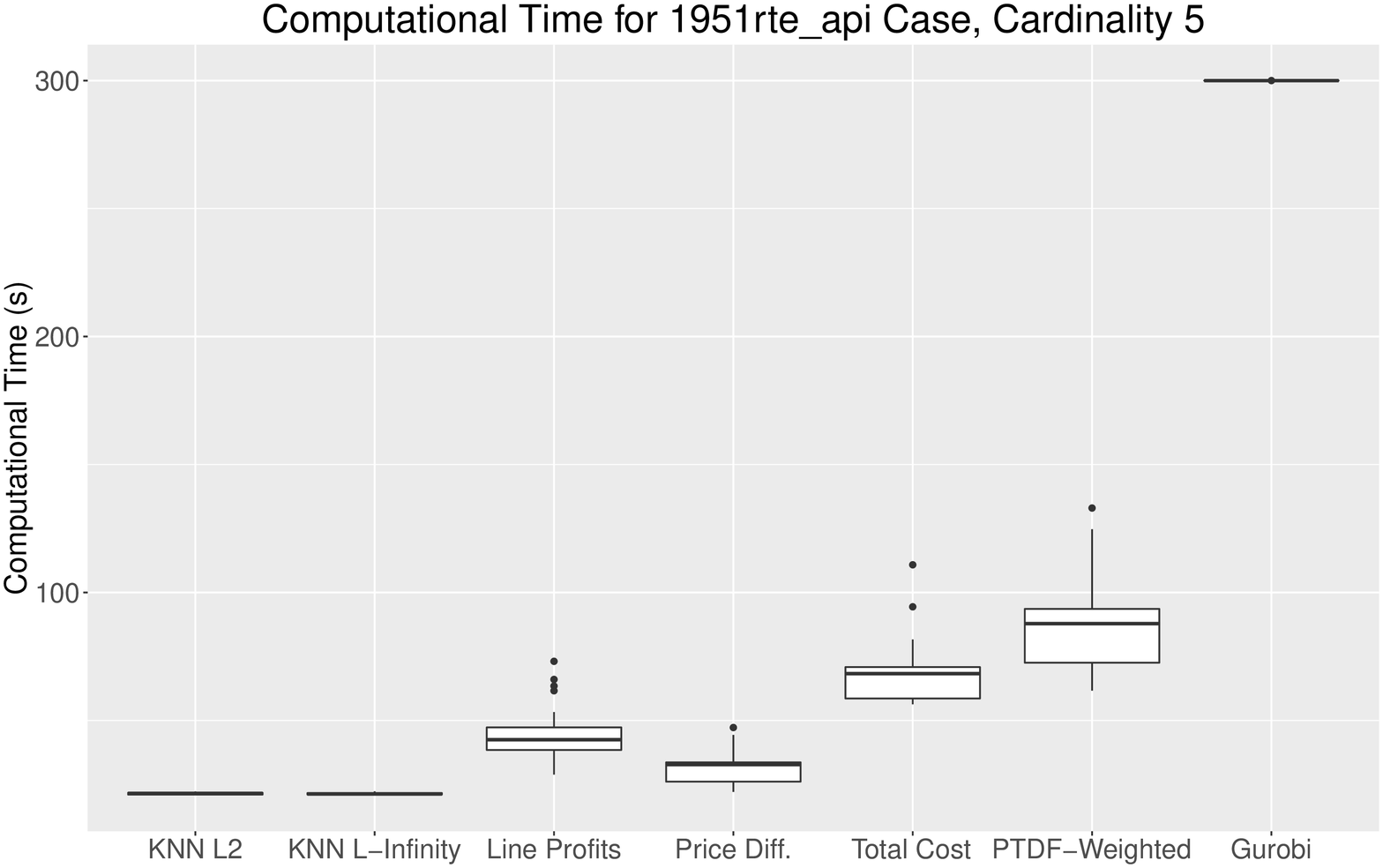} \\
		\includegraphics[width=0.5\linewidth]{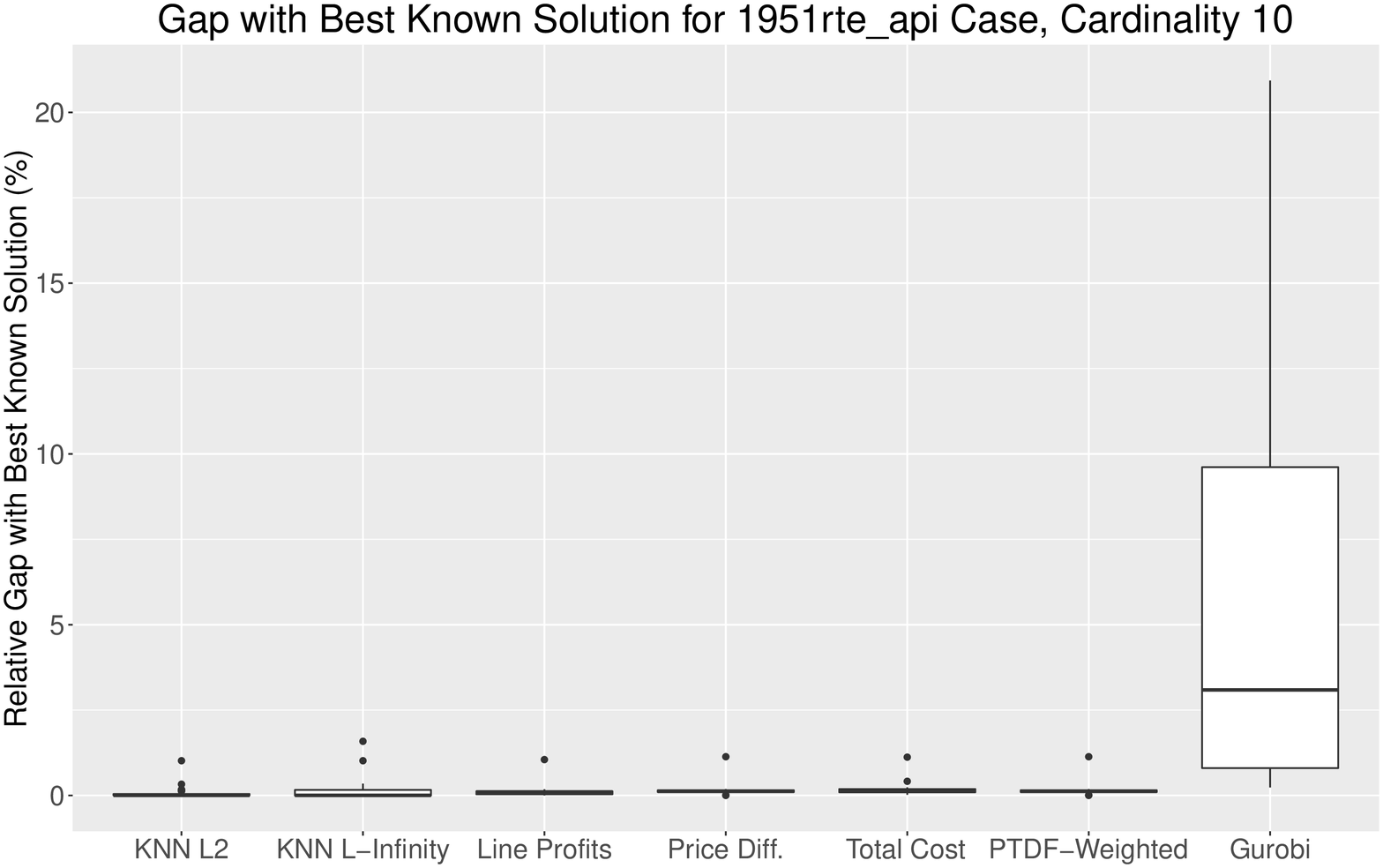} &
		\includegraphics[width=0.5\linewidth]{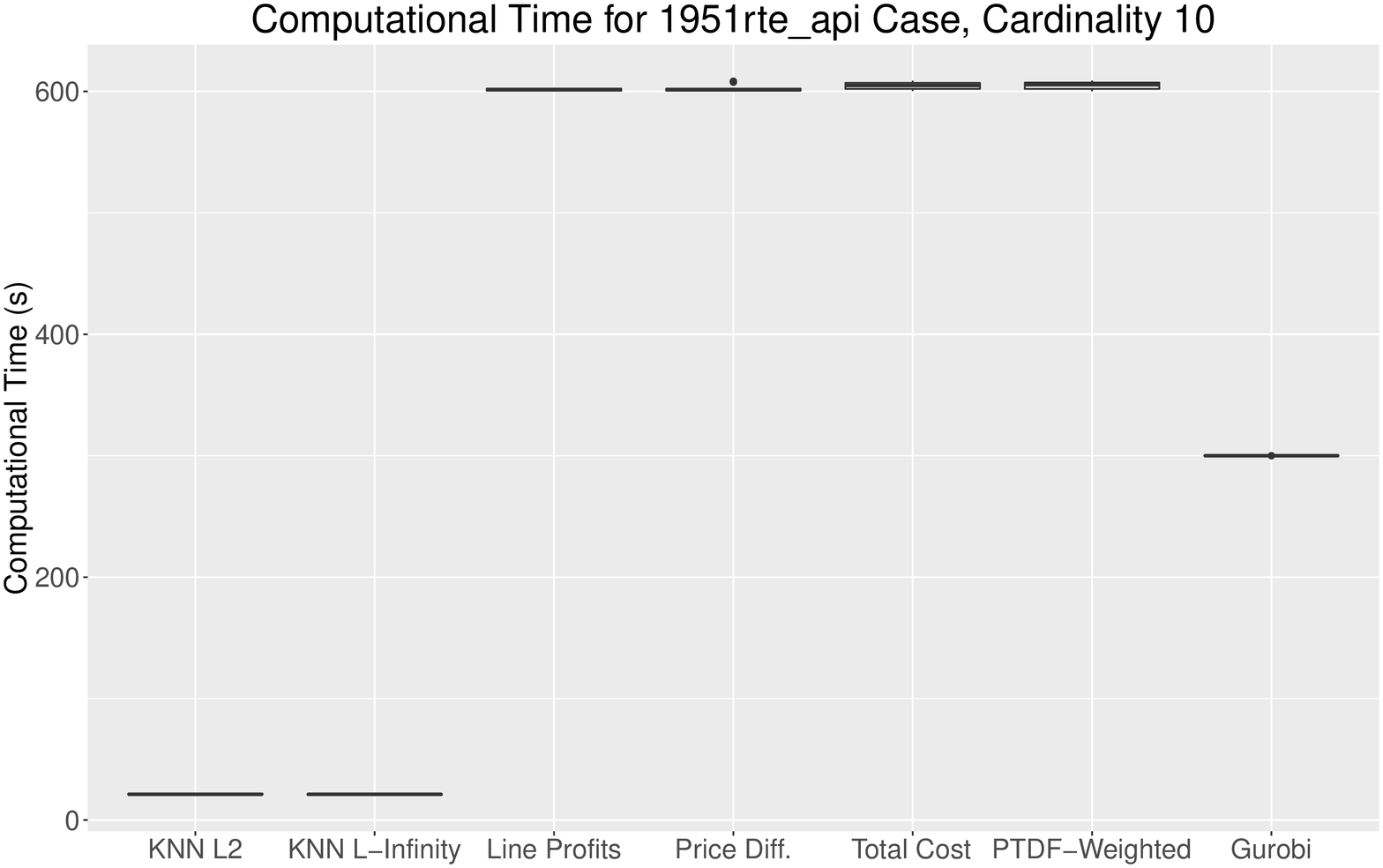} \\
		\includegraphics[width=0.5\linewidth]{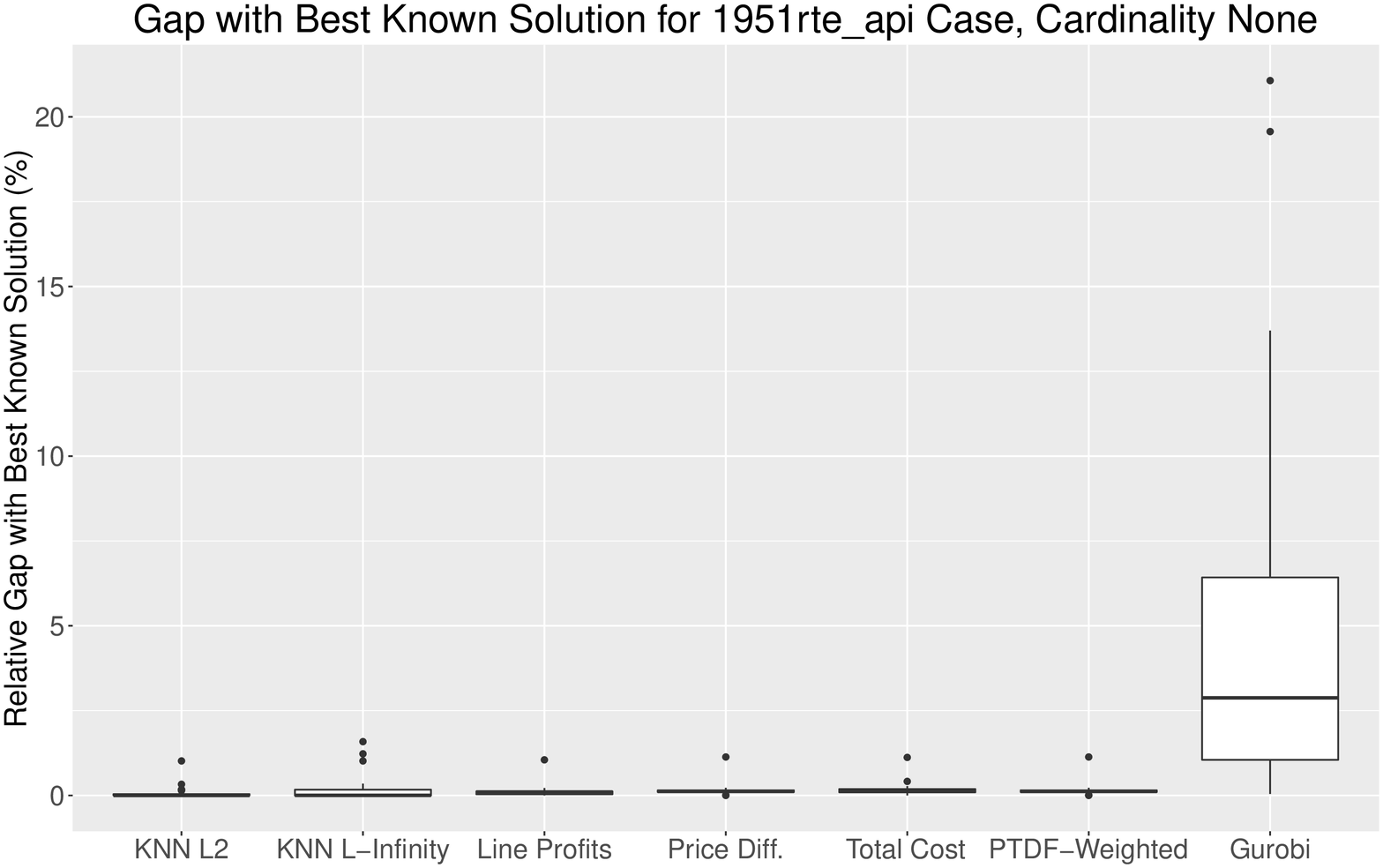} &
		\includegraphics[width=0.5\linewidth]{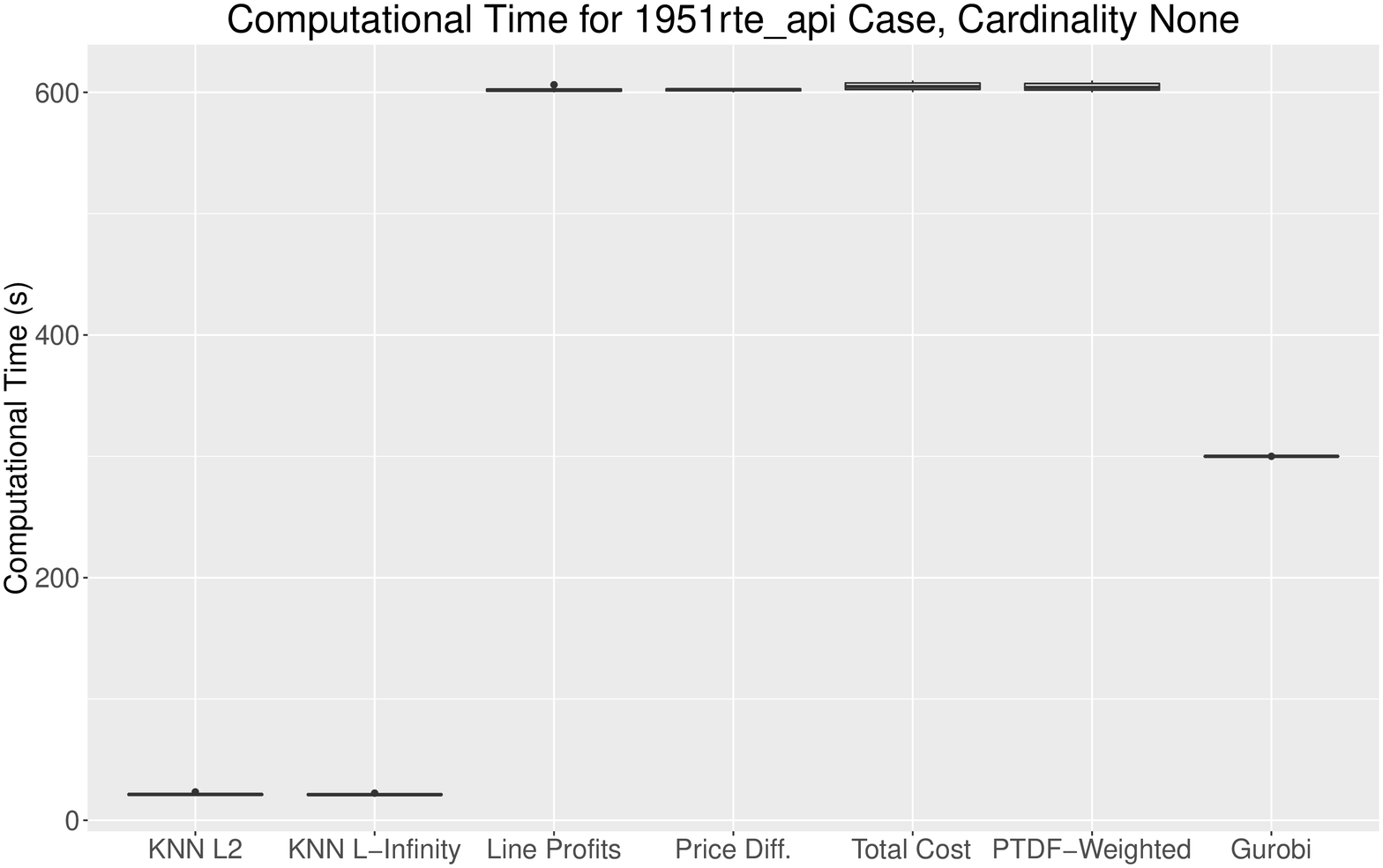}
	\end{tabular}
	\vspace{-0.4cm}
	\caption{Solution quality and computational time results for the 1951 bus test case for the three different cardinality options.}
	\label{fig:1951-results}
\end{figure*}

\begin{figure*}
	\begin{tabular}{c@{\hskip 0in}c}
		\includegraphics[width=0.5\linewidth]{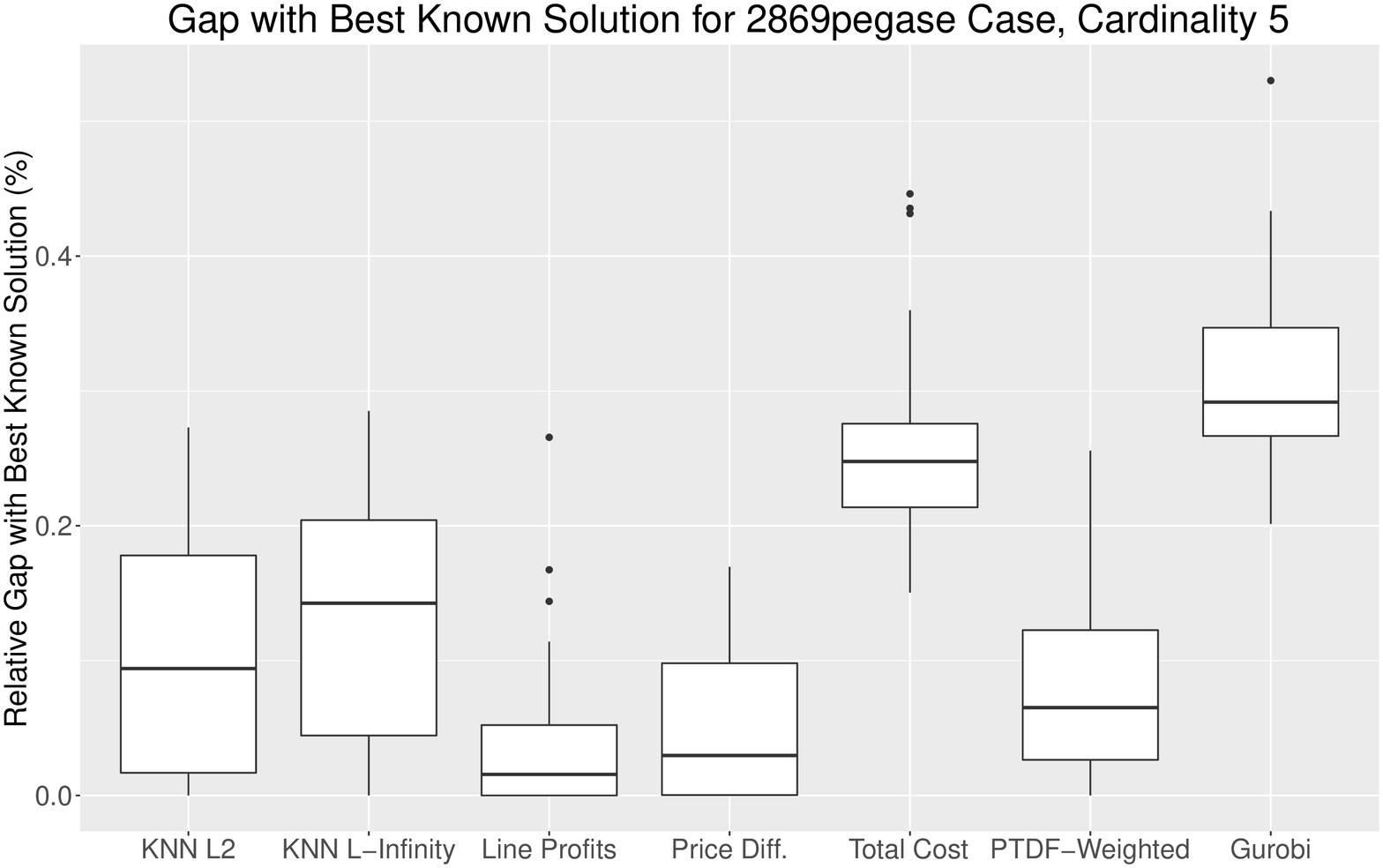} &
		\includegraphics[width=0.5\linewidth]{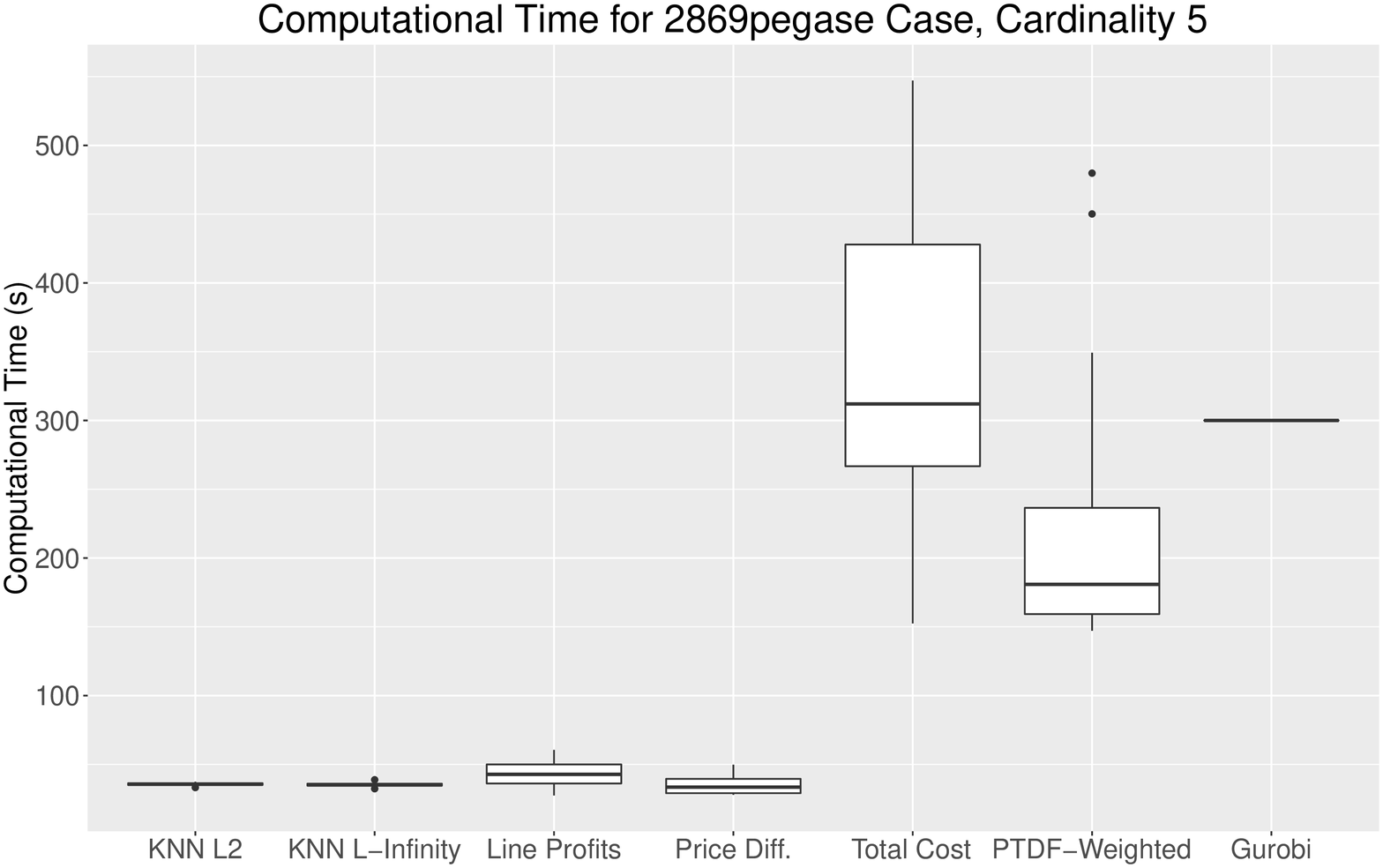} \\
		\includegraphics[width=0.5\linewidth]{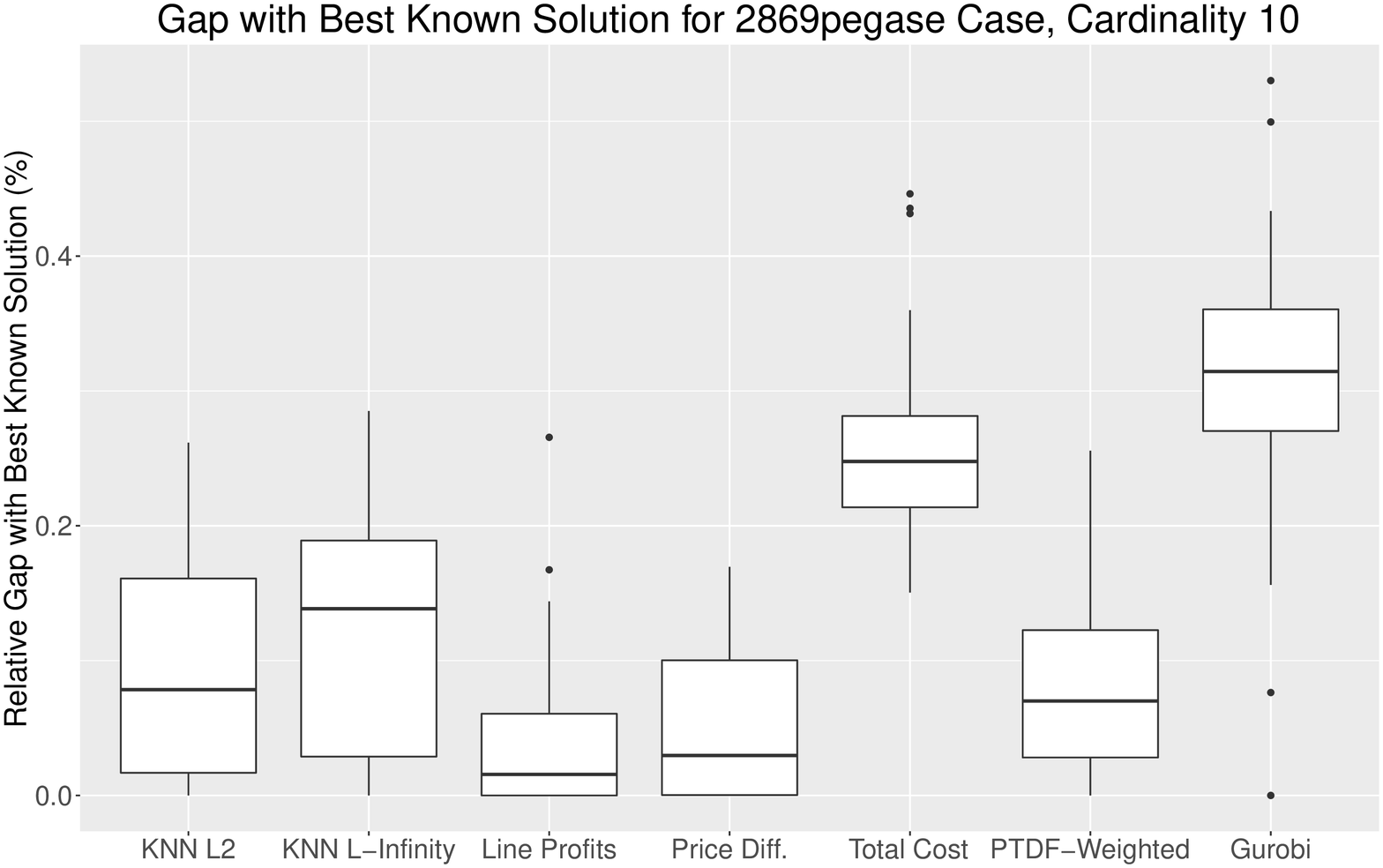} &
		\includegraphics[width=0.5\linewidth]{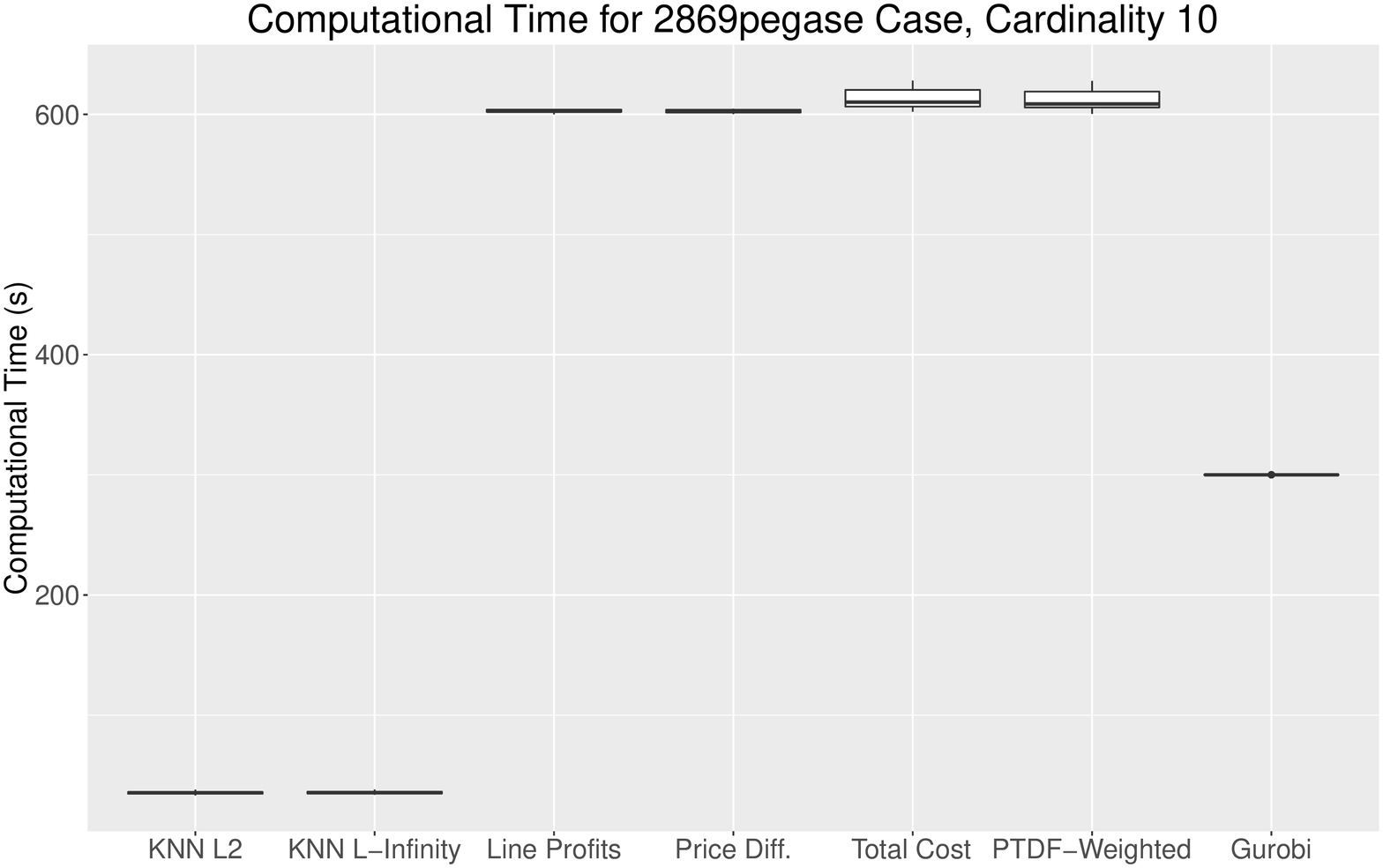} \\
		\includegraphics[width=0.5\linewidth]{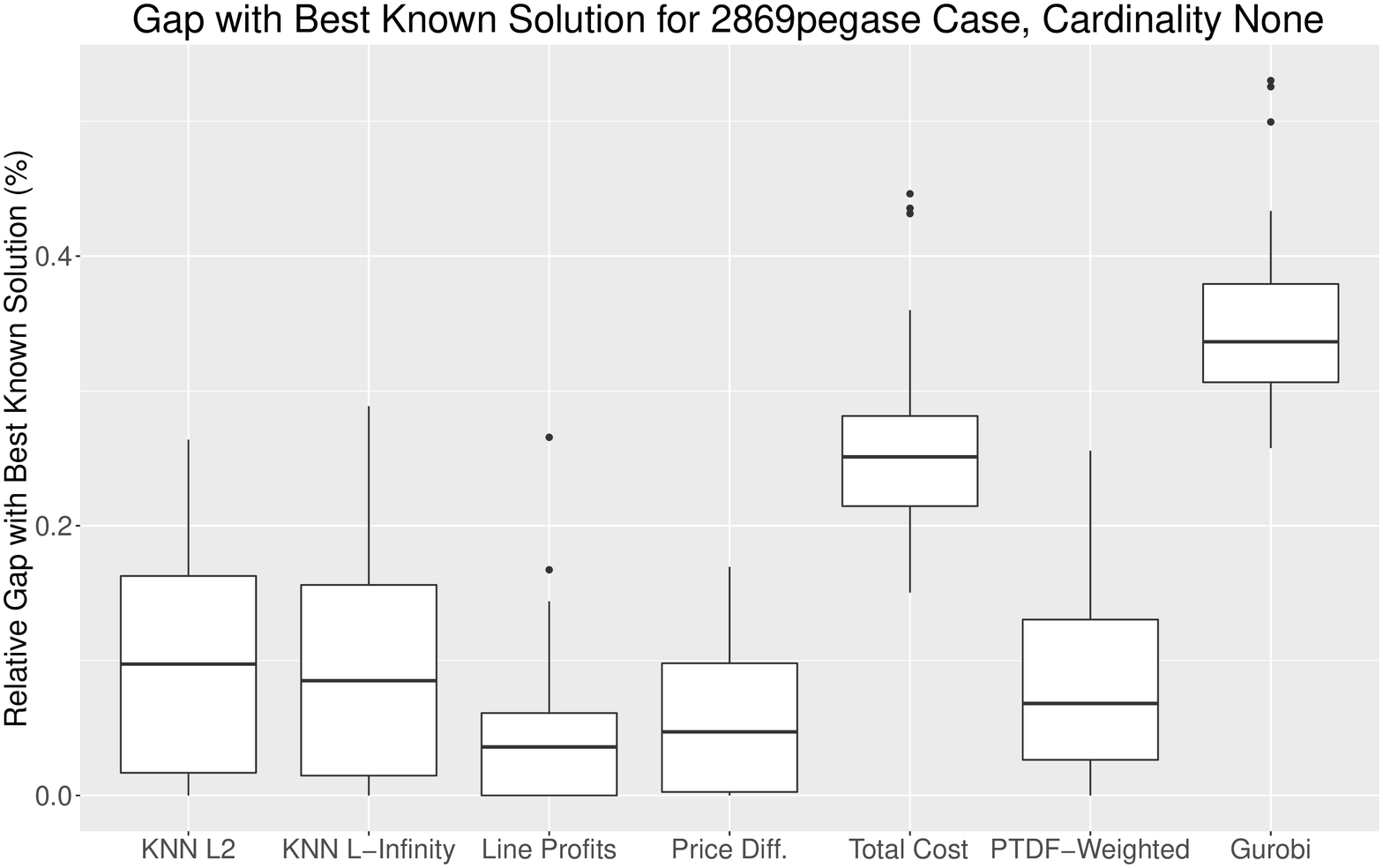} &
		\includegraphics[width=0.5\linewidth]{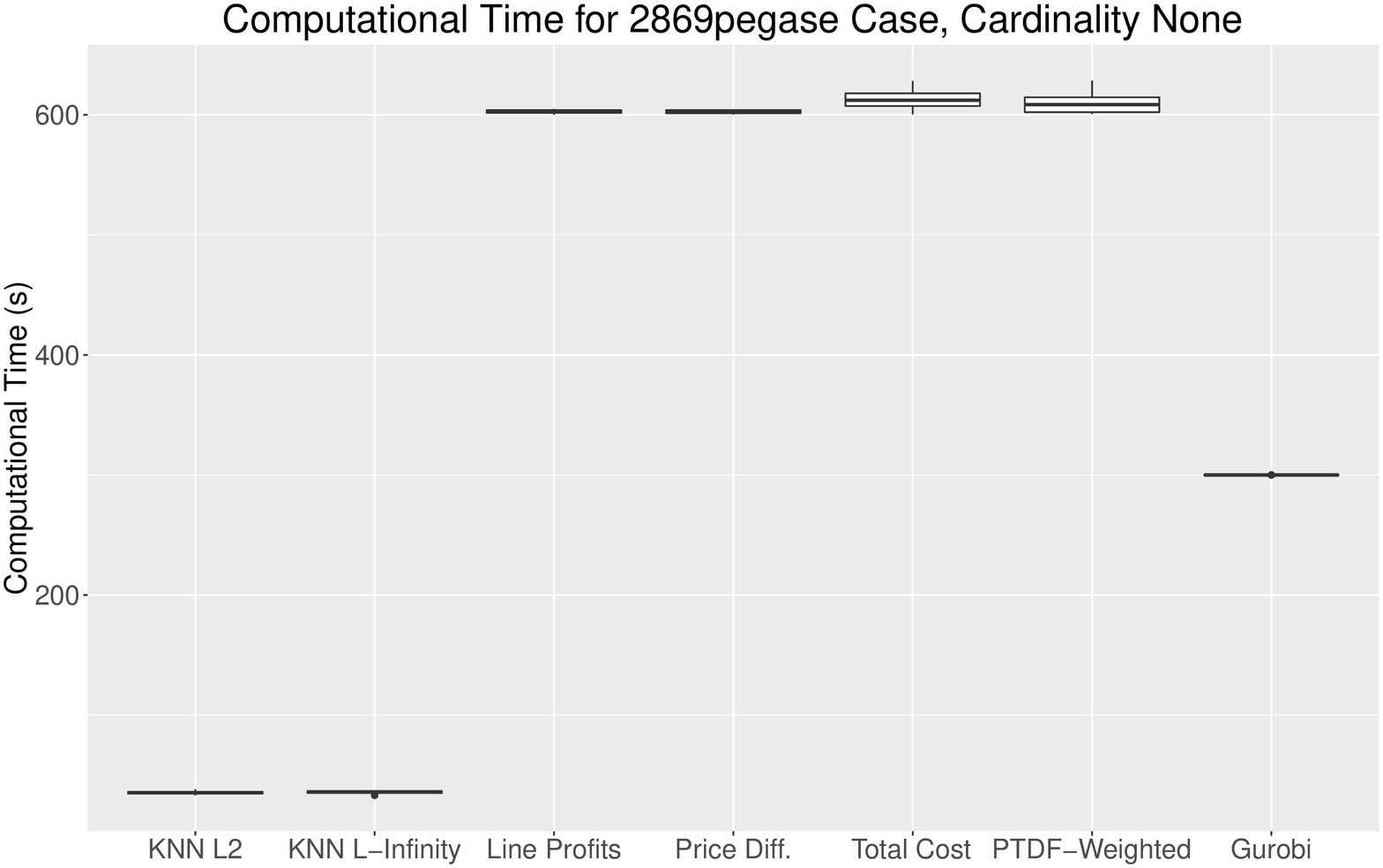}
	\end{tabular}
	\vspace{-0.4cm}
	\caption{Solution quality and computational time results for the 2869 bus test case for the three different cardinality options.}
	\label{fig:2869-results}
\end{figure*}

\begin{figure*}
	\begin{tabular}{c@{\hskip 0in}c}
		\includegraphics[width=0.5\linewidth]{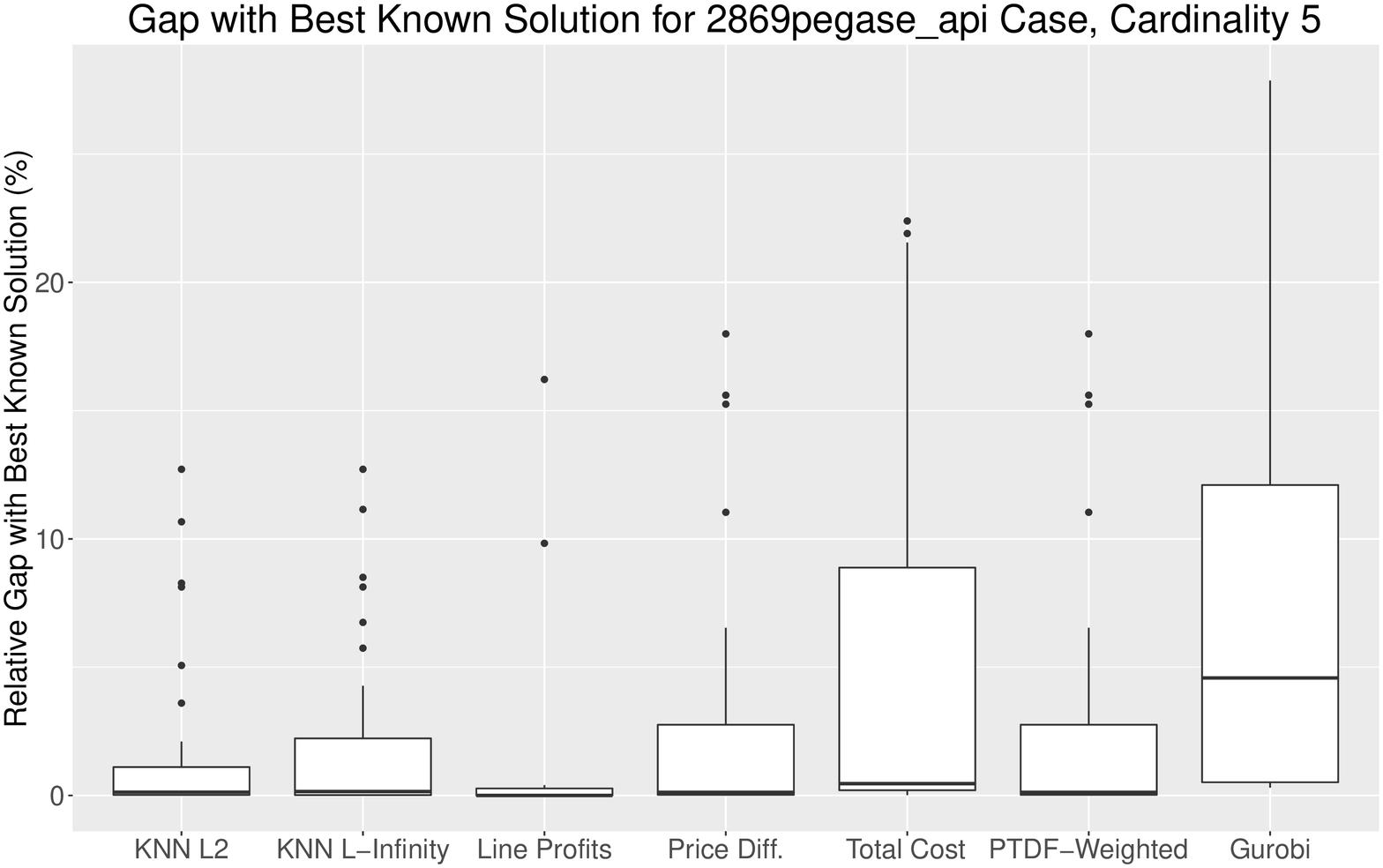} &
		\includegraphics[width=0.5\linewidth]{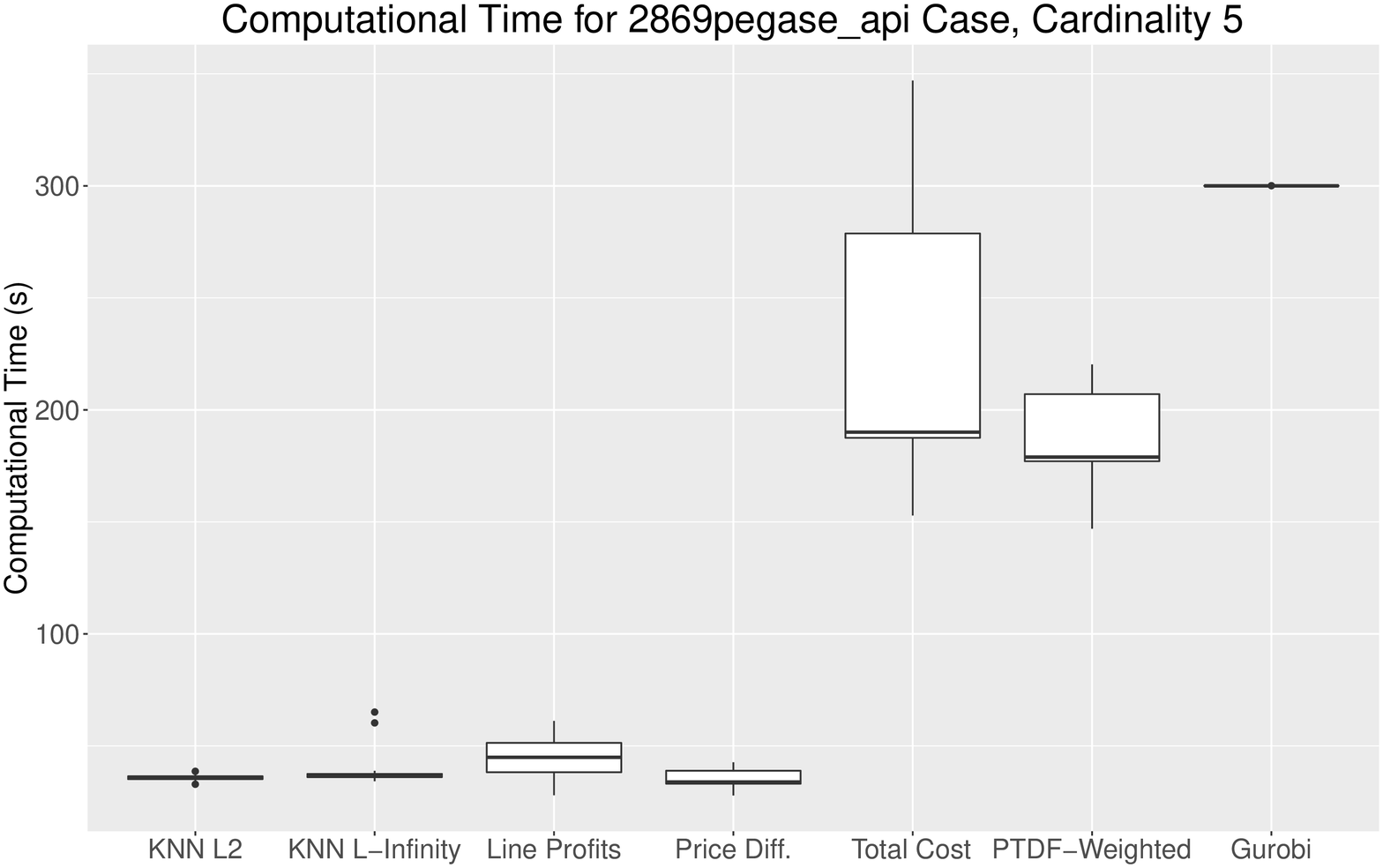} \\
		\includegraphics[width=0.5\linewidth]{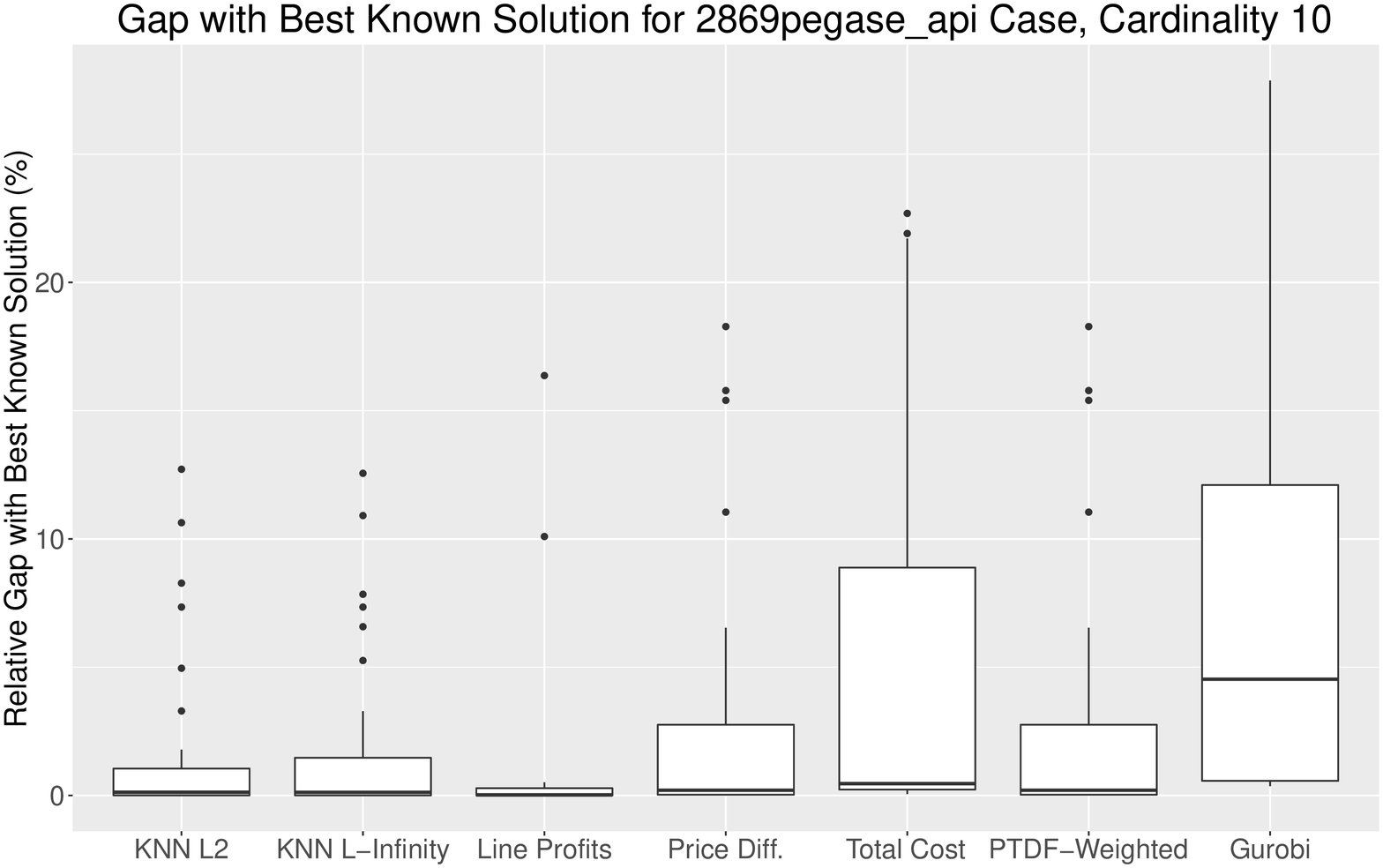} &
		\includegraphics[width=0.5\linewidth]{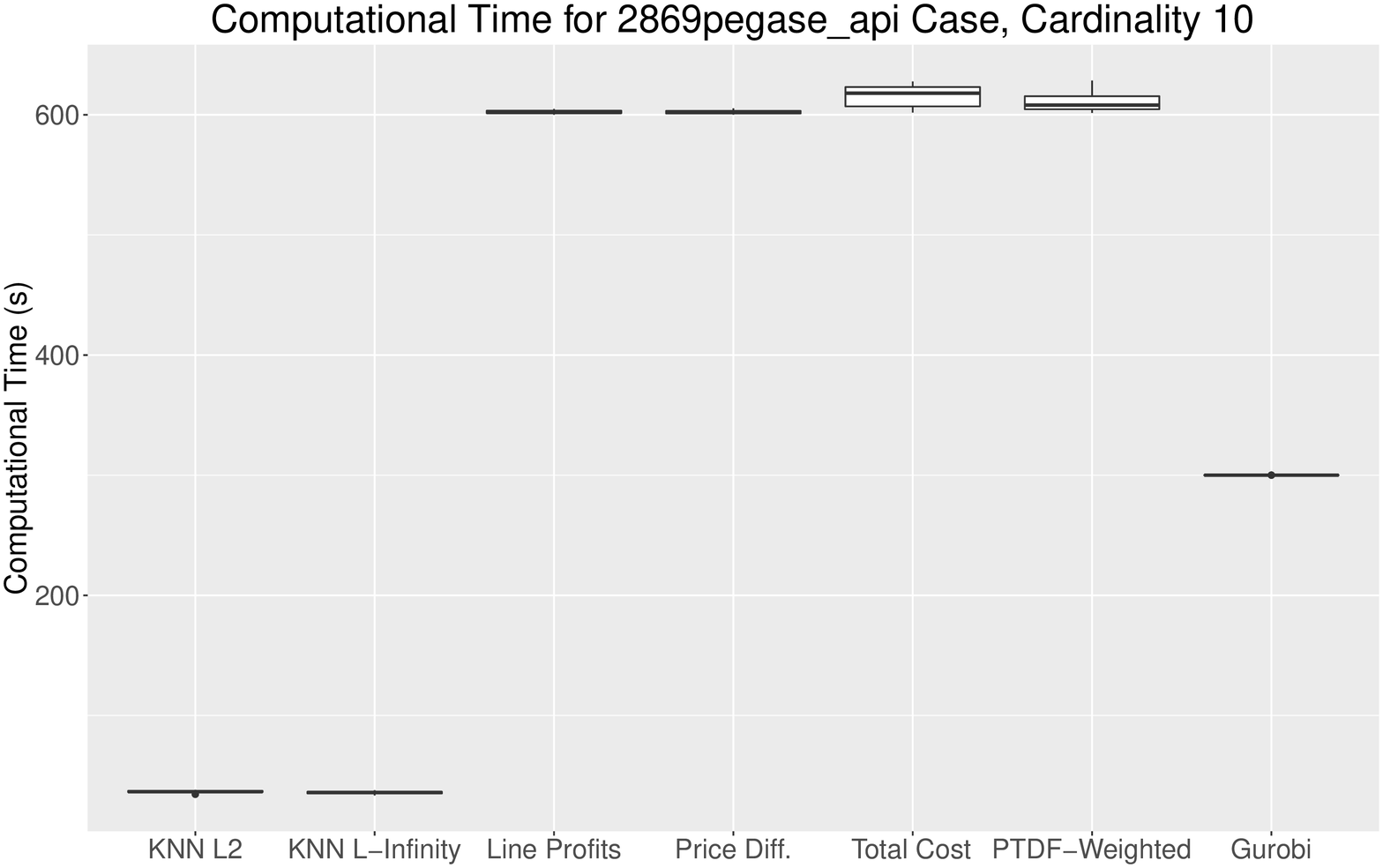} \\
		\includegraphics[width=0.5\linewidth]{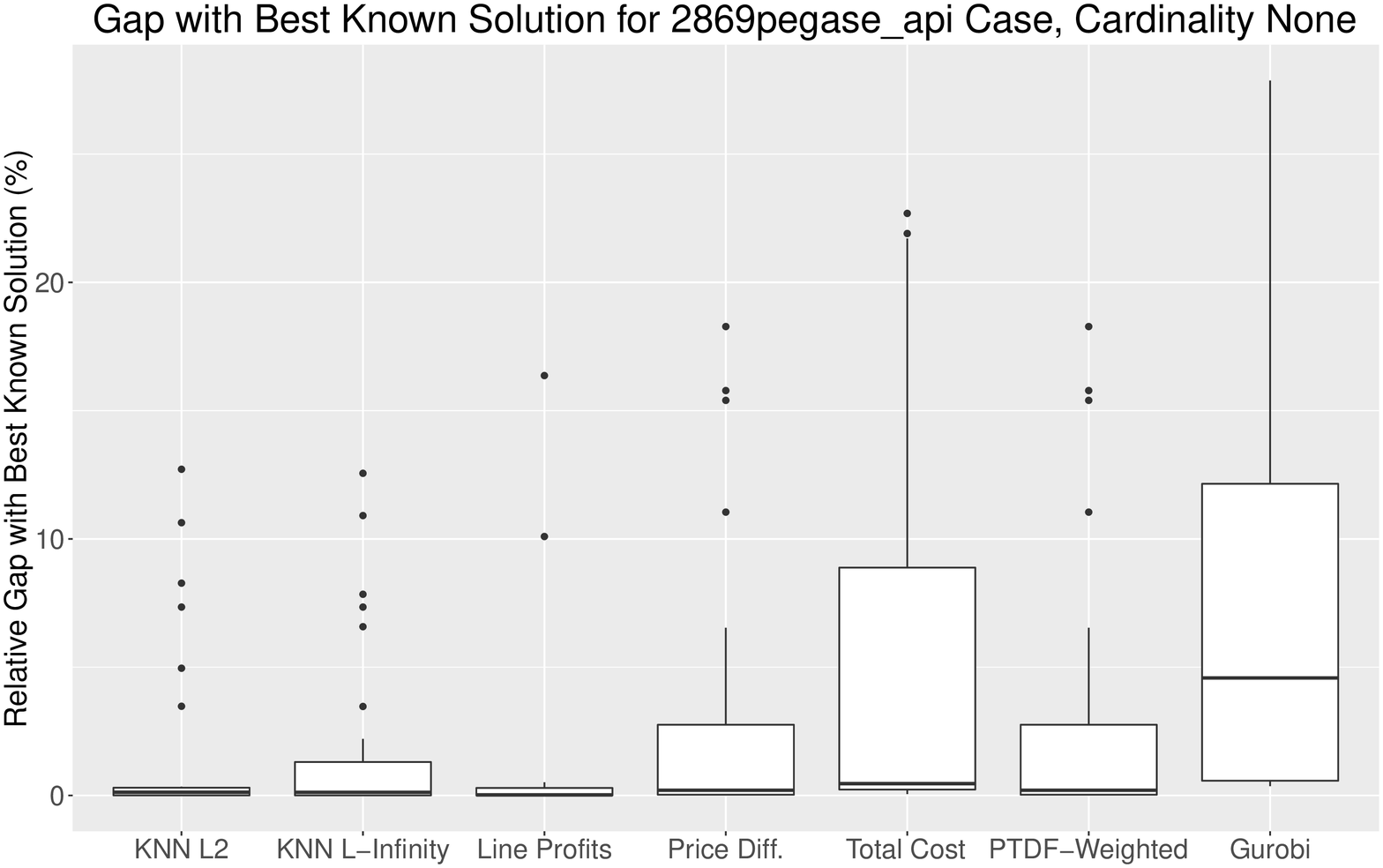} &
		\includegraphics[width=0.5\linewidth]{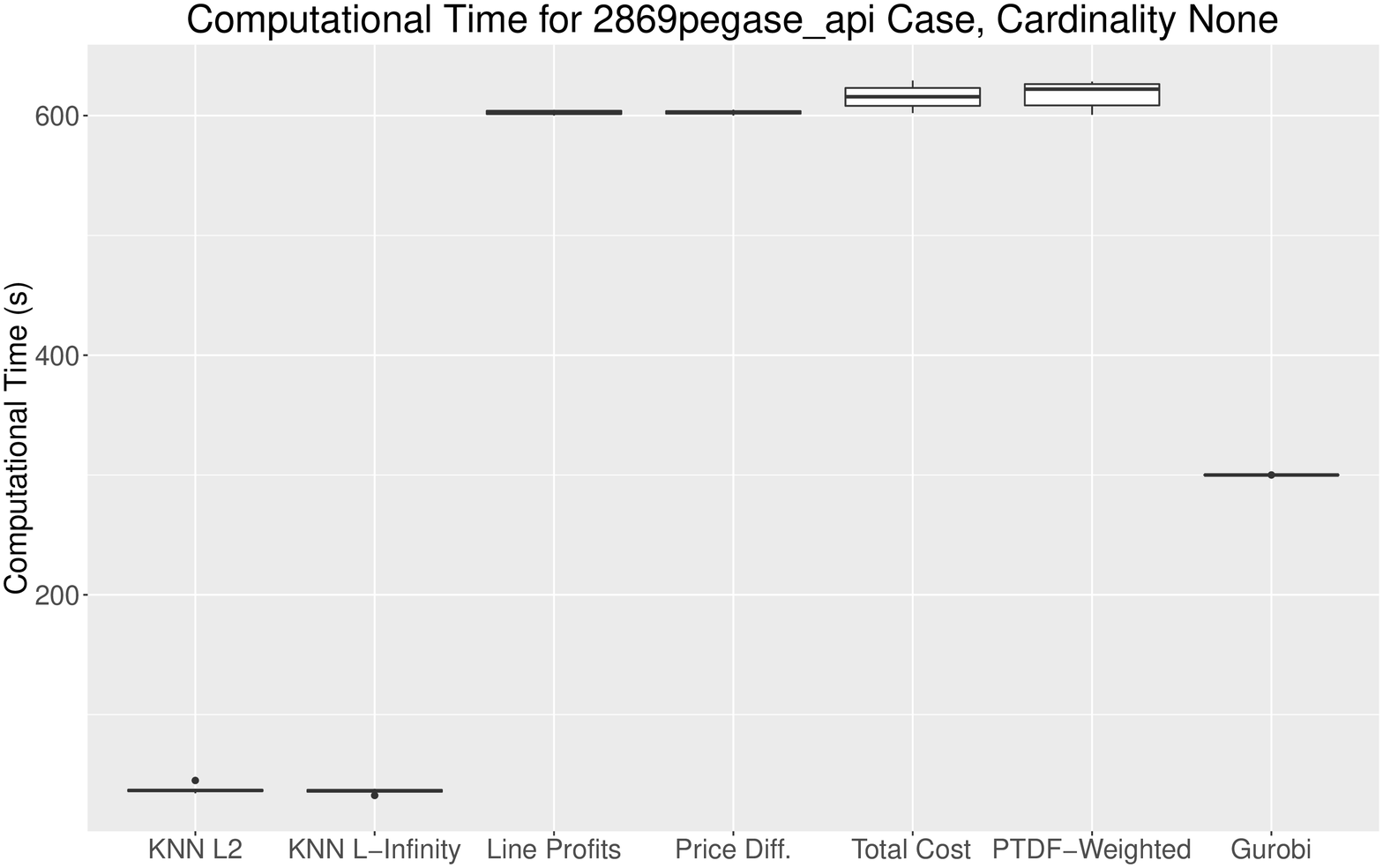}
	\end{tabular}
	\vspace{-0.4cm}
	\caption{Solution quality and computational time results for the 2869\_api bus test case for the three different cardinality options.}
	\label{fig:2869-api-results}
\end{figure*}

\begin{figure*}
	\begin{tabular}{c@{\hskip 0in}c}
		\includegraphics[width=0.5\linewidth]{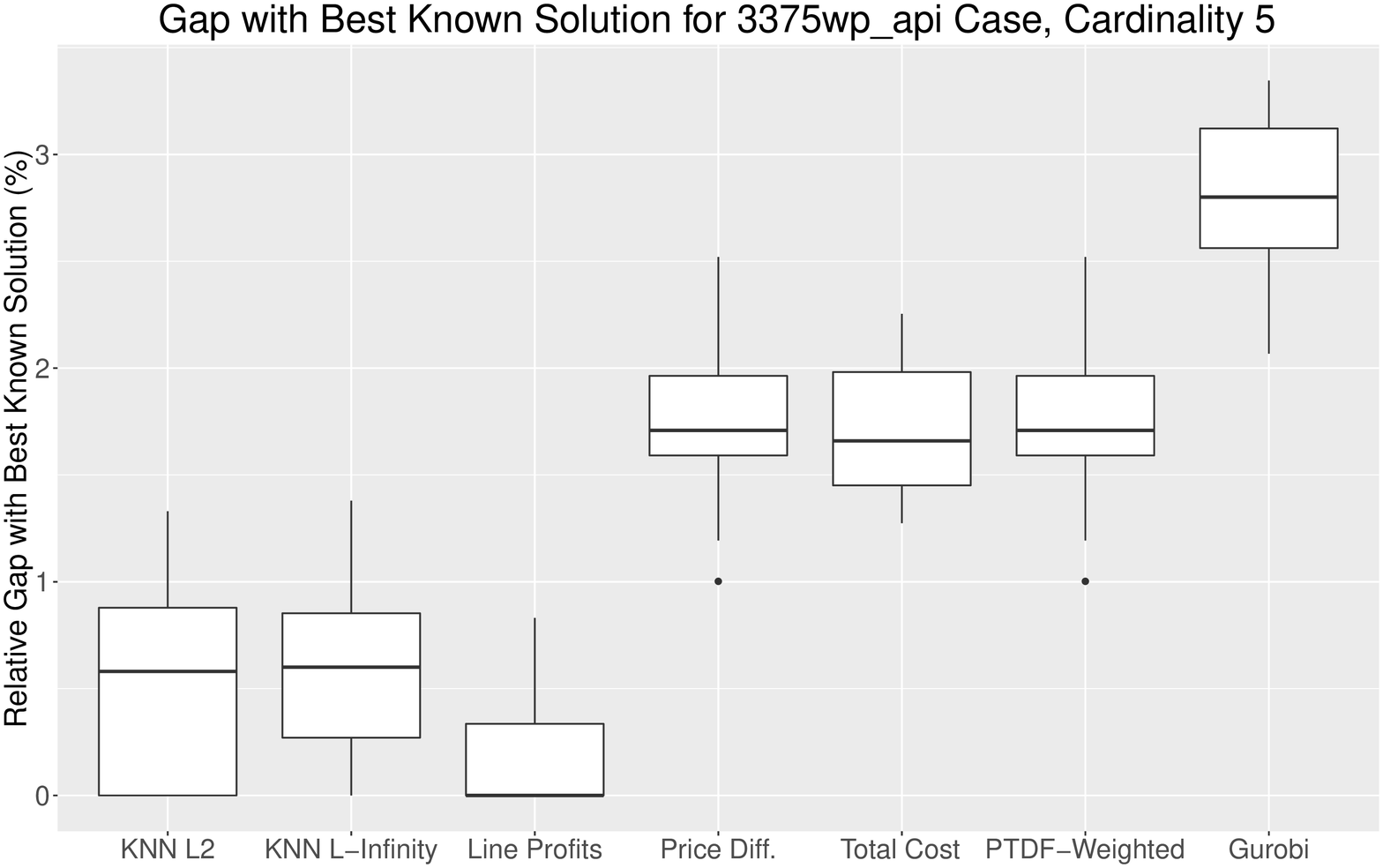} &
		\includegraphics[width=0.5\linewidth]{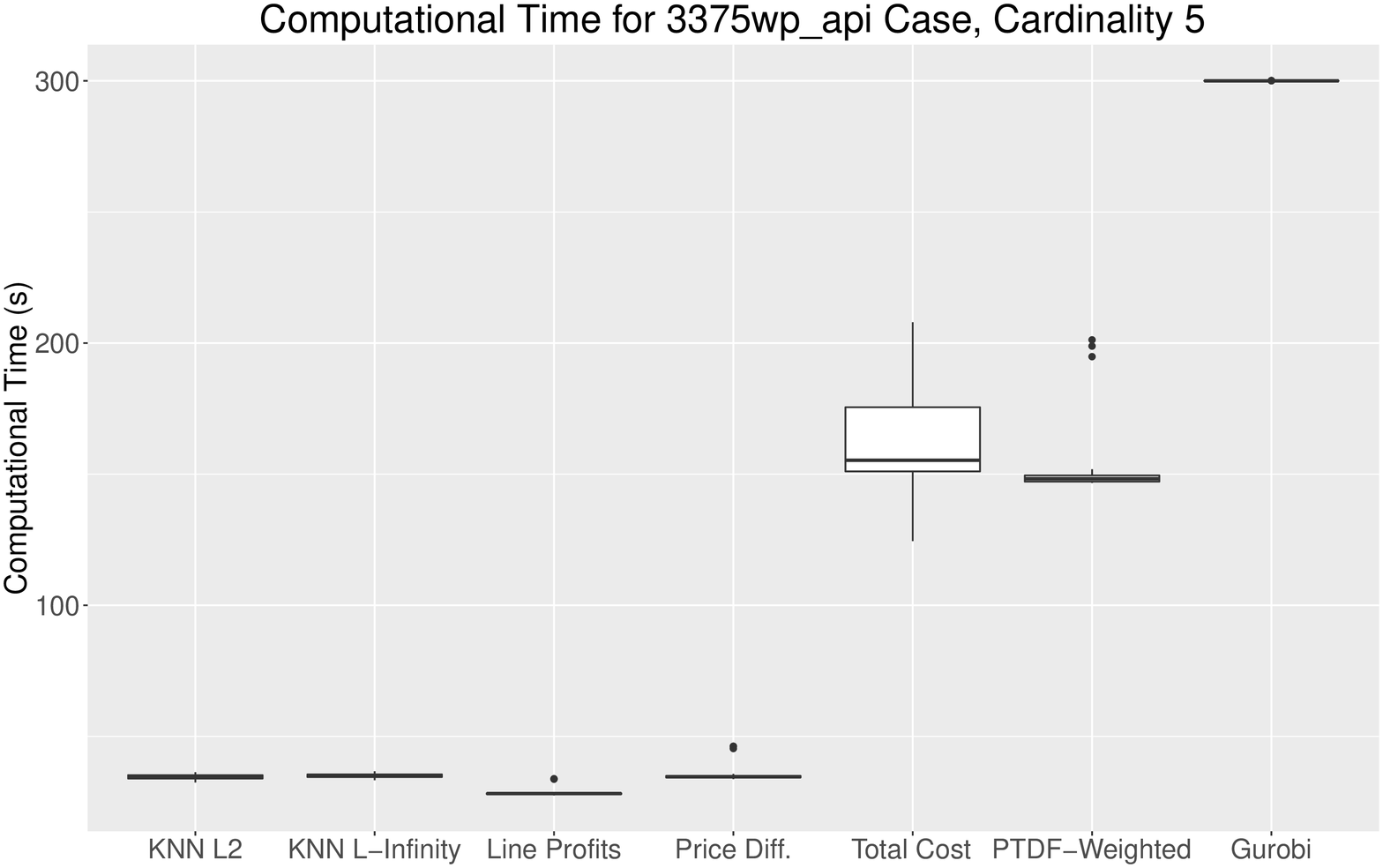} \\
		\includegraphics[width=0.5\linewidth]{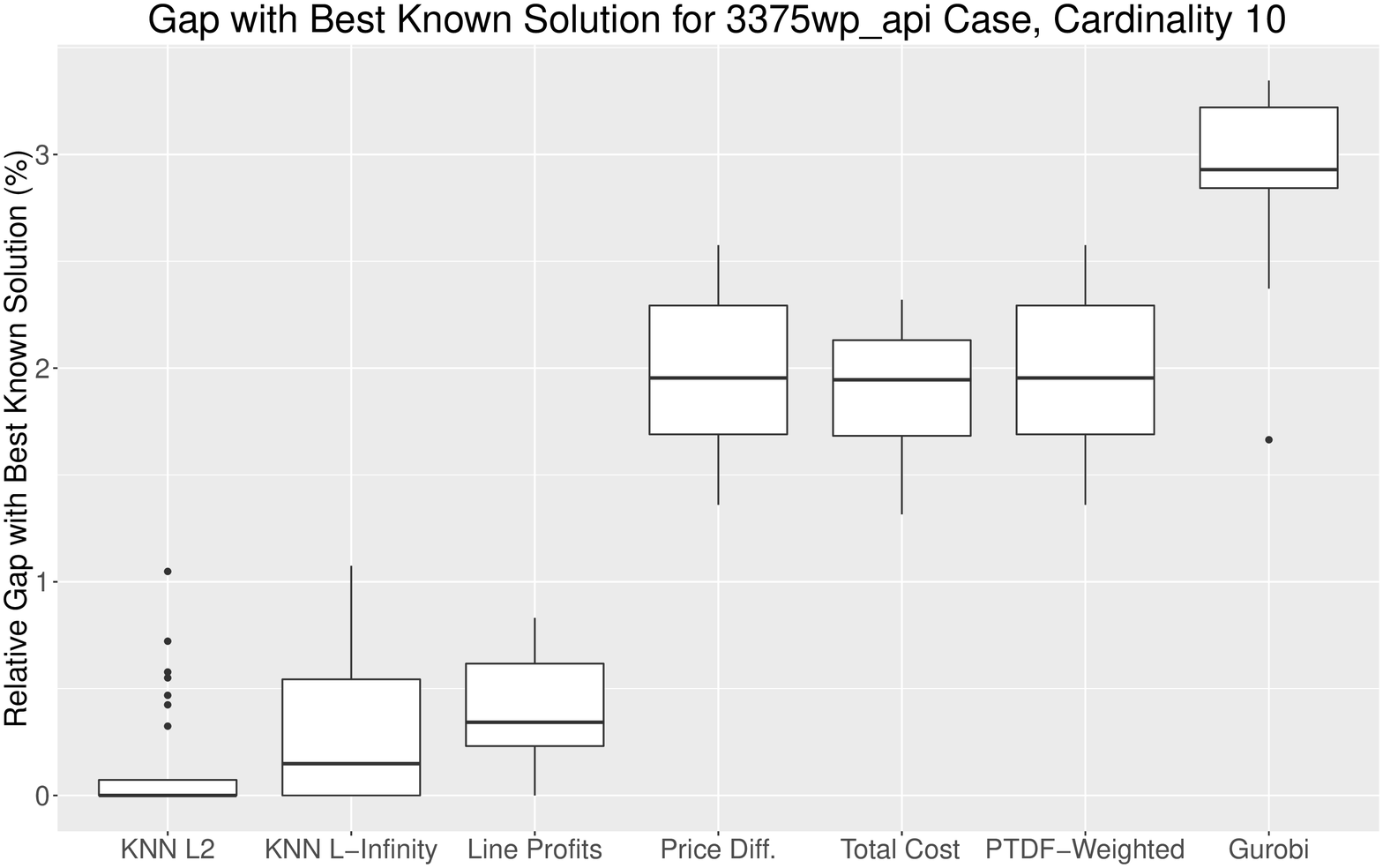} &
		\includegraphics[width=0.5\linewidth]{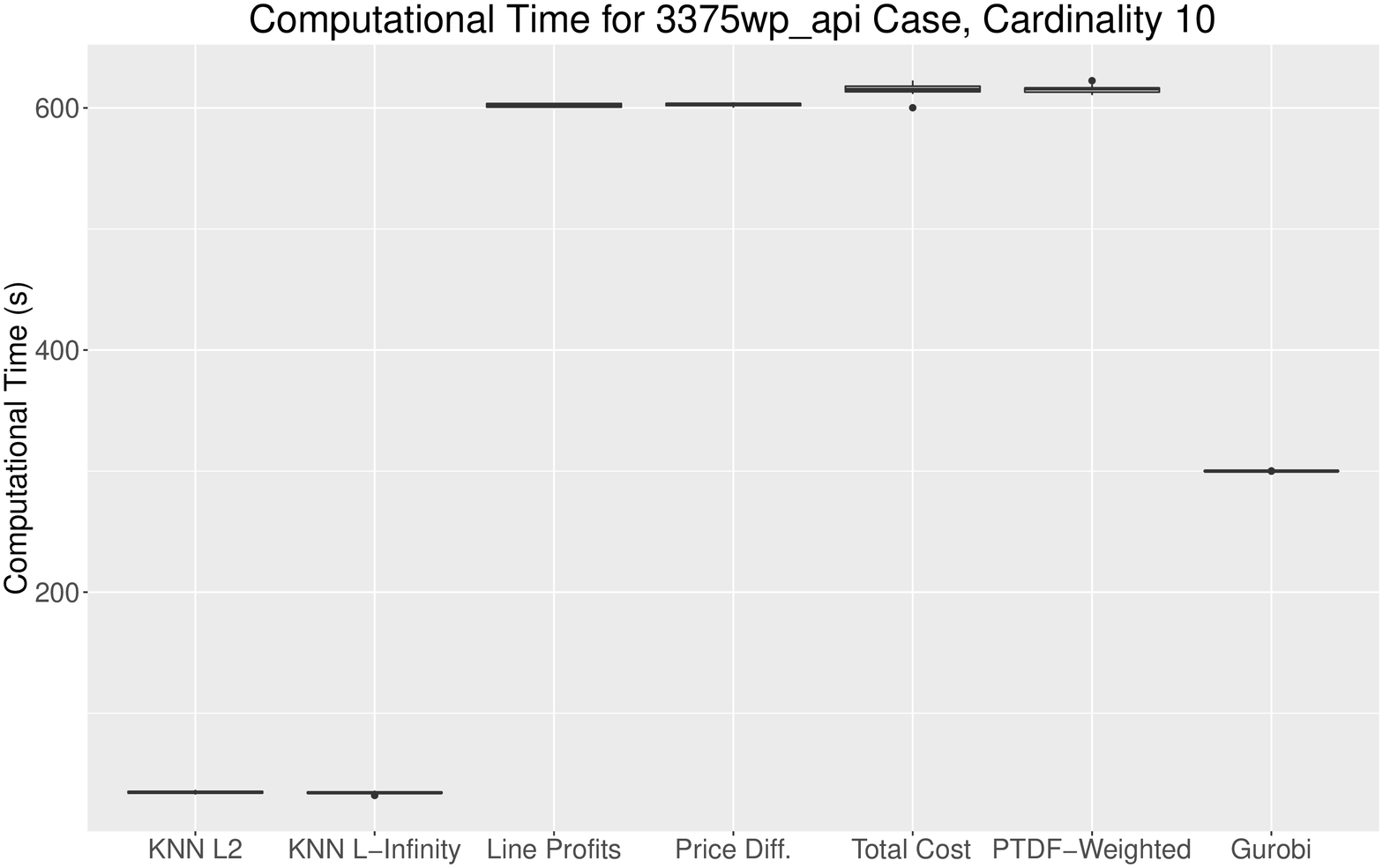} \\
		\includegraphics[width=0.5\linewidth]{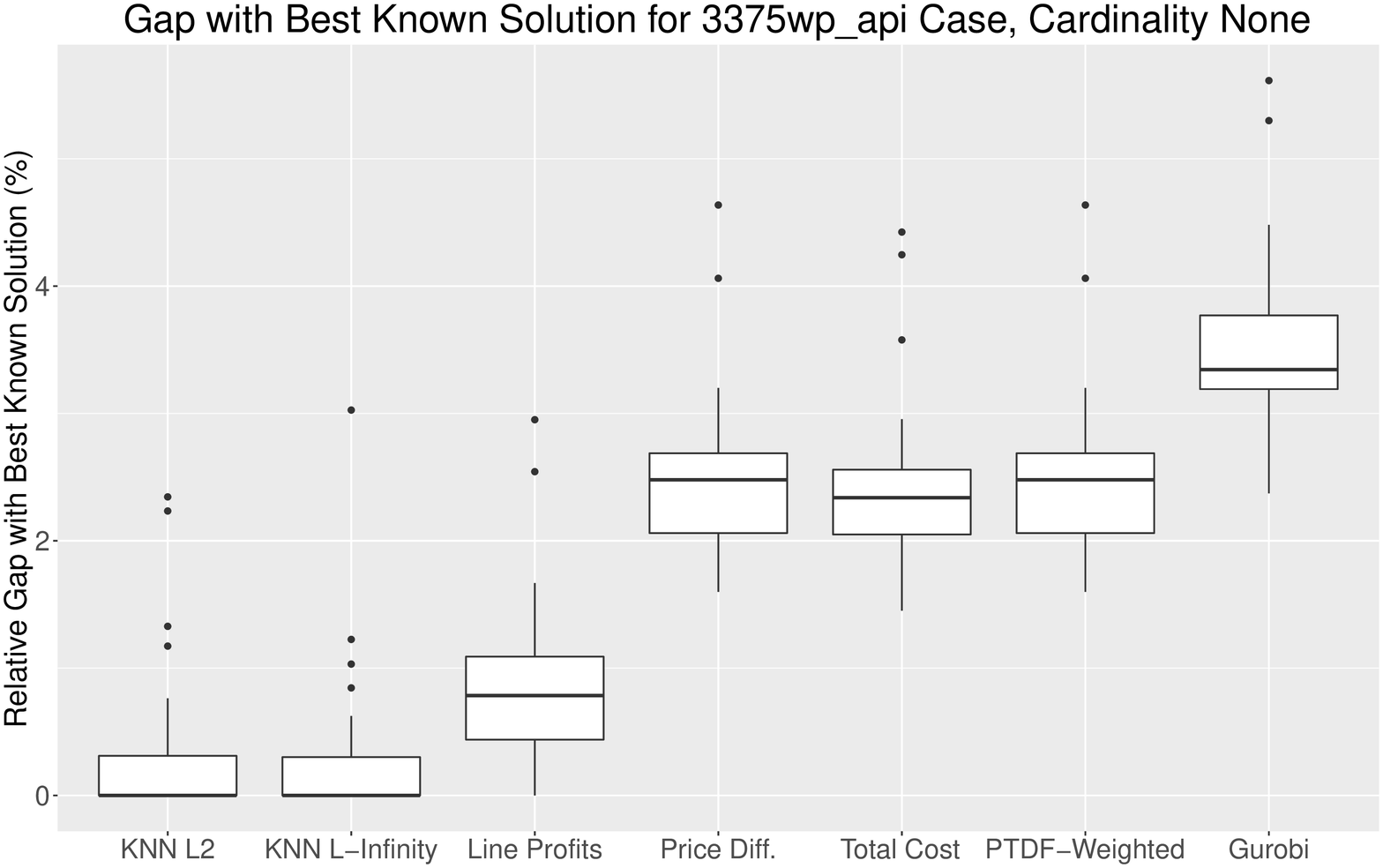} &
		\includegraphics[width=0.5\linewidth]{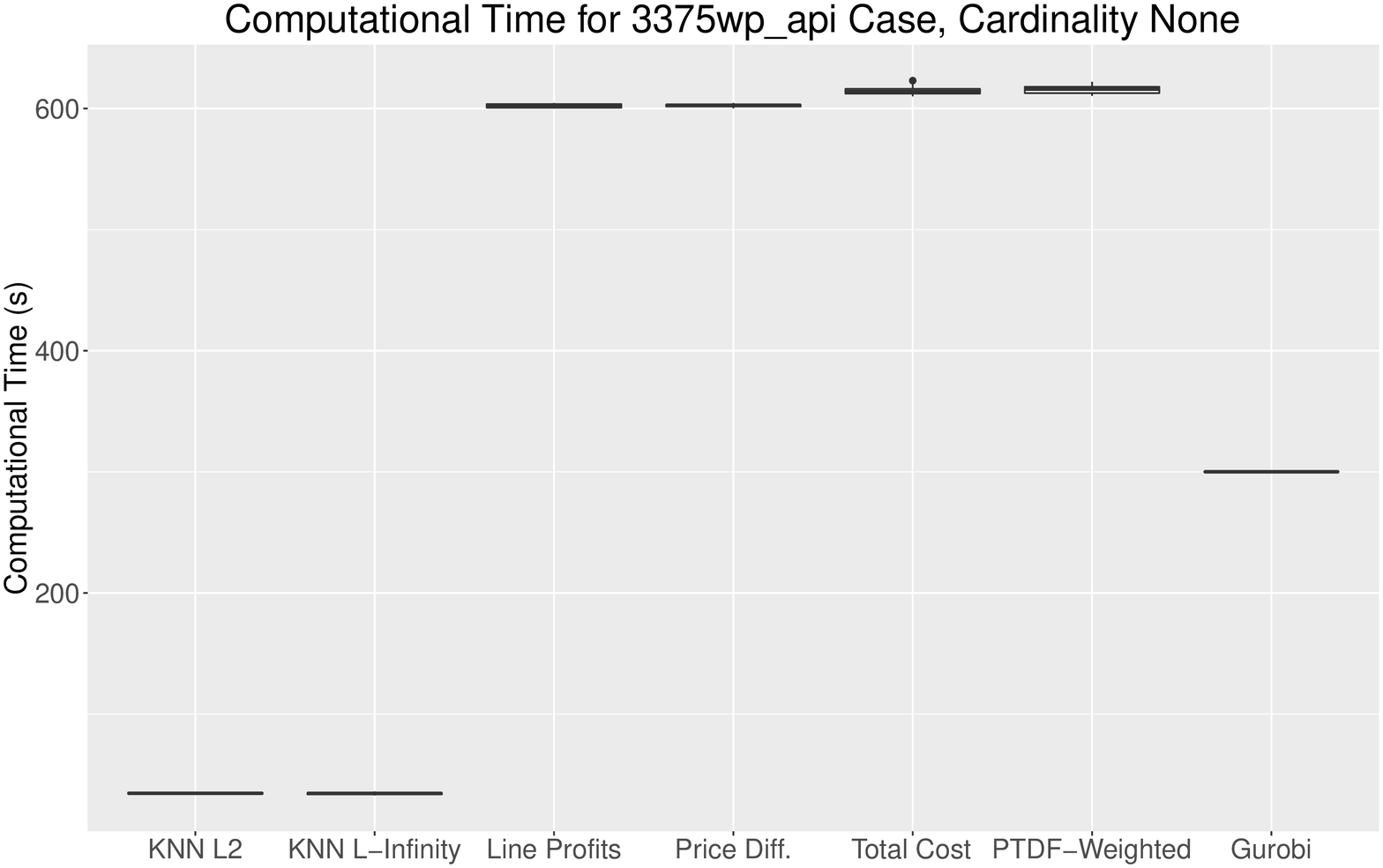}
	\end{tabular}
	\vspace{-0.4cm}
	\caption{Solution quality and computational time results for the 3375 bus test case for the three different cardinality options.}
	\label{fig:3375-results}
\end{figure*}

Our results for the 30 test instances are shown in the boxplots in Figures \ref{fig:118-results} through \ref{fig:3375-results}.
In all plots, the middle line is the median value of the statistic shown, the bottom and top of the box are the 1st and 3rd quartiles respectively, and the ``whiskers" extend to the largest (in absolute value) value from the box which is no further than 1.5 times the interquartile range. Outliers beyond that range are shown as individual points.

The plots in the left-hand column of each figure show the distributions of solution quality for each of the heuristics, where we measure solution quality as the relative gap of the objective from the best-known solution for the test instance.
That is, the gap is given by
$
\frac{z - \hat z}{\hat z},
$
where $z$ is the cost of the heuristic solution and $\hat z$ is the cost of the best known solution for the instance.
Note that, since we only solved the training instances to 1\%-optimality and we stopped at a time limit, it is possible to find a new best-known solution via any of the tested heuristics.
Thus, our best-known solution is not necessarily the solution we calculated during training: It is the best known solution across training and all of the presented heuristics.

The right-hand column of plots shows the distributions of computational times for each of the heuristics.
Note that, for the Total Cost and PTDF-Weighted Cost heuristics, we do not include the time needed to calculate the Load Outage Distribution Factor (LODF) matrix since that only needs to be done once for a network and so can be treated as data.

\subsection{Computational Time}
In all of our test cases and for all the variations in the cardinality constraint, the KNN heuristic has one of the least computational times.
Even on the 3375wp\_api case, it runs for 35 seconds on average, and always takes less than 40 seconds.
This is because the time is completely independent of the cardinality of the set of opened lines and is also relatively agnostic to the size of the network.
Regardless of these, the time to calculate the solution from the KNN heuristic is the time to find the $k$ closest instances from among the 270 training instances, and the time to solve the $k$ DCOPF problems fixing the solutions from these $k$ instances.
Thus, the time only scales up slowly with network size, since network size increases the DCOPF solve times slightly.
In contrast, for the sensitivity-based heuristics and for the greedy local search algorithm, the time scales up with increased cardinality and with the number of lines in the system.
The Gurobi heuristics also do not appear to scale well for increases in the cardinality constraint right-hand side and for increased network size.

In several cases, the Line Profits and Price Difference sensitivity-based heuristics have comparable computational times in the cardinality 5 variant.
However, even in the 118blumsack case, these heuristics are no longer competitive for a cardinality of 10.
Gurobi is never competitive in terms of computational time, often reaching its 5-minute time limit.
The Local Search algorithm was only slightly worse than the KNN heuristic and the sensitivity-based heuristics for case 118blumsack, though its computational time increases with cardinality, like the sensitivity-based heuristics.

\subsection{Solution Quality}
\begin{table}
	\centering
	\caption{Number of test instances for which each heuristic finds the best solution out of all the heuristic solutions. Ties are counted as both of the heuristics finding the best solution, explaining why the sums of the rows can exceed 30. \label{winners}}
	\begin{tabular}{l | r |  r  r  r  r r r  r }
			\hline
			Test Case & Cardinality & KNN L2 & \makecell{KNN\\L-Infinity} & \makecell{Line\\ Profits} & \makecell{Price\\Diff.} & \makecell{Total\\Cost} & \makecell{PTDF-\\Weighted} & Gurobi \\
			\hline
			\multirow{3}{2.2cm}{118blumsack} & 5 & 4 & 3 & 0 & 0 & 0 & 0 & 26 \\
			& 10 & 3 & 3& 0 & 0 & 0 & 0 & 26 \\
			& None & 0 & 0 & 0 & 0 & 0 & 0 & 30 \\
			\hline
			 \multirow{3}{2.2cm}{300kocuk} & 5 & 10 & 10 & 0 & 1 & 0 & 0 & 19 \\
			& 10 & 11 & 9 & 0 & 0 & 0 & 0 & 19 \\
			 & None & 2 & 1 & 0 & 0 & 0 & 0 & 28 \\
			\hline
			\multirow{3}{2.2cm}{1354pegase} & 5 & 22 & 20 & 0 & 0 & 0 & 0 & 0 \\
			& 10 & 22 & 17 & 0 & 0 & 0 & 0 & 3 \\
			 & None & 22 & 20 & 0 & 0 & 0 & 0 & 0 \\
			\hline
			\multirow{3}{2.2cm}{1951rte\_api} & 5 & 12 & 15 & 8 & 0 & 1 & 0 & 0 \\
			& 10 & 17 & 16 & 4 & 0 & 0 & 0 & 0 \\
			& None & 18 & 16 & 4 & 0 & 0 & 0 & 0 \\
			\hline
			 \multirow{3}{2.2cm}{2869pegase} & 5 & 5 & 2 & 10 & 8 & 0 & 7 & 0 \\
			& 10 & 5 & 2& 10 & 8 & 0 & 6 & 1 \\
			& None & 5 & 4 & 9 & 7 & 0 & 7 & 0 \\
			\hline
			\multirow{3}{2.2cm}{2869pegase\_api} & 5 & 6 & 6 & 19 & 0 & 0 & 0 & 0 \\
			& 10 & 8 & 8& 16 & 0 & 0 & 0 & 0 \\
			 & None & 8 & 8 & 16 & 0 & 0 & 0 & 0 \\
			\hline
			\multirow{3}{2.2cm}{3375wp\_api} & 5 & 9 & 4 & 19 & 0 & 0 & 0 & 0 \\
			& 10 & 22 & 12& 5 & 0 & 0 & 0 & 0 \\
			& None & 20 & 17 & 2 & 0 & 0 & 0 & 0 \\
			\hline
			\multirow{3}{2.2cm}{\makecell{Total Over All\\Test Cases}} & 5 & 68 & 60 & 56 & 9 & 1 & 7 & 45 \\
			& 10 & 88 & 67 & 35 & 8 & 0 & 6 & 49 \\
			& None & 75 &66 & 31 & 7 & 0 & 7 & 58 \\
			\hline
			Total &  & 231 & 193 & 122 & 24 & 1 & 20 & 152 \\
			\hline
	\end{tabular}
	{}
\end{table}

In Table \ref{winners}, we show, for each test case, cardinality variant, and heuristic, the number of times that the heuristic finds the best solution out of all the heuristics in our experiments.
(The Local Search algorithm is omitted since it never finds the best solution.)
The second-to-last three rows of Table \ref{winners} show the total number of times each heuristic found the best solution across all the test cases, but still subdivided by cardinality.
In the last row, we show the overall totals, that is, the sum of the previous three rows.
Both the $\ell_2$-norm and $\ell_{\infty}$-norm variants of the KNN heuristic find the best solution for more test cases than any of the other heuristics, and the $\ell_2$-norm is the best by this measure.
When we break these numbers down by cardinality, while the KNN heuristics still find the best solution for more test instances, we see that the margin by which they do so increases as cardinality increases.

In the 118blumsack test case for all cardinalities, Gurobi's heuristics are competitive with the KNN heuristic in terms of solution quality, and usually find a better solution.
However, this is at the cost of computational time: The KNN heuristic runs for at most 3.5 seconds whereas Gurobi runs for at least 24 seconds, but usually for 5 minutes.
For the 300kocuk case, both the Price Difference heuristic and Gurobi are competitive with the KNN heuristic, but this is again at the cost of time, except in the Cardinality 5 version, where the Price Difference heuristic is also very fast.
For the 1354pegase case also, the KNN heuristic yields higher quality solutions than all the comparisons, and it does so in less time in all but the cardinality 5 case.

In the 1951rte\_api, 2869pegase, 2869pegase\_api, and 3375wp\_api cases, one or more of the heuristics from \cite{RuizRCGNP2012} match the KNN heuristic in terms of solution quality, despite the fact that they are terminated early by the time limit for the cardinality 10 case and the case without a cardinality constraint.
In particular the Line Profits and Price Difference heuristics tend to perform well.
However, for all of these cases, only the KNN heuristic achieves this quality in less than 5 minutes for the larger cardinalities, and, in the case of 1951rte\_api, KNN is even slightly faster for cardinality 5.
Last, note that, while some of the sensitivity-based heuristics do sometimes find better solutions than the KNN heuristic, the monetary savings from these are limited.
For example, even in the cardinality 5 version of the 2869pegase case, where perhaps this difference is most extreme, the Line Profits heuristic is on average a 0.06\% improvement over the $\ell_2$-norm version of the KNN heuristic, which translates to an additional \$1,451 in cost savings.
While not insignificant, this shows that the difference in solution quality may be of less importance for less-congested systems where the potential cost savings from transmission switching is limited.
In contrast, for the 1354pegase case with cardinality 5, the average savings of the $\ell_2$-norm version of the KNN heuristic over the Line Profits heuristic is \$3,700.
In summary, how much solution quality matters is partly a function of congestion, since some systems have the potential for much more impact from transmission switching.

In general, we see that the solution quality from the KNN heuristic is fairly constant across different sized networks and different cardinality constraints, usually achieving a relative gap between 0.01\% and 1.0\%.
In contrast, it seems that the performance of the sensitivity-based heuristics is more variable for different power networks. 
With a limit on switching only 5 lines, the Line Profits and Price Difference heuristics also tend to achieve good quality solutions in very little computational time, but in Table \ref{winners}, we see that overall, across all our test systems and cardinality variants, the KNN heuristics find the best solution more often.

\subsection{Feasibility}

As mentioned in Section \ref{sec:learning-heuristic}, it is possible that the KNN heuristic returns a topology which requires load shed or over-generation.
Since using the DC power flow approximation already means that post-processing is required to ensure true feasibility, we allow this to happen.
We found that, for the 118blumsack, 1354pegase, and 2869pegase test cases, our heuristic never returned a solution which had a positive value for load shed or over-generation.
The other four networks all did have infeasible solutions returned, despite our relatively high value of $M$.
However, this is also the case for the other heuristics: For some of our generated instances, it is very difficult to find a feasible solution or the instance is slightly infeasible, and none of the tested heuristics guarantee feasibility.

\begin{figure}
	\includegraphics[width=\linewidth]{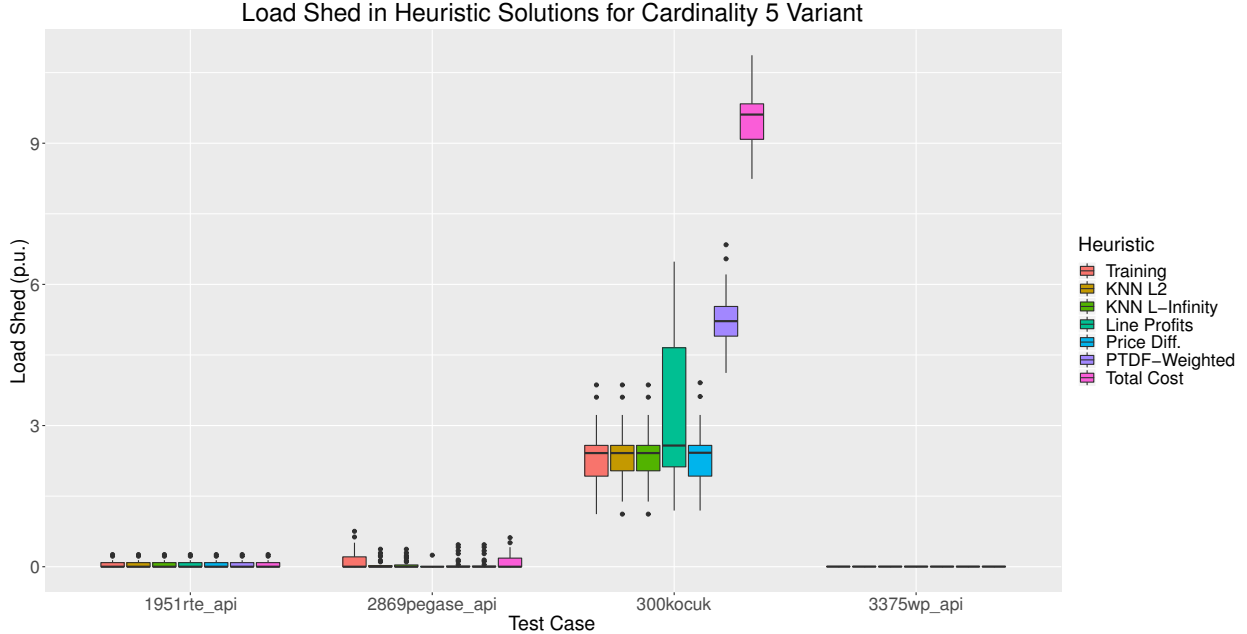}
	\caption{The distribution of the load shed (in per unit) for each of the heuristic solutions to the test instances which are ever infeasible.}
	\label{fig:infeasibility}
\end{figure}

In Figure \ref{fig:infeasibility}, we show the distribution of load shed for the test instances for each of the test cases where infeasibility occurs.
We show this for the cardinality 5 variant, since for instances that require switching to be feasible, this will be the hardest to find feasibility. 
The same plot for the other two cardinality variants looks similar.
We also plot the distributions of load shed in the training solutions, which is about the same as those of the heuristics.
Our instance generating procedure is not guaranteed to give feasible instances, and we did not filter the instances for feasibility.
Some of them are not feasible without transmission switching and some of them are not feasible even with transmission switching allowed.
While this would not be a problem in routine operations, transmission switching has been proposed for contingency response and for disaster response, so it is worthwhile to test our heuristic on cases which require load shedding.

Note that, in all cases where the KNN heuristic returned an infeasible solution, all $k$ topologies we tested were infeasible: This means it is not the value of $M$ which was at fault.
While it may seem that increasing the value of $k$ could mitigate this issue, we are achieving similar levels of infeasibility as the training set, implying that $k$ is also not the reason for this behavior.
In all but the 300kocuk case, the worst load shed for any heuristic solution from the Euclidean norm version of the KNN heuristic is less than 38MW, and on average is less than 5MW.
This is actually a slight improvement over the load shed in the training solutions, which is at worst 76MW and on average no worse than 14MW.
In general the KNN heuristics and the heuristics from \cite{RuizFRC2012} have the same levels of infeasibility, though the Line Profits heuristic has less load shed on average for the 2869pegase\_api case.
The 300kocuk case is not feasible even in the nominal case, explaining the much higher amounts of load shed for all the heuristic solutions.
Lastly, we only saw positive over-generation for two test instances in the 300kocuk case and in two training solutions in the 2869\_api case.
In all of these solutions, the over-generation was less than 2MW.

\subsection{Parameter Tuning}

The difference between the Euclidean and $\ell_\infty$-norms is minimal in most of the experiments.
Unsurprisingly, the two are essentially indistinguishable in terms of computational time.
Overall, the $\ell_\infty$-norm version of the heuristic has a more variable performance and usually is slightly worse on average.
It also tends to have more extreme outliers, occasionally performing much worse than the Euclidean norm on a test instance.
It therefore seems that the Euclidean norm is preferable, as it has both a more consistent and slightly improved performance on the test instances.

Lastly, though the results presented here all have 270 training instances and 10 as the value of $k$, we also experimented with the KNN heuristic's sensitivity to these parameters.
Even for the large networks, 270 training instances appears sufficient: 
There was very little benefit in increasing to 500, though reducing to 100 did hurt the heuristic's performance.
We also experimented with setting $k$ to 3 and 5.
Both of these perform well, though for smaller $k$ we saw small amounts of load shed for the version of the 118blumsack case without a cardinality constraint.
Overall, since the computational time is manageable when $k$ is 10 and the quality and feasibility are improved, the larger value seems preferable.

\subsection{Leave-one-out Cross Validation}

We ran leave-one-out cross validation on the 300 instances for each of the cardinalities with $k = 10$ and using the Euclidean norm.
Results are shown in Table \ref{kfold}.
On average, across all of the test systems and cardinalities, the heuristic returns a result within 0.75\% of the best-known solution. 
In addition, the heuristic is consistent: It appears that some instances of the 2869pegase\_api network could benefit from a larger value of $k$, but with the exception of this network, the worst-case is within 6\% of the best known solution, and is usually even better.

\begin{table}
	\centering
	\caption{Leave-one-out cross validation results: Mean, minimum, and maximum relative gaps of the KNN heuristic solution compared to the best known solution.\label{kfold}}
	\begin{tabular}{l |  r  r  r | r r r | r r r }
			\hline
			\multirow{3}{2cm}{Test Case} & \multicolumn{9}{c}{Cardinality} \\
			&  \multicolumn{3}{c}{5} &  \multicolumn{3}{c}{10} & \multicolumn{3}{c}{None} \\
			 & mean & min & max & mean & min & max & mean & min & max  \\
			\hline
			 118blumsack & 0.10\% & 0.00\% & 1.48\% & 0.24\% & 0.00\% & 3.06\% & 0.71\% & 0.00\% & 12.00\%  \\
			300kocuk & 0.02\% & 0.00\% & 4.49\% & 0.10\% & 0.00\% & 5.24\% & 0.09\% & 0.00\% & 5.23\% \\
			1354pegase & 0.00\% & 0.00\% & 0.07\% & 0.00\% & 0.00\% & 0.07\%& 0.00\% & 0.00\% & 0.07\% \\
			1951rte\_api & 0.06\% & 0.00\% & 3.68\% & 0.06\% & 0.00\% & 3.68\% & 0.06\% & 0.00\% & 3.68\% \\
			2869pegase & 0.01\% & 0.00\% & 0.19\% & 0.01\% & 0.00\%  & 0.21\% & 0.01\% & 0.00\% & 0.21\% \\
			2869pegase\_api & 0.31\% & 0.00\% & 16.9\% & 0.43\% & 0.00\% & 17.4\% & 0.45\% & 0.00\% & 17.4\% \\
			3375wp\_api & 0.03\%& 0.00\% & 1.46\% & 0.03\% & 0.00\% & 1.11\% & 0.08\% & 0.00\% & 3.40\% \\
			\hline
	\end{tabular}
	{}
\end{table}

\section{Reasons the KNN Heuristic Works}\label{sec:why}

While it is difficult to give theoretical justification for the good performance of the KNN heuristic with the small training set, it is possible to prove that KNN can find the correct integer solution for general MIPs with a sufficiently large training set. Formally, we have the following proposition:
\begin{proposition}\label{prop:nieghborhood}
	Suppose we have the following problem:
	\begin{equation}\label{bmip-original}
	\begin{aligned}
	g(b) := \min_{x, y} \quad & c_x^T x + c_y^Ty \\
	\text{s.t.} \quad & Ax + Gy \geq b \\
	& x \in \{0,1\}^n \\
	& y \in \mathbb R^m
	\end{aligned}
	\end{equation}
	Let $(x^*, y^*)$ be an optimal solution which is unique in the integer component and let there be a neighborhood of $b$ such that there exists a feasible solution of (\ref{bmip-original}) everywhere in this neighborhood, with the integer component being $x^*$.
	Then there exists a neighborhood of $b$ such that for all $\hat b$ in this neighborhood, $x^*$ is the integer part of the optimal solution to the corresponding problem.
\end{proposition}

A proof of Proposition \ref{prop:nieghborhood} is given in Appendix \ref{sec:proof-of-prop}.
Note that, according to Proposition \ref{prop:nieghborhood}, if we had a training set such that the union of corresponding neighborhoods covered the space of all possible testing demands, then the KNN heuristic would find the optimal solution. 
In practice, this amount of training data is not realistic, but as we saw in Section \ref{sec:results}, it does not appear necessary. 
We suspect that the success of the KNN heuristic on DCOTS is due to the properties of the problem. 
Our arguments about why the KNN heuristic works are:
\begin{enumerate}
	\item The switching solutions which are optimal within normal variance in demands tend to be composed of similar sets of lines. 
	That is, for a given network, certain lines appear good to open in general. We discuss this further in Section \ref{subsec:few-lines}. \label{few-lines}
	\item There are many near-optimal solutions for a given instance of DCOTS, and optimal solutions for close instances in parameter space tend to be near-optimal for their neighbors. We give results related to this in Section \ref{subsec:distance}. \label{distance-matters}
	\item Changes in demand at certain nodes behave nearly identically to changes in demand at other nodes when viewed in terms of their effect on the optimal topology.
	In essence, changes in demand at individual buses are less important than the overall effect on the congestion patterns in the network. This is discussed further in Section \ref{subsec:dimensionality-reduction}. We find that aggregating certain regions of the network into a net injection leads to no loss in performance for the KNN heuristic.  \label{node-classes}
\end{enumerate}

Points \ref{few-lines} and \ref{node-classes} suggest that, while the KNN heuristic treats the optimal DCOTS solution as a function of the demands and the generation costs, it is really a function of the congestion pattern of the network. 
Somewhat pathologically, we will show in Section \ref{subsec:distance} that, while near solutions in parameter space tend to be high-quality, the best solutions are further.
Thus, while distance is a good proxy for clustering instances with the same DCOTS solution, there could be room for improvement in modeling congestion patterns and their effect on transmission switching decisions.

\subsection{Few Lines are Ever Opened in DCOTS Solutions}\label{subsec:few-lines}

In the following, we will present some results analyzing our training sets for the instances presented in Section \ref{sec:test-cases}.
Though there are many unique DCOTS solutions in our training set, the number of lines ever opened in any training solution is small, and furthermore, some sets of lines tend to be opened together.
This does not mean, however, that every solution in the training set is near-optimal for our test instances.
The fact that the KNN heuristic considers solutions of problems close in the parameter space still appears important to finding high-quality solutions.

\begin{table}
	\footnotesize
	\centering
	\caption{Number of unique topologies and number of lines which are opened in the 300-instance training sets for different test cases. In each entry, the raw number is shown first, and then in parentheses we give that number as percentage of the 300 test instances in the case of unique solutions and as a percentage of the total number of lines in the system in the case of the number of lines opened.\label{table:few-solutions}}
	\begin{tabular}{l |  r  l | r l | r  l | r l | r l | r l }
			\hline
			\multirow{3}{2cm}{Test Case} & \multicolumn{6}{c}{Number of Unique Solutions} & \multicolumn{6}{c}{Number of Lines Opened in Any Solution} \\
			& \multicolumn{6}{c}{Cardinality} & \multicolumn{6}{c}{Cardinality} \\
			 &  \multicolumn{2}{c}{5} &  \multicolumn{2}{c}{10} & \multicolumn{2}{c}{None} &  \multicolumn{2}{c}{5} &  \multicolumn{2}{c}{10} & \multicolumn{2}{c}{None} \\
			\hline
			118blumsack & 55& (18.3\%) & 161& (53.7\%) & 300& (100.0\%) & 27& (14.5\%) & 69& (37.1\%) & 145& (78.0\%) \\
			300kocuk & 6& (2.0\%) & 139& (46.3\%) & 300& (100.0\%) & 7 & (1.7\%) & 65 & (15.8\%) & 277 & (67.4\%) \\
			1354pegase & 247& (82.3\%) & 247& (82.3\%) & 247& (82.3\%) & 181& (9.1\%) & 197& (9.9\%) & 569& (28.6\%) \\
			1951rte\_api & 280& (93.3\%) & 298& (99.3\%) & 300& (100.0\%) & 158 & (3.4\%) & 220 & (4.8\%) & 948 & (20.7\%) \\
			2869pegase & 283& (94.3\%) & 284& (94.7\%) & 292& (97.3\%) & 226& (3.4\%) & 299& (6.5\%) & 1217& (26.6\%) \\
			2869pegase\_api & 294 & (98.0\%) & 299 & (99.7\%)& 300 & (100.0\%) & 240 & (5.2\%) & 352 & (7.7\%) & 2816 & (61.5\%) \\
			3375wp\_api & 276& (92.0\%) & 294 & (98.0\%) & 296 & (98.7\%) & 124 & (3.0\%) & 179 & (4.3\%) & 2689 & (64.6\%) \\
			\hline
	\end{tabular}
	{}
\end{table}
In the left section of Table \ref{table:few-solutions}, we show the number of unique DCOTS solutions in our set of 300 generated instances for each cardinality variant. 
This number is given as a percentage of 300 in parentheses.
Notice that with few exceptions, we tend to have nearly 300 unique DCOTS solutions in this set.
For the lower-cardinality variants of 118blumsack and 300kocuk, this is not true, perhaps because these cases are so congested that there are fewer topologies which make sense, despite variations in demand and generation costs.

In the right section of Table \ref{table:few-solutions}, we show the number of lines which are opened in any of the 300 generated topologies, both as a raw number and as a percentage of the total lines in the system. 
In general, and especially for the lower-cardinality variants, few lines are ever beneficial to open. 
In addition, certain lines are almost always beneficial to open, and certain sets of lines are often opened together.
For example, in the training set for 118blumsack with cardinality 5, one line is switched off in 78\% of the unique solutions.
Another is switched off in 72\% of the solutions. 
However, a third line, which is opened in 28\% of the solutions, is only opened in combination with the first line mentioned in 16\% of the solutions and with the second line mentioned in 12\% of the solutions. 
These patterns suggest that, not only do the same lines get switched in many different solutions, certain sets of lines get switched off together. 
This is good evidence suggesting that heuristics common in the literature which limit the set of switchable lines are finding high-quality solutions if this set is chosen wisely.

\begin{figure}
	\includegraphics[width=\linewidth]{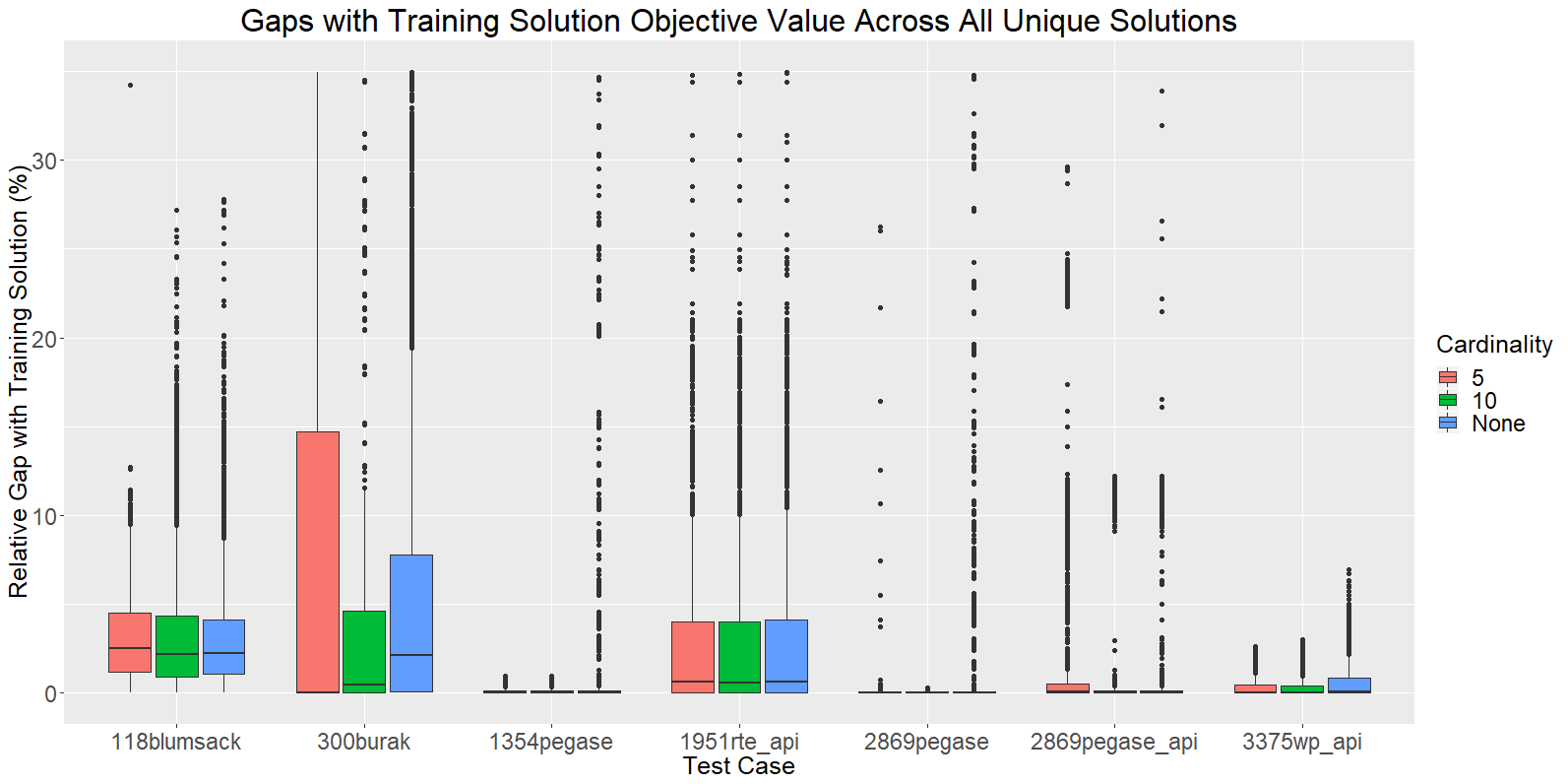}
	\caption{The distribution of the quality of all unique feasible training solutions for all 30 test instances. Note that many outliers (often corresponding to solutions with high load shed or over-generation) are omitted.}
	\label{fig:random-not-good}
\end{figure}

Because certain lines get opened very often in all the training solutions, it seems possible that all of the training solutions are near-optimal. 
This is, however, not the case.
In Figure \ref{fig:random-not-good}, we plot the distribution of relative gaps with the training solution, fixing every unique topology from the training set in each of the test instances (i.e., if there are 50 unique solutions, since there are 30 test instances, this is a distribution of 1,500 gaps).
We can see that, while there is a lot of variation among the test cases, most of the training solutions are at least 5\%-optimal for the test instances.
However, there are many outliers which are far worse, usually because they are actually infeasible (i.e., have positive load shed or over-generation).
On average, these gaps are higher than the gaps of less than 0.75\% which we achieve from the KNN heuristic, so despite the sparsity of the set of lines which are opened, the distance between instances in parameter space still plays a role in finding high-quality primal solutions.

\subsection{Is Distance the Correct Metric?}\label{subsec:distance}

\begin{figure}
	\includegraphics[width=\linewidth]{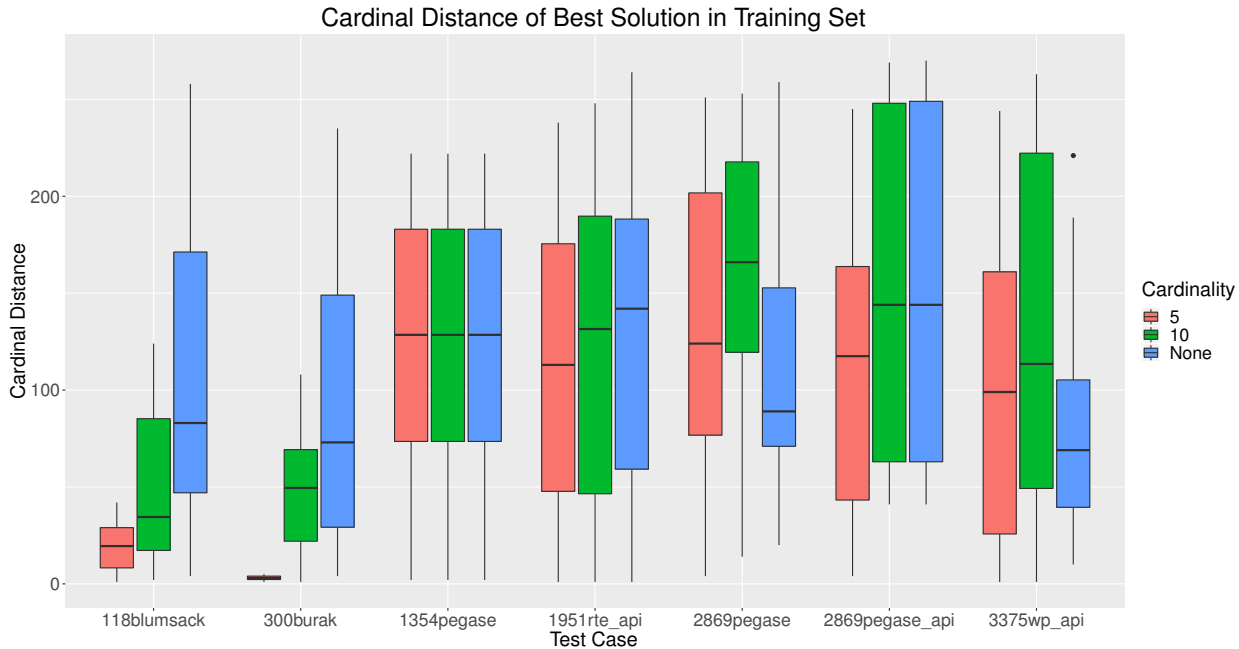}
	\caption{The distribution of the cardinal distances of the best solution in the training set for the 30 test instances in each test case and cardinality variant. The best solution does not tend to be close to the test instance.}
	\label{fig:best-solns-are-far}
\end{figure}

To explore the role of distance in determining instances with the same DCOTS solution, we again considered the set of all unique solutions for a fixed test case and cardinality variation.
This time, using the Euclidean norm, we calculated the cardinal distance of the best training solution from among all the training solutions.
(By cardinal distance, we mean the integer number where 1 is the closest of the instances, 2 is the second closest, and the furthest is the number of unique solutions.)
The distributions of the cardinal distance of the best solution is shown in Figure \ref{fig:best-solns-are-far}.
Surprisingly, the best solutions do not tend to be among the closest, with the exception of the 300kocuk case with cardinality 5.
For most of the instances, the best solution is fairly evenly distributed (keeping in mind that the maximum value for an individual box plot in the figure is the number of unique solutions for that instance).  
This figure suggests that KNN is not the right heuristic for finding the best solution in the training set: There is no reasonable value of $k$ which could consistently find the best solution.

\begin{figure}
	\includegraphics[width=\linewidth]{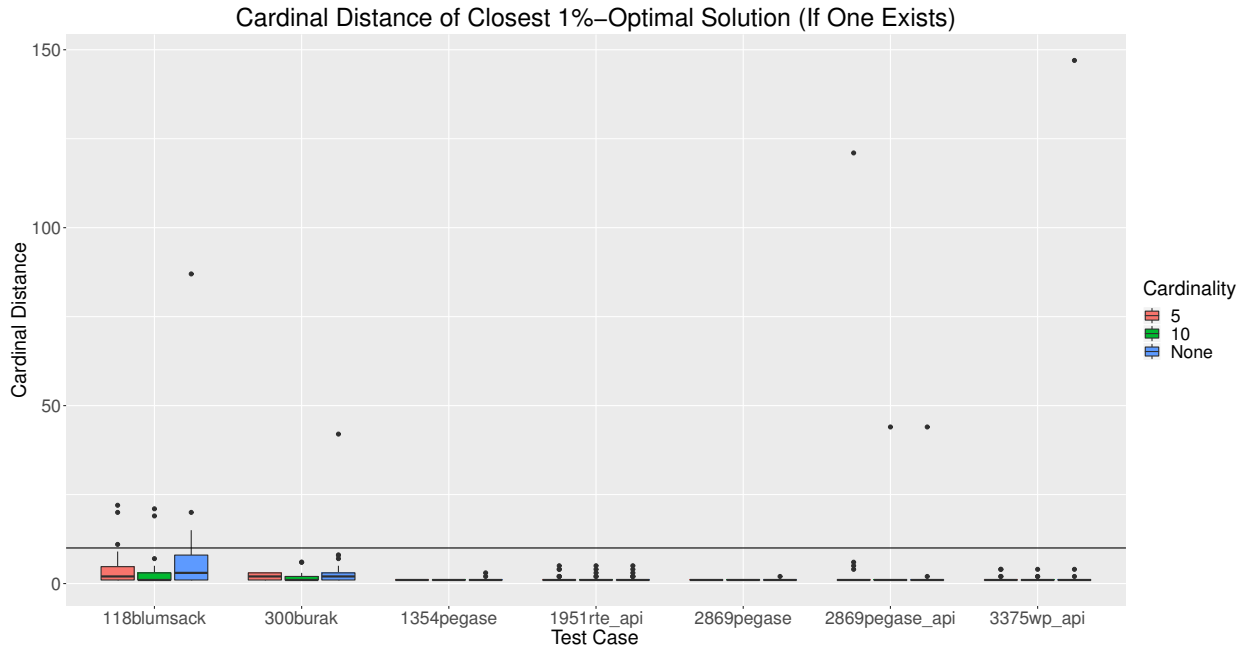}
	\caption{The distribution of the closest 1\%-optimal solution in the training set for the 30 test instances in each test case and cardinality variant. The black horizontal line denotes a cardinal distance of 10.}
	\label{fig:1perc-optimal-are-close-enough}
\end{figure}

However, Figure \ref{fig:best-solns-are-far} obscures the fact that, just because the best solution is far does not mean that near solutions are not good.
In Figure \ref{fig:1perc-optimal-are-close-enough}, we plot the cardinal distance of the closest solution in the training set which is at least 1\%-optimal.
The horizontal line is drawn at the cardinal distance of 10.
This figure illuminates why the KNN heuristic is successful despite what we show in Figure \ref{fig:best-solns-are-far}: Within the neighborhood defined by setting $k=10$, with few exceptions, we have a 1\%-optimal solution.
This stems partly from what we observed in Figure \ref{distance-matters}: In four of our test systems, more than half of the solutions in the training set are within 1\% of optimality for the instances in the test set. 
Note, however, that this is not true for 118blumsack, 300kocuk, or 1951rte\_api, and for these test systems also, $k=10$ is sufficient to find high-quality solutions.


\subsection{Buses with Similar Effect on Congestion}\label{subsec:dimensionality-reduction}

While we have seen that distance does work as a metric to find high-quality solutions within the training set, we also suspect that insight into which changes in demand lead to congestion in the network and where this congestion is would be more powerful in determining a DCOTS solution from the training set. 
Essentially, understanding what directions in the parameter space exacerbate or reduce certain congestion patterns might mean that the distance between a test instance and a training instance is not as important in determining the quality of the training DCOTS solution as the vector between the two instances is. 
If the way the demands have shifted does not change where the network is congested, the same topology may be near-optimal. 
While we do not have a precise characterization of these directions, in this section we show that changes in demand in some parts of the network have identical effects on the DCOTS solution as if those changes had been in other parts of the network.

\begin{figure}
	\includegraphics[width=\linewidth]{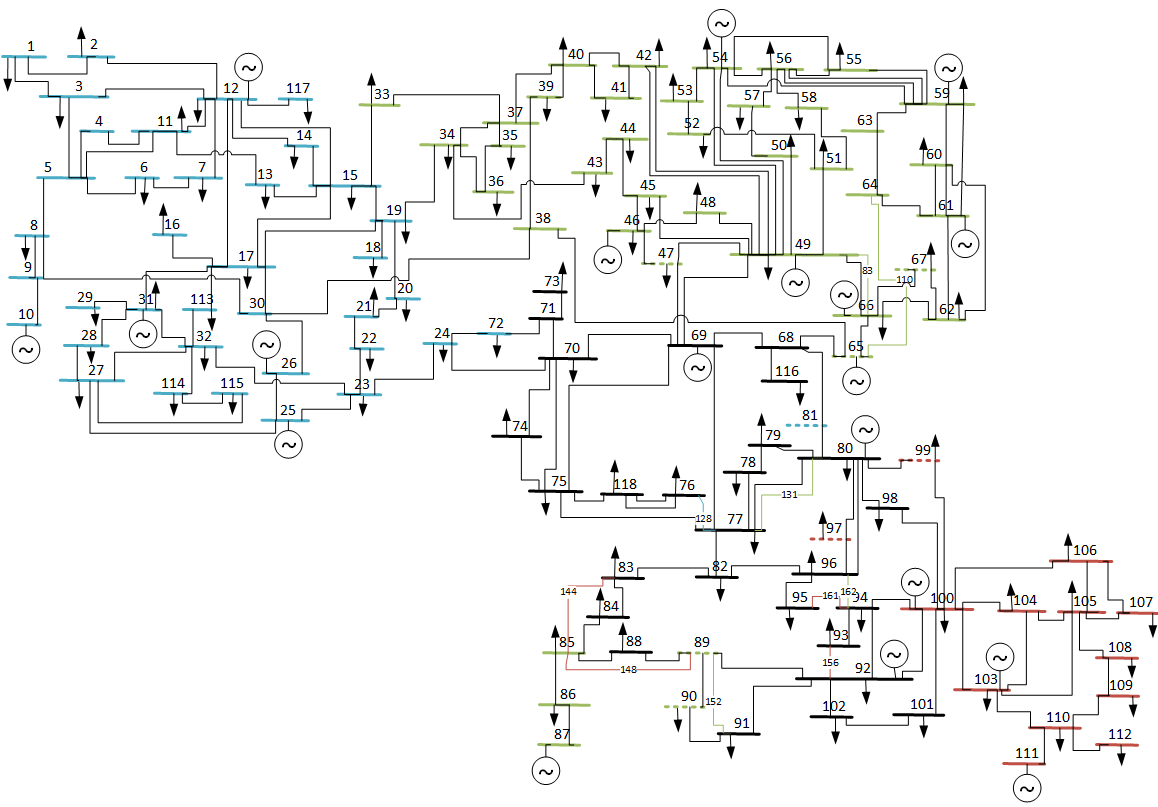}
	\caption{Buses of the 118blumsack network labeled by their optimal switching solution when their demand is increased by 1 p.u. The 35 buses colored green (in the upper-right region) all switch lines 83, 110, 131, 152, and 162 (labeled and also colored green) in their corresponding optimal solution. The 37 buses labeled blue switch lines 110, 128, 131, 152, and 162 in their corresponding optimal solution. The 11 buses colored red switch lines 144, 148, 156, 161, and 162. The buses represented with dotted lines correspond to optimal switching solutions which differ from their color by only one line. Uncolored buses belong to either 1-, 2-, or 3-bus classes not shown in this diagram.}
	\label{img:118blumsack-classes}
\end{figure}

We present results from an experiment where, for each bus in the 118blumsack network, we increased its demand by 1 p.u. above the nominal demands and found a 1\%-optimal transmission switching solution with cardinality 5.
We then sorted the buses of the network into classes that shared the same 1\%-optimal transmission switching solution in the above experiment.
The results are visualized in Figure \ref{img:118blumsack-classes}.
The largest three classes are colored, and buses indicated by dotted lines have a corresponding switching solution which differs by only one line from the class sharing their color.
The lines which are switched in any of these solutions are also labeled and colored according to one of the solutions in which they are opened.

When we sort the buses in this way, they divide regionally.
In this network, the lower right-hand region's demand far exceeds its generation capacity, meaning that power must move downward through buses 77 and 80.
This creates congestion in the part of the network which divides the blue and green classes from the red one.
This diagram leads us to hypothesize that the optimal DCOTS solution is not sensitive to which buses have increased demand within a class, only to whether there is an increase regionally.
In other words, the buses within a class behave like each other in terms of their effect on the overall congestion pattern and hence on the DCOTS solution.

\begin{figure}
	\includegraphics[width=\linewidth]{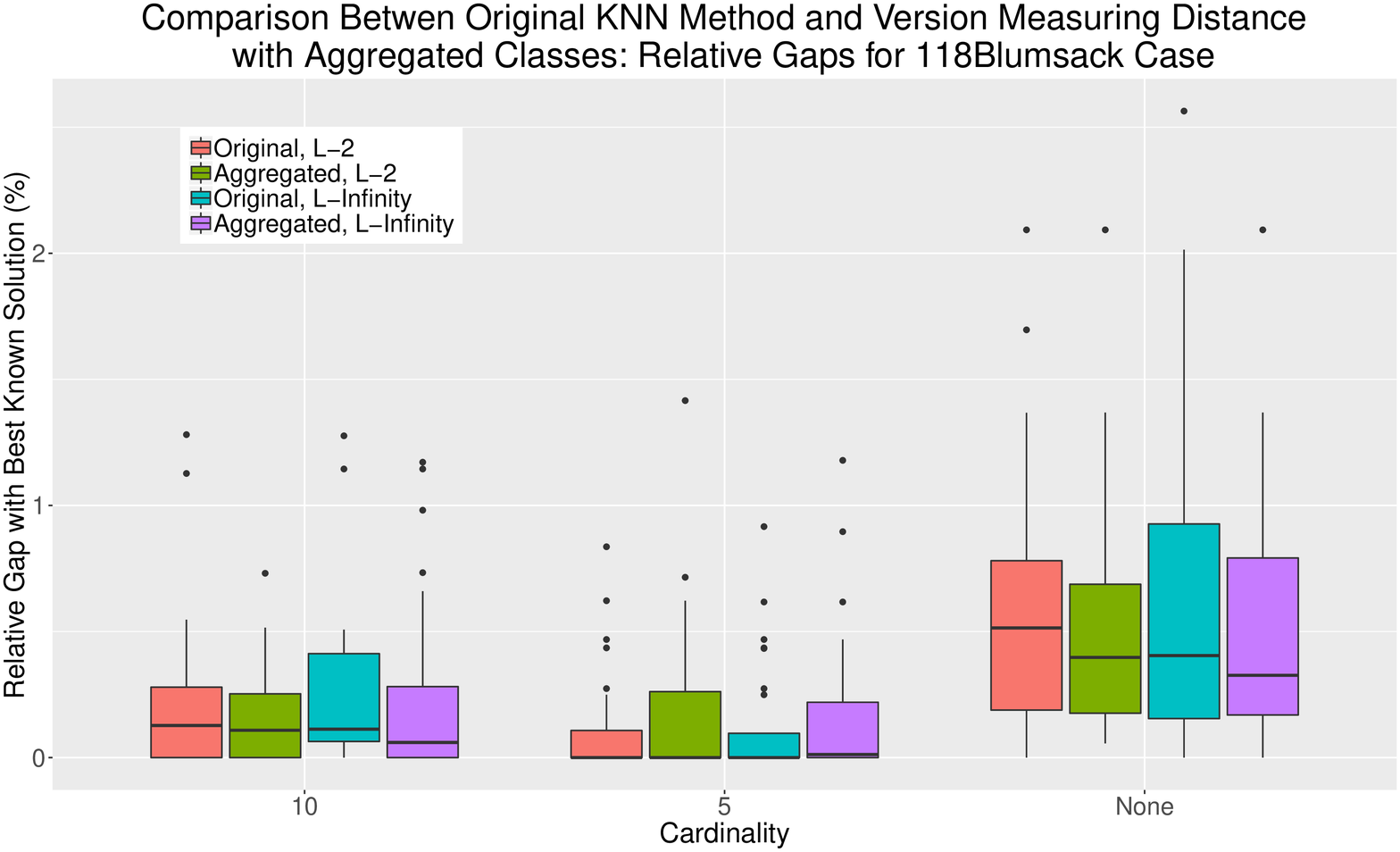}
	\caption{Results comparing solution quality of the KNN heuristic using the original method and using a reduced-dimension vector of net demands in each of the classes from Figure \ref{img:118blumsack-classes}. \label{fig:118-aggregate-class-results}}
\end{figure}

To test this hypothesis, we reran the KNN heuristic on data where only the demand varied and the generation costs were constant.
Instead of characterizing an instance by its vector of demands, we aggregated the net demand within each class (including 2- and 3-bus classes not pictured in Figure \ref{img:118blumsack-classes}).
For the 118blumsack instance, this means that instead of a parameter vector of length 118, we used one of length 30.
The results are compared with the original method (still with only variance in demand, not generation costs) in Figure \ref{fig:118-aggregate-class-results}.
For this instance, the reduced-dimension version of KNN performs on average as well as the original version, and, for the larger cardinalities, has slightly reduced variance.
In experiments not reported here, we tried the same procedure on the 300kocuk case, with similar results.

We do not show results for this experiment for larger test systems because the training to find the classes is much more computationally intensive, as it scales with the number of buses in the network.
We have no reason to believe that this method is an improvement over the original method, given the difficulty of training and the fact that the computational time of the heuristic is not significantly reduced.
(The distance calculations do not dominate the time: The LP solves do.)
However, we present these results because they are revealing of what information matters in determining an optimal topology: Many buses behave like each other when it comes to how their demands affect the DCOTS solution.
In particular, it appears that changes in less-congested regions of the network can be treated in aggregate.
What is important in determining the optimal DCOTS solution is changes in the overall patterns of congestion.

\section{Conclusion}\label{sec:conclusions}

We presented a KNN-based heuristic for DCOTS which finds high-quality solutions within the time limit imposed by real-time.
We showed through a case study on seven test instances that this heuristic yields solutions competitive with heuristics from the literature, and in less computational time.
In particular, the heuristic scales up well with the size of the network since it has only a weak dependence on the number of lines in the system.
Last, through an analysis of our training data, we observe that relatively few lines are ever opened in any DCOTS solutions for a fixed network. 
Despite this, there is variation in the quality of the optimal topologies in the training solutions for the test instances.
While distance is an imperfect metric for finding the best topology, it is proficient at finding high-quality topologies.
However, it would be of interest to be able to characterize the effects of changing demand on congestion, as this relationship appears to drive the optimal DCOTS solution.

Since transmission switching is a tool to reduce generation costs given real-time fluctuation in demand, it is a problem which in practice should be solved quickly, but one for which data from past solves is plentiful.
This makes it an ideal problem for machine learning techniques.
In the future, this same heuristic could be applied to AC optimal transmission switching.
In addition, it could be used to find good warm starts for planning problems or day-ahead operational problems in power systems.

\section*{Acknowledgments}
We would like to thank Cynthia Phillips for helpful discussions during this work.

This work was supported by the Department of Energy, Office of
Electricity Delivery \& Energy Reliability, Advanced Grid Modeling Program
led by Dr. Ali Ghassemian.
Sandia National Laboratories is a multimission laboratory
managed and operated by National Technology and Engineering
Solutions of Sandia, LLC, a wholly owned subsidiary
of Honeywell International, Inc., for the U.S. Department
of Energy’s National Nuclear Security Administration under
contract DE-NA0003525. SAND NO. 2021-5927 O.
The views expressed in the article do not necessarily represent the views of the U.S. Department of Energy or the United States Government.

This work was performed under the auspices of the U.S. Department of Energy by Lawrence Livermore National Laboratory under Contract DE-AC52-07NA27344. This document was prepared as an account of work sponsored by an agency of the United States government. Neither the United States government nor Lawrence Livermore National Security, LLC, nor any of their employees makes any warranty, expressed or implied, or assumes any legal liability or responsibility for the accuracy, completeness, or usefulness of any information, apparatus, product, or process disclosed, or represents that its use would not infringe privately owned rights. Reference herein to any specific commercial product, process, or service by trade name, trademark, manufacturer, or otherwise does not necessarily constitute or imply its endorsement, recommendation, or favoring by the United States government or Lawrence Livermore National Security, LLC. The views and opinions of authors expressed herein do not necessarily state or reflect those of the United States government or Lawrence Livermore National Security, LLC, and shall not be used for advertising or product endorsement purposes.

\bibliography{switching_paper_neutral}

\begin{appendices}
	\section{DCOTS Model}\label{sec:model}
	
	\subsection{Notation}
	
	We will use the following notation to describe the DCOTS model.
	
	\subsubsection{Sets}
	
	\begin{description}
		\item [$\mathcal{L}$] Set of transmission lines
		\item [$\mathcal{G}$] Set of generators
		\item [$\mathcal{B}$] Set of buses
		\item [$\mathcal{G}_b$] Set of generators at bus $b$
		\item [$\mathcal{L}_b^{\text{from}}$] Set of transmission lines leaving bus $b$
		\item [$\mathcal{L}_b^{\text{to}}$] Set of transmission lines entering bus $b$
	\end{description}
	
	\subsubsection{Parameters}
	
	\begin{description}
		\item [$o(l)$] Origin bus of transmission line $l$
		\item [$d(l)$] Destination bus of transmission line $l$
		\item [$B_l$] Susceptance of transmission line $l$
		\item [$\overline{F}_l$] Thermal limit for transmission line $l$
		\item [$\overline{P}_g$] Upper limit of generator $g$ dispatch level
		\item [$\underline{P}_g$] Lower limit of generator $g$ dispatch level
		\item [$D_b$] Demand at bus $b$
		\item [$C_g$] Per-unit generation cost of generator $g$
		\item [$M$] Infeasibility cost
		\item [$K$] Maximum number of lines that can be opened (not to be confused with $k$, the parameter for the KNN heuristic)
	\end{description}

	\subsubsection{Variables}
	
	\begin{description}
		\item [$f_l$] Power flow through transmission line $l$
		\item [$p_g$] Generator dispatch level for generator $g$
		\item [$u_b$] Load shed at bus $b$
		\item [$v_b$] Over-generation at bus $b$
		\item [$\theta_b$] Phase angle for bus $b$
		\item [$y_l$] Indicator of whether or not line $l$ is closed
	\end{description}
	
	We use the following single-period formulation for DCOTS:
	\begin{align}
	\min \quad &\sum_{g \in \mathcal G} C_gp_g + M\sum_{b \in \mathcal B} (u_b + v_b)\label{cost-obj}\\
	\text{s.t.} \quad 
	&f_l \leq B_l(\theta_{o(l)} - \theta_{d(l)}) + 2\pi B_l (1 - y_l) & \forall l \in \mathcal L \label{ohms-lb} \\
	& f_l \geq B_l(\theta_{o(l)} - \theta_{d(l)}) - 2\pi B_l (1 - y_l) & \forall l \in \mathcal L \label{ohms-ub} \\
	&\sum_{g \in \mathcal G_b} p_g + \sum_{l \in \mathcal L_b^{\text{to}}} f_l - \sum_{l \in \mathcal L_b^{\text{from}}} f_l = D_b - u_b + v_b \label{power-balance} & \forall b \in \mathcal B \\
	&-\overline F_ly_l \leq f_l \leq \overline F_ly_l & \forall l \in \mathcal L \label{flow-bounds}\\
	&\theta_{o(l)} - \theta_{d(l)} \geq -\frac{\pi} 6 - 2\pi(1 - y_l) & \forall l \in \mathcal L \label{angle-difference-bounds-lb}\\
	&\theta_{o(l)} - \theta_{d(l)} \leq \frac{\pi} 6 + 2\pi(1 - y_l) & \forall l \in \mathcal L \label{angle-difference-bounds-ub} \\
	&\sum_{l \in \mathcal L} (1 - y_l) \leq K \label{cardinality}\\
	& \underline P_g \leq p_g \leq \overline P_g & \forall g \in \mathcal G \label{generation-bounds}\\
	& -\pi \leq \theta_b \leq \pi & \forall b \in \mathcal B \label{phase-angle-bounds}\\
	&u_b, v_b \geq 0 & \forall b \in \mathcal B \label{slack-bounds}\\
	&y_l \in \{0,1\} & \forall l \in \mathcal L \label{y-binary}
	\end{align}
	This is the model used in \cite{FisherOF2008} with the addition of load shed and over-generation so that the model is feasible for all topologies $y$.
	Generation costs are minimized in the first term of (\ref{cost-obj}), and the second term penalizes infeasibility in the balance constraint.
	We set $M$ to be $10^6$, several orders of magnitude higher than the maximum generation cost.
	Equations (\ref{ohms-lb}) and (\ref{ohms-ub}) are the McCormick relaxation of Ohm's law.
	The nodal balance constraint with slacks included on the right-hand side is given in (\ref{power-balance}).
	Equation (\ref{flow-bounds}) sets the line flow to 0 when the line is opened and enforces transmission limits when the line is closed.
	In (\ref{angle-difference-bounds-lb}) and (\ref{angle-difference-bounds-ub}), we bound the phase angle differences when the line is closed. 
	This is to help the accuracy of the DC approximation.
	We enforce a maximum cardinality of lines which can be opened in constraint (\ref{cardinality}).
	In practice, we expect the value of $K$ to be 5 or 10.
	This is because switching large numbers of lines at once, even in large networks, makes finding an AC feasible dispatch difficult (\cite{CoffrinHLV2014}).
	In addition, most of the cost savings from switching can be attained by switching a small number of lines (\cite{KocukJDLLS2016}).
	Equations (\ref{generation-bounds})-(\ref{y-binary}) enforce variable bounds.
	Note that when the transmission switching binaries $y_l$ are fixed, this formulation is the B-$\theta$ formulation for DCOPF, a linear program. 
	
	\section{Heuristic Algorithms}\label{sec:algorithms-formally}
	
	In this appendix, we give formal descriptions of first the KNN Heuristic algorithm and then the Local Search algorithm.
	
	\begin{algorithm}[H]
		\SetAlgoNoLine
		\DontPrintSemicolon
		\KwIn{Set of solved instances $\mathcal Q = \{I(q^1), I(q^2), \dots, I(q^n)\}$ with $\epsilon$-optimal transmission switching solutions $\{y^1, y^2, \dots, y^n\}$, $p \geq 1$, $k \in \mathbb Z^+$, and an unsolved instance $I(q^{n+1})$}
		\KwOut{Heuristic solution $y^{n+1}$ for instance $I(q^{n+1})$}
		$T \leftarrow \emptyset$ \;
		let $\hat q^{n+1} = \frac{1}{||q^{n+1}||_2}q^{n+1}$ \;
		\For{$I(q^i) \in \mathcal Q$}{\label{dist-loop}
			let $\hat q^i = \frac{1}{||q^i||_2}q^i$ \;
			let $d_i = ||\hat q^{n+1} - \hat q^i||_p$ \;
			\If{$|T| < k$}{
				$T \leftarrow T \cup \{I(q^i)\}$ \;
			}
			\Else{
				\For{$I(q^j) \in T$}{
					\If{$d_i < d_j$}{
						$T \leftarrow (T \setminus \{I(q^j)\}) \cup \{I(q^i)\}$ \;
					}	
				}
			}
		}
		LB $\leftarrow \infty$ \;
		\For{$I(q^{j}) \in T$}{\label{lp-loop}
			Fix the transmission switching solution $y^j$ in $I(q^{n+1})$ and solve. \;
			Let the optimal value be $v_j$ \;
			\If{$v_j <$ LB}{
				LB $\leftarrow v_j$ \;
				$y^{n+1} \leftarrow y^j$ \;
			}
		}
		\caption{KNN for Transmission Switching}
		\label{learning-alg}
	\end{algorithm}
	
	\begin{algorithm}[H]
		\SetAlgoNoLine
		\DontPrintSemicolon
		\KwIn{DCOTS instance $I(q)$, maximum number of lines that can be opened, $K$}
		\KwOut{Heuristic solution $y$ for instance $I(q)$}
		$S \leftarrow \mathcal L$ \;
		$\kappa \leftarrow 0$ \;
		Solve $I(q)$ with all lines closed. Let UB be the optimal value. \;
		\While{$\kappa < K$}{
			$x \leftarrow \infty$ \;
			\For{$l \in S$}{\label{cost-loop}
				In $I(q)$, fix $y_l = 0$ and $y_k = 1$ for all $k \in S \setminus \{l\}$ \;
				Solve $I(q)$ and let $v$ be the optimal value. \;
				\If{$v < x$}{
					$x \leftarrow v$ \;
					$m \leftarrow l$ \;	
				}
			}
			\If{$x \geq UB$}{
				STOP \;
			}
			\Else{
				Fix $y_m = 0$ \;
				UB $\leftarrow x$ \;
				$S \leftarrow S \setminus \{m\}$ \;
				$\kappa \leftarrow \kappa + 1$ \;
			}
		}
		\caption{Greedy Local Search}
		\label{local-search-alg}
	\end{algorithm}
	
	\section{Congestion of the Test Instances}\label{sec:test-instance-congestion}
	
	In this appendix, we give further detail on our seven test instances with regards to their congestion and hence their potential cost benefit from transmission switching.
	
	In Table \ref{all-ones}, we show the average relative gap of the solution where no lines are opened for the 30 test instances on each of our test systems.
	Note that, for the 300kocuk, 1951rte\_api, 2869\_api, and 3375wp\_api cases, many of our training instances are not feasible with all lines closed: They require switching for feasibility, so part of the gap reported in the table is the infeasibility cost in our model.
	Also note that, for many instances, there are no additional savings for higher-cardinality switching budgets.
	This is unsurprising in light of the observation in \cite{KocukJDLLS2016} that few lines are required to be switched in order to get most of the cost benefit from transmission switching.
	The 118blumsack and 3375wp\_api networks are so congested they will switch more than 10 lines if allowed. 
	However, note that, even in these cases, most of the cost savings can still be attained by switching only 5 lines.
	
	\begin{table}
		\centering
		\caption{Average relative gap of the solution with no lines switched open compared to the best known solution.\label{all-ones}}
		\begin{tabular}{l | r  r  r }
				\hline 
				Test Case & Cardinality 5 & Cardinality 10 & No Cardinality Constraint \\
				\hline
				118blumsack*  & 1118.00\% & 1234.00\% & 1248.00\% \\
				300kocuk* & 365.00\% & 371.00\% & 371.00\% \\
				1354pegase & 1.11\% & 1.12\% & 1.11\% \\
				1951rte\_api & 5.86\% & 5.89\% & 5.89\% \\
				2869pegase & 0.35\% & 0.35\% & 0.35\% \\
				2869pegase\_api &  7.63\% & 7.67\% & 7.68\% \\
				3375wp\_api & 2.84\% & 3.03\% & 3.56 \% \\
				\hline
		\end{tabular}
		{Note that for the instances marked with a ``*", the gaps are above 100\% because they are not always feasible with all lines on, so the DCOPF problem is paying an infeasibility cost for load shed. Also note that, in the 1354pegase case, we would expect the relative gap in the no-cardinality case to be the same or better than with cardinality 10. However, since the best known solutions from among all the heuristics are not all optimal, we see a slight error here.}
	\end{table}
	\begin{figure}
		\centering
		\includegraphics[width=0.75\linewidth]{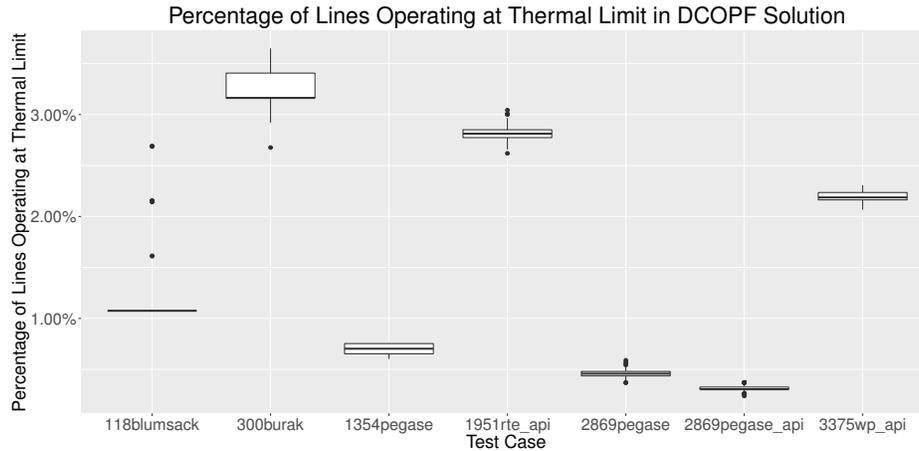}
		\caption{Percentage of lines in each test system which have flow equal to the transmission limits when switching is not allowed for 300 instances on each network each with different demands and generation costs.}
		\label{fig:congestion}
	\end{figure}
	In Figure \ref{fig:congestion} we show the percentage of lines which have flow equal to the transmission limit, which is another indication of the congestion of the system.
	While we expect that more congestion implies more cost savings, note that this is not always the case, since Figure 1 does not take generation costs into account: It is possible to have high numbers of congested lines but still not have as much cost benefit from transmission switching. 
	For example, the 300 bus test case is most extreme by this measure, and, surprisingly, the 2869pegase\_api does not have a high percentage of tight lines, though it does get relatively high cost benefit from switching in Table \ref{all-ones}.
	
	\section{Proof of Proposition \ref{prop:nieghborhood}}\label{sec:proof-of-prop}
	Note that $\{0,1\}^n$ is a finite set, so we can rewrite (\ref{bmip-original}) in the following way:
	\begin{align}\label{bmip}
	g(b) = \min_{x \in \{0, 1\}^n} f_x(b),
	\end{align}
	where
	\begin{equation}\label{fixed-x-lp}
	\begin{aligned}
	f_x(b) = \min_y \quad & c_x^T x + c_y^Ty \\
	\text{s.t.} \quad & Gy \geq b - Ax \\
	& y \in \mathbb R^m
	\end{aligned}
	\end{equation}
	Note that (\ref{fixed-x-lp}) is a linear program, so we know that $f_x(b)$ is piecewise linear and convex in $b$. (We can assume that $f_{\hat x}(b) = + \infty$ if (\ref{fixed-x-lp}) is infeasible when the integer variables are fixed to $\hat {x}$.)
	This means that $g(b)$ is the minimum of finitely many piecewise linear convex functions.
	Fix some $b$ and suppose that $(x^*, y^*)$ is an optimal solution to (\ref{bmip}) which is unique in the integer component.
	
	Pick any other binary vector $\hat{x}$. We will show that there is a neighborhood $N(\hat x)$ of $b$ such that \begin{align}\label{eq:toprove}
	f_{x^*}(\hat b) \leq f_{\hat x}(\hat b) \ \forall \hat{b} \in N(\hat x). 
	\end{align}
	Proving this completes the proof since there are only a finite number of binary vectors, i.e., the final neighborhood we obtain will be the intersection of the neighborhoods corresponding to different binary variables. 
	
	There are two cases to consider:
	\begin{enumerate}
		\item $f_{\hat x}(b) = + \infty$: Since the set of right-hand sides for which an LP is feasible forms a polyhedron, i.e. a closed set, this implies there is a neighborhood of $b$ where $f_{\hat x}(b) = + \infty$. Pick this neighborhood as $N(\hat x)$.
		\item $f_{\hat x}(b) < + \infty$: Note that $f_{\hat x}(b) $ is bounded from below since we have assumed that $f_{x^{*}}(b)$ is the global minimum.  This implies that:
		\begin{enumerate}
			\item There exists a neighborhood around $b$ such that both the functions $f_{\hat x}(\cdot)$ and $f_{x^*}(\cdot)$ are continuous in this neighborhood. (This follows from the fact that these functions are convex and convex functions are continous in the relative interior of their domain.) 
			\item There exists $\epsilon >0$ such that $f_{x^*}(b) \leq f_{\hat x}(b) + \epsilon$.
		\end{enumerate}
		The above two points indicate that there is a neighborhood around $b$ such that (\ref{eq:toprove}) holds. \qed
	\end{enumerate}
\end{appendices}

\end{document}